\documentclass[11pt]{article}
\usepackage{a4}
\usepackage{amsmath}
\usepackage{amsthm}
\usepackage{amssymb}
\usepackage[dvips]{graphicx}
\usepackage{array}

\title{Resurgent Deformations for an Ordinary Differential Equation 
  of Order 2.}

\author{
{\em Eric~Delabaere}
\thanks{D\'epartement de Math\'ematiques, UMR CNRS 6093,
Universit\'e d'Angers, 2 Boulevard Lavoisier, 49045 Angers Cedex 01,
France~/
tel. +33 (0)2.41.73.54.94~/ fax. +33 
(0)2.41.73.54.54~/
e-mail. eric.delabaere@univ-angers.fr}, 
{\em Jean-Marc~Rasoamanana}
\thanks{D\'epartement de Math\'ematiques, UMR CNRS 6093,
Universit\'e d'Angers, 2 Boulevard Lavoisier, 49045 Angers Cedex 01,
France~/
tel. +33 (0)2.41.73.54.55~/ fax. +33 
(0)2.41.73.54.54~/
e-mail. jean-marc.rasoamanana@univ-angers.fr}
}

\newtheorem{thm}{Theorem}[section]
\newtheorem{prop}[thm]{Proposition}
\newtheorem{lem}[thm]{Lemma}
\newtheorem{cor}[thm]{Corollary}
 
\theoremstyle{definition}
 
\newtheorem{defn}[thm]{Definition}
\newtheorem{notation}[thm]{Notation}

\theoremstyle{remark}
 
\newtheorem{rem}[thm]{Remark}

\def\a{{\underline a }}
\def\ap{{\underline a }^\prime}

\def\an{{\underline a }_n}
\def\anp{{\underline a }_n^\prime}
\def\ant{\widetilde{\underline a }_n}
\def\antp{\widetilde{\underline a }_n^\prime}

\begin{document}

\maketitle
 
\begin{abstract}

We consider in the complex field 
the differential equation $\displaystyle 
\frac{d^2}{d  x^2} \Phi(x) = \frac{P_m(x,\a)}{x^2}\Phi(x)$, where $P_m$ is
a monic polynomial function of order $m$ 
with coefficients $\a=(a_1, \cdots , a_m)$. We
investigate the asymptotic, resurgent, properties of the solutions at infinity,
focusing in particular on the analytic dependence on $\a$ 
of the Stokes-Sibuya multipliers. 
Taking into account the non trivial 
monodromy at the origin, we derive a set of functional equations for
the Stokes-Sibuya multipliers. We show how these functional relations
can be used to compute the Stokes multipliers for a class of
polynomials $P_m$. In particular, we obtain conditions for isomonodromic 
deformations when $m=3$.

\end{abstract}

\section{Introduction}

This article is the first of a series of three papers to come. 
The motivation stems  from the well-known theory of Sibuya  \cite{Sib75}
and its Gevrey-resurgent extensions, and their applications 
in spectral analysis.
 
In \cite{Sib75}, Sibuya gives an exhaustive description of the
asymptotic properties when $|x| \rightarrow \infty$ of the solutions of
 the ordinary differential  equation $\displaystyle
 -\frac{d^2 \Phi}{dx^2}+P(x)\Phi=0$, where $P(x) = x^m + a_{1}x^{m-1} +
 \cdots +a_m$ is a complex polynomial function of order $m$. Among
 various results, he shows the existence of a set of fundamental functional
 relations between the Stokes connection matrices, when viewed as functions
 of the coefficients of $P$. The asymptotic behavior
 of the Stokes-Sibuya coefficients when the constant term $a_m$ of $P$ tends
 to infinity is also provided. Later,
these results have been clarified and extended in the framework of the Gevrey
and resurgence theories \cite{R93, L95, Ec85, CNP2, Vo83, DDP93, DP99,
  Costin98, Trinh1}.

One of the main applications of the above results is both the 
qualitative and quantitative description of the spectral set of the 
Schr\"odinger operator $\displaystyle -\frac{d^2}{dx^2}+P(x)$, e.g.,
\cite{DDP97, DP97, KT94, Vo000}. 
Recently, the exact asymptotic analysis has been applied
with success to describe the spectral properties  of a class of $PT$-symmetric
oscillators, i.e.,  when the potential function $P$ satisfies 
$P(x) = \overline{P(-\overline{x})}$. As a rule, such $PT$-symmetric
operators are not Hermitian, therefore the existence of a real
spectrum is a non obvious but interesting question 
(see, e.g.,  \cite{BenBoe99, BenBerry002} for the motivations and
applications in physics).  
 In \cite{DP98, DT00}, the authors consider the
$PT$-symmetric complex cubic oscillator, and prove  the
reality of the spectrum and the appearance of spontaneous symmetry breaking, 
this depending on the values of the coupling constant. Meanwhile, the
reality of the spectrum of the $PT$-symmetric Schr\"odinger operator 
$\displaystyle -\frac{d^2}{dx^2}+P(x)$ with  $P$ polynomial 
was proved by Shin \cite{Shin002} under convenient hypotheses. 
His work relies on a clever examination of the results
of Sibuya, and on ideas and tools usually used  in the context of
integrable models in quantum field theory. \\
As a matter of fact, apart from the key-results of Sibuya, the 
strategy  followed by Shin 
was previously developed by Dorey et al \cite{Dorey} to prove 
the reality of the spectrum of the $PT$-symmetric operator 
$\displaystyle -\frac{d^2}{dx^2}-(ix)^{2m} -
\alpha(ix)^{m-1} +\frac{l(l+1)}{x^2}$.

Our programme is to generalize all the above results to the 
 one dimensional Schr\"odinger operator 
$\displaystyle H= -\frac{d^2}{dx^2}+\frac{P(x)}{x^2}$ with 
$P(x)$  a complex polynomial function of order $m \in
 \mathbb{N}^\star$.

In the  present article, we consider the ordinary differential equation 
$$ \frac{d^2}{d x^2}\Phi(x)=\frac{P_m(x,\a)}{x^2} \Phi(x) 
\leqno{(\mathfrak{E}_m)}$$
in the complex $x$ plane, where $\a := (a_1, \cdots, a_m) \in
\mathbb{C}^m$, 
$m \in  \mathbb{N}^\star$, and~:
$$ 
P_m(x,\a)=x^m+a_1x^{m-1}+ \ldots +a_m \in \mathbb{C}[x] .
$$
This equation admits a unique irregular singular point located at
infinity, and our aim is to concentrate on the asymptotic behaviors of
the solutions of $(\mathfrak{E}_m)$ at this point, and their
deformations in the parameter $\a=(a_1,\ldots,a_m)$. Compared with
Sibuya's work \cite{Sib75}, 
the main novelty comes from the existence  (as a rule)
of another singular point, 
a regular singular one at the origin. The Stokes-Sibuya
coefficients, when considered as functions of the coefficients 
$\a$ of the
polynomial $P_m$, are still governed by a set of independent
functional relations, but the non trivial monodromy at
the origin has now to be taken into account. As we shall see, 
this translates in term of an interesting $\a$-dependent equational resurgence
structure.  \\
The paper is organized as follows. In section \ref{sec2} we study the
solutions near infinity, introducing a well-behaved family of systems
of fundamental solutions.  The associated Stokes-Sibuya coefficients are defined, and
their analytic dependence on $\a$ are analyzed. The main existence
theorem given in \ref{sec2} is proved in section \ref{Resurgentview}
by resurgent methods, and we compare the Stokes-Sibuya coefficients
with 
the Stokes multipliers given by the resurgence
viewpoint. \\ 
In section \ref{sec3}, we introduce
and study a convenient  system of fundamental solutions near the origin, 
by means of the Fuchs theory. By comparing, in section \ref{sec4}, these
different families of fundamental solutions, this yields a set of
functional relations which are detailed in section
\ref{sec5}. Besides describing these general properties, 
we provide different examples in section \ref{sec6} and 
appendix 
\ref{sec7}, 
which will serve as guide lines in 
our next papers. In particular, for $m=3$, we provide a family of
isomonodromic  deformations conditions.

In a second paper, we shall investigate the asymptotics of the
solutions and  Stokes multipliers with respect to the parameter
$\a$, when $\|\a\| \rightarrow \infty$. Roughly speaking, this
corresponds to the second part of Sibuya's work \cite{Sib75}. However, this
work will be done in the framework of the exact WKB analysis, thus taking
advantage of the tools and ideas of the (quantum) resurgence theory,
in the spirit of \cite{Vo83, DP99}.

A third paper will be devoted to applications to  spectral
analysis, especially for $PT$-symmetric operators 
$\displaystyle H= -\frac{d^2}{dx^2}+\frac{P_m(x,\a)}{x^2}$, with a
generalization of the result of
Dorey et al \cite{Dorey} as an interesting by-product.

\section{Solutions of $(\mathfrak{E}_m)$ in the neighbourhood of
  infinity: classical asymptotics }\label{sec2} 

We begin in the framework of the usual Poincar\'e asymptotic analysis,
(see, e.g., \cite{Dieu68, Fed83, Olv74, W65}). We are interested in
the question of the existence of solutions (formal or not) at infinity
for equation $(\mathfrak{E}_m)$.
The starting point of our study will be the main existence 
 theorem \ref{Sibuya}, which can be
seen as
an adaptation of a classical theorem due to Sibuya  (\cite{Sib75},
p. 15).
It asserts the existence and the unicity of a solution of  $(\mathfrak{E}_m)$
defined by its asymptotic expansion at infinity.

\subsection{The main existence theorem}

In the sequel, it will be convenient to think of $x$ as an element
of the universal covering of $\mathbb{C}^\star$ with base point
$1$. Since this covering can be identified to $\mathbb{C}$ provided with the
projection $t \mapsto x = e^t$, we shall associated to $x$ its
argument $\arg (x) \in \mathbb{R}$.

In what follows, 
$$\frac{\sqrt{P_m(x,\a)}}{x} = x^{\frac{m}{2}-1}+\sum_{k=1}^{N}
b_{\frac{m}{2}-k}(\a)x^{\frac{m}{2}-k-1} + O(x^{\frac{m}{2}-N})
$$
stands for the asymptotic expansion  at infinity in $x$ of 
$\displaystyle \frac{\sqrt{P_m(x,\a)}}{x}$.

\begin{thm}\label{Sibuya}

The differential equation  $(\mathfrak{E}_m)$ admits a unique solution
$\Phi_0(x,\a)$  satisfying the following
condition 1.~:
\begin{itemize}

\item 1. $\Phi_0$ is an analytic function in  $x$ in the sector 
$\Sigma_0 = \{|x|>0,\, \mid \arg(x) \mid < \frac{3\pi}{m} \}$ such that,
in  any strict sub-sector of $\Sigma_0$, 
$\Phi_0$ admits an asymptotic expansion at infinity of the following 
form~\footnote{Throughout this theorem,  
    $x^\alpha = \exp\left(\alpha \ln(x)\right)$ with $\ln(x)$ real for
    $\arg(x)=0$.}
$$ T\Phi_0(x,\a)=x^{r(\a)} e^{-S(x,\a)}\phi_0(x,\a),$$
uniformly with respect to $\a$ in
any compact set of $\mathbb{C}^m$, where:
\begin{itemize}
\item i). if $m$ is odd, 
$\displaystyle
\left\{ 
\begin{array}{l}
\displaystyle
S(x,\a)=\frac{2}{m}x^{\frac{m}{2}}+\sum_{k=1}^{\frac{m-1}{2}}
\frac{b_{\frac{m}{2}-k}(\a)}{\frac{m}{2}-k} 
x^{\frac{m}{2}-k} \in \mathbb{C}[\a][x^\frac{1}{2}]\\
\\
\displaystyle r(\a)=\frac{1}{2}-\frac{m}{4}\\
\\
\displaystyle 
\phi_0 \in \mathbb{C}[\a][[x^{-\frac{1}{2}}]] \mbox{  with constant
  term } 1.
\end{array}
\right. 
$
\item ii). if $m$ is even,
$\displaystyle
\left\{ 
\begin{array}{l}
\displaystyle S(x,\a)=\frac{2}{m}x^{\frac{m}{2}}+\sum_{k=1}^{\frac{m}{2}-1}
\frac{b_{\frac{m}{2}-k}(\a)}{\frac{m}{2}-k}
x^{\frac{m}{2}-k} \in \mathbb{C}[\a][x]\\
\\
\displaystyle
r(\a)=\frac{1}{2}-\frac{m}{4}-b_0(\a)\\
\\
\displaystyle 
\phi_0 \in \mathbb{C}[\a][[x^{-1}]] \mbox{  with constant
  term } 1.
\end{array}
\right. 
$
\end{itemize}
Moreover:

\item 2. $\Phi_0$ extends analytically in $x$
to the  universal covering of $\mathbb{C}^\star$, and  
is an entire function  in $\a$.

\item 3. The derivative $\Phi_0'$ of $\Phi_0$ with respect to $x$
  admits an asymptotic expansion at infinity of the form:
$$ T \left( \frac{d}{dx} \Phi_0(x,\a) \right) = 
\frac{d}{dx}\left( T\Phi_0 (x,\a) \right)= 
x^{r(\a)+\frac{m}{2}-1} e^{-S(x,\a)}
(-1+ o(1))$$
when $x$ tends to infinity in  any strict sub-sector of $\Sigma_0$, 
uniformly with respect to $\a$.
 
\end{itemize}
\end{thm} 

Needless to say, the asymptotic expansion $
T\Phi_0(x,\a)$ of $\Phi_0$ at infinity in $\Sigma_0$ can be computed
algorithmically. For instance, for $m=3$ one gets (with Maple)
$$
\begin{array}{c}
T\Phi_0(x,\a) =
e^{-\frac{2}{3}x^{3/2} - a_{1}x^{1/2}}x^{- \frac{1}{4}} \times\\
\\
 \left( 1 + (a_{2} - \frac{1}{4}a_{1}^{2})x^{-1/2}  + 
( - \frac{1}{4}a_{1}^{2}a_{2} + \frac{1}{32}a_{1}^{4} -  
\frac{1}{4}a_{1} + \frac{1}{2}a_{2}^{2})x^{-1} +O(x^{-3/2}) \right),
\end{array}
$$
while for $m=4$:
$$
\begin{array}{c}
T\Phi_0(x,\a) =
e^{( - \frac{1}{2}x^{2} - \frac{1}{2}a_{1}x)}\,x^{( - \frac{1}{2} - 
\frac{1}{2}a_{2} + \frac{1}{8} a_{1}^{2})} \times \\
\\\left(  1 + 
 ( \frac{1}{16}a_{1}^{3} - \frac {1}{4}a_{1}a_{2} - 
\frac{1}{4}a_{1} +  \frac{1}{2}a_{3})x^{-1}  + 
(\frac{5}{32}a_{1}^{2} - \frac{1}{16}a_{2}^{2} - 
 \frac{1}{64}a_{1}^{4}a_{2}  + 
\frac{1}{32}a_{1}^{2}a_{2}^2
 + \frac{5}{32}a_{1}^{2}a_{2} -\right. \\
\\
\left. 
\frac{1}{8}a_{1}a_{2}a_{3} + 
\frac{1}{4}a_{4} -  \frac{1}{4}a_{2} - \frac{9}{256}a_{1}^{4} - 
\frac{1}{4}a_{1}a_{3} -  
\frac{3}{16}  + \frac{1}{512}a_{1}^{6}
+ \frac{1}{32}a_{1}^{3}a_{3} + 
\frac{1}{8}a_{3}^{2})x^{-2} +O(x^{-3}) \right).
\end{array}
$$

We shall discuss the proof of theorem \ref{Sibuya} in a moment 
(\S \ref{Resurgentview}). Here, we would
like to show how one can derive 
fundamental systems of solutions of $(\mathfrak{E}_m)$ 
from  $\Phi_0$ only.  
This is the subject of the following subsection.

\subsection{Stokes-Sibuya coefficients}

In the sequel, it will be convenient to introduce the following
notations:

\begin{notation} 
For all  $\lambda \in \mathbb{C}$ and all  $\a=(a_1, \cdots, a_m) \in
\mathbb{C}^m$, we note
$$\lambda.\a := (\lambda a_1, \cdots , \lambda^m
a_m).$$
We set:
$$
\omega=e^{\pm \frac{2i\pi}{m}}
$$
and we introduce:
$$
\forall k \in \mathbb{Z}, \, \, \, \, 
\Phi_k(x,\a)=\Phi_0(\omega^k x , \omega^k . \a).
$$
where $\Phi_0$ is given by theorem \ref{Sibuya}.
\end{notation}

We bring into play a quasi-homogeneity property of equation
$(\mathfrak{E}_m)$. 
We note that equation $(\mathfrak{E}_m)$ is invariant under the
transformation $(x,a) \mapsto (\omega x , \omega . \a)$ so that, with
the above notations,
theorem \ref{Sibuya} easily translates into the following lemma~: 

\begin{lem}\label{lemme1} 
For any
$k \in \mathbb{Z}$, $\Phi_k$ is a solution of $(\mathfrak{E}_m)$, and is 
entire  in $a$. Its 
asymptotic expansion when $x$ tends to infinity in the sector
$\Sigma_k =\{|x|>0, \, \mid arg(x)+k. arg(\omega)
\mid < \frac{3\pi}{m} \}$, uniformly in $\a$ in
any compact set of $\mathbb{C}^m$, is given by~:
$$T\Phi_k(x, \a)=T\Phi_0(\omega^k x , \omega^k .\a)$$
 where $T\Phi_0$ is the  asymptotic expansion of $\Phi_0$ in $\Sigma_0$
described  in theorem \ref{Sibuya}.
\end{lem}

We deduce the following corollary~:

\begin{cor} For every
$k \in \mathbb{Z}$, the solution $\Phi_k$ is exponentially
decreasing (``subdominant function" in \cite{Sib75}, p. 19) in the sector 
$\Lambda_k=\{ \mid arg(x)+k .arg(\omega)
\mid < \frac{\pi}{m} \}$.
\end{cor}

We note that the sectors $\Lambda_{k-1}$, $\Lambda_k$ and 
$\Lambda_{k+1}$ are included in $\Sigma_k$
and, by the previous lemma \ref{lemme1}, each solution $\Phi_k$ has an
exponential growth of order not greater than $\displaystyle
\frac{m}{2}$
 in the two sectors
$\Lambda_{k-1}$ and $\Lambda_{k+1}$ adjacent to $\Lambda_k$. This
allows to show the following lemma:

\begin{lem}\label{lemme2}
For every $k \in \mathbb{Z}$, $\{ \Phi_k,\Phi_{k+1} \}$ constitutes a
fundamental system of solutions of $(\mathfrak{E}_m)$ and moreover,
$$W(\Phi_k,\Phi_{k+1}) =  2(-1)^k
\omega^{k(1-\frac{m}{2})+r(\omega^{k+1}.a)}$$ 
where $W$ denotes the Wronskian,
while $r$ is given by theorem \ref{Sibuya} .
\end{lem}

\begin{proof} By quasi-homogeneity of
$P_m(x,\a)$, we note that
$$S(\omega x, \omega.\a)=-S(x,\a).$$
Thus, by lemma \ref{lemme1}, for $x \in \Sigma_k$,
$$\displaystyle T\Phi_k(x,\a)=\omega^{kr(\omega^k.\a)} 
x^{r(\omega^k.\a)}  e^{(-1)^{k-1}S(x,\a)}
(1+o(1)).$$
Using part 3. of Theorem \ref{Sibuya} we have also, for $x \in \Sigma_k$,
$$\displaystyle T\Phi'_k(x,\a)=
(-1)^{k-1} \omega^{kr(\omega^k.\a)} 
x^{r(\omega^k.\a)+\frac{m}{2}-1} 
e^{(-1)^{k-1}S(x,\a)}(1+o(1)).$$
Moreover, the coefficient $b_0$ of theorem 1 is a quasi-homogeneous polynomial
in $\a$ of order $\frac{m}{2}$ so that~:
$$r(\omega^k.\a)+r(\omega^{k+1}.\a) = 1-\frac{m}{2}.$$
As a result, we get the equalities~: for $x \in \Sigma_k \cap
\Sigma_{k+1}$,
$$
\displaystyle W(\Phi_k,\Phi_{k+1})  = 
 \displaystyle \Phi_k\Phi'_{k+1}-\Phi'_k \Phi_{k+1}$$
$$ = 
 \displaystyle 2(-1)^k \omega^{k(r(\omega^k. \a)+
r(\omega^{k+1}.a))+r(\omega^{k+1}.\a)} 
x^{r(\omega^k.\a)+r(\omega^{k+1}.\a)+\frac{m}{2}-1}
(1+o(1))$$
$$= 
 \displaystyle 2(-1)^k \omega^{k(1-\frac{m}{2})+r(\omega^{k+1}.\a)}
(1+o(1)).
$$
The Wronskian  $W(\Phi_k,\Phi_{k+1})$ being independent of $x$, this
completes the proof.
\end{proof}

Since each system $\{ \Phi_k,\Phi_{k+1} \}$ constitutes a
fundamental system of solutions of $(\mathfrak{E}_m)$ we deduce, from 
the classical theory on linear differential
equations, the existence of functions $C_k(\a)$,
$\widetilde{C}_k(\a)$   depending
only on the variable $\a$,  such that~:
\begin{equation}\label{uneequation}
\forall k \in \mathbb{Z}, \Phi_{k-1} = C_k(\a) \Phi_k +
\widetilde{C}_k(\a) \Phi_{k+1}.
\end{equation}

\begin{defn}
The functions $C_k(\a)$ and $\widetilde{C}_k(\a)$  defined by 
(\ref{uneequation})
are called the {\sl Stokes-Sibuya coefficients} of $\Phi_{k-1}$
associated respectively with $\Phi_k$ and $\Phi_{k+1}$. The matrices
$\displaystyle \mathfrak{S}_k (\a) :=
\left(
\begin{array}{cc} 
C_k(\a) & \widetilde{C}_k (\a) \\ 
1 & 0 
\end{array} 
\right)$ are called the {\sl Stokes-Sibuya connection matrices}.
\end{defn}

Differentiating the previous equalities (\ref{uneequation}) 
with respect to the variable $x$,
we obtain the following formulas~:
\begin{equation}\label{cramer1}
 C_k(\a) = \frac{W(\Phi_{k-1},\Phi_{k+1})}{W(\Phi_k,\Phi_{k+1})} 
\end{equation}
and
\begin{equation}\label{cramer2}
 \widetilde{C}_k(\a) = \frac{W(\Phi_{k-1},\Phi_k)}{W(\Phi_{k+1},\Phi_k)}.
\end{equation}

Using the fact that the $\Phi_k$'s are entire functions in $\a$, we
deduce from  (\ref{cramer1}), (\ref{cramer2}) and lemma \ref{lemme2}
that
the Stokes-Sibuya coefficients are entire functions in $\a$. Also,
it follows from the very definition of the $\Phi_k$'s,  from  
(\ref{cramer2}) and lemma \ref{lemme2} that
$\displaystyle C_k(\a) = C_0(\omega^k .\a)$, while  
 $\displaystyle \widetilde{C}_k(\a) = \widetilde{C}_0(\omega^k .\a) = 
\omega^{m-2+2r(\omega^k.\a)}$. In particular, 
since $\omega^m = e^{\pm 2i\pi}$, we get: for all $k \in \mathbb{Z}$, 
$C_k = C_{k \!\! \mod m}$ and $\widetilde{C}_k = \widetilde{C}_{k
  \!\! \mod m}$.

We summarize our results.

\begin{thm}\label{SibuyaCoeff}
For all $k \in \mathbb{Z}$ we note
\begin{equation}\label{lesPhik}
\Phi_k(x,\a)=\Phi_0(\omega^k x, \omega^k.\a),
\end{equation}
where $\Phi_0$ is the solution of $(\mathfrak{E}_m)$
 defined in theorem \ref{Sibuya}. 
Then, for every $k\in \mathbb{Z}$, 
\begin{itemize}
\item 
$\Phi_k(x,\a)$ is analytic in $x$ on
the universal covering of $\mathbb{C}^\star$ and entire in $\a$.
\item  The system $\{ \Phi_k,\Phi_{k+1} \}$ constitutes a
fundamental system of solutions of $(\mathfrak{E}_m)$.  
\item  We have
\begin{equation}\label{lesCk}
\left(\begin{array}{c} 
\Phi_{k-1} \\ 
\Phi_k \end{array} 
\right) (x, \a) =
\mathfrak{S}_k (\a)
\left(\begin{array}{c} 
\Phi_k \\
\Phi_{k+1} 
\end{array} 
\right)(x, \a),
\end{equation}
where the Stokes-Sibuya connection matrix $\mathfrak{S}_k (\a)$ is 
invertible, and entire in $\a$. Moreover, 
\begin{equation}
\mathfrak{S}_k (\a) = \mathfrak{S}_{k-1} (\omega .\a), \hspace{5mm}
\mathfrak{S}_k (\a) = \mathfrak{S}_0 (\omega^k .\a).
\end{equation}
\item 
The Stokes-Sibuya coefficients  $ C_k(\a)$ and $\widetilde{C}_k(\a)$
associated respectively with $\Phi_k$ and $\Phi_{k+1}$  are entire
functions in $\a$ and, 
\begin{equation}\label{qhCk}
\left\{
\begin{array}{l}
\displaystyle C_k(\a) = C_0( \omega^k.\a),
\hspace{5mm} C_k = C_{k \!\!\! \mod m}\\
\\
\displaystyle \widetilde{C}_k(\a) = 
\widetilde{C}_0( \omega^k.\a) = \omega^{m-2+2r(\omega^k.\a)}, \hspace{5mm} 
\widetilde{C}_k = \widetilde{C}_{k \!\!\! \mod m}
\end{array}
\right. .
\end{equation}
\end{itemize}
\end{thm}

For any $k \in \mathbb{Z}$, the analytic continuation of 
 $\{ \Phi_{k-1},\Phi_{k} \}$  constitutes a
fundamental system of solutions of $(\mathfrak{E}_m)$ (by lemma
 \ref{lemme2}, the Wronskian $W(\Phi_{k-1},\Phi_{k})$ does not
 vanish). 
In particular, there exists a unique invertible
$2 \times 2$  matrix $\mathfrak{M}_k^\infty (\a)$, entire in $\a$, such
 that  
$\displaystyle 
\left(\begin{array}{c} 
\Phi_{k-1} \\ 
\Phi_k \end{array} 
\right) (\omega^m x, \a)  = \mathfrak{M}_k^\infty (\a)
\left(\begin{array}{c} 
\Phi_{k-1} \\
\Phi_{k} 
\end{array} 
\right) (x, \a)$.  

\begin{defn}
The $2 \times 2$ matrices  
$\mathfrak{M}_k^\infty (\a)$, $k \in \mathbb{Z}$,  defined by
\begin{equation}\label{Monoinfty1}
\displaystyle 
\left(\begin{array}{c} 
\Phi_{k-1} \\ 
\Phi_k \end{array} 
\right) (\omega^m x, \a)  = \mathfrak{M}_k^\infty (\a)
\left(\begin{array}{c} 
\Phi_{k-1} \\
\Phi_{k} 
\end{array} 
\right) (x, \a),
\end{equation}
are called the $\infty$-monodromy matrices.
\end{defn}

From the very definition of the $\Phi_k$'s, we note
that  $\displaystyle \mathfrak{M}_k^\infty (\a) = \mathfrak{M}_0^\infty
(\omega^k .\a)$. Also, 
$$\displaystyle 
\left(\begin{array}{c} 
\Phi_{-1} \\ 
\Phi_0 \end{array} 
\right) ( x, \a) = \mathfrak{S}_0 (\a) \cdots  \mathfrak{S}_{m-1} (\a)
\left(\begin{array}{c} 
\Phi_{m-1} \\ 
\Phi_m \end{array} 
\right) (x, \a)$$
$$=
\mathfrak{S}_0 (\a) \cdots  \mathfrak{S}_{m-1} (\a)
\left(\begin{array}{c} 
\Phi_{-1} \\ 
\Phi_0 \end{array} 
\right) (\omega^m x, \omega^m.\a).$$
Since $\displaystyle \left(\begin{array}{c} 
\Phi_{-1} \\ 
\Phi_0 \end{array} 
\right)$ is entire in $\a$, we obtain
$$\left(\begin{array}{c} 
\Phi_{-1} \\ 
\Phi_0 \end{array} 
\right) ( x, \a)= \mathfrak{S}_0 (\a) \cdots  \mathfrak{S}_{m-1} (\a)
\left(\begin{array}{c} 
\Phi_{-1} \\ 
\Phi_0 \end{array} 
\right) (\omega^m x, \a)
$$
$$ = \mathfrak{S}_0 (\a) \cdots  \mathfrak{S}_{m-1} (\a)
\mathfrak{M}_0^\infty (\a)
\left(\begin{array}{c} 
\Phi_{-1} \\ 
\Phi_0 \end{array} 
\right) (x, \a),$$
and  $\{ \Phi_{-1},\Phi_{0} \}$ being a fundamental system, this
yields:

\begin{thm}\label{Monoinfty}
For every $k \in \mathbb{Z}$, 
the $\infty$-monodromy matrix $\mathfrak{M}_k^\infty (\a)$ is invertible, 
entire in $\a$, and 
\begin{equation}\label{Monoinfty2}
\mathfrak{M}_k^\infty (\a) = \mathfrak{M}_0^\infty (\omega^k .\a).
\end{equation}
Furthermore, the Stokes-Sibuya connection matrices satisfy the
functional relation:
\begin{equation}\label{Monoinfty3}
\mathfrak{S}_0 (\a) \cdots  \mathfrak{S}_{m-1} (\a) =
\left({\mathfrak{M}_0^\infty}\right(\a))^{-1}. 
\end{equation}
\end{thm}

Relation (\ref{Monoinfty3}) generalized a functional relation due to
Sibuya \cite{Sib75}, p. 85. Unfortunately, as a rule, the $\infty$- monodromy
matrix $\mathfrak{M}_0^\infty$ is difficult to compute. We 
return to this question in section \ref{sec5}.

\section{Solutions of $(\mathfrak{E}_m)$ in the neighbourhood of
  infinity: resurgent point of view }\label{Resurgentview}

Theorem \ref{Sibuya} can be shown with the methods developed in
Sibuya's book \cite{Sib75} and, in  fact, 
theorem \ref{Sibuya} is actually proved in \cite{Mu68, Bakk}.
 However, using the resurgent viewpoint, 
one can get a stronger result in a simpler way.

\subsection{Basic notions in resurgence theory}\label{baseres}

Since the terminologies we shall use in this section are likely to be least familiar to
the readers, we first introduce the necessary definitions,
cf. \cite{CNP1, CNP2, DP97, Ec81-1} for more details. We mention that
our notations differ from those usually used by Ecalle. \\
As usual in this article, we identify an element of the  universal covering of
$\mathbb{C}^\star$ (with base point 1) by specifying its argument in $\mathbb{R}$.

\begin{defn}
A {\sl sectorial neighbourhood of infinity of aperture $I=]\alpha, \beta[
\subset \mathbb{R}$} is an open set $U$ of the universal covering of
$\mathbb{C}^\star$ such that for any open interval $J \subset I$, there
is $z \in U$ such that $zJ \subset U$, where 
$$\displaystyle zJ :=
\{z+re^{i\theta}, \, r>0, \, \theta \in J\}.$$ 
\end{defn}

\begin{defn}
If $U$ is a sectorial neighbourhood of infinity of aperture $I$ and if 
$\Psi$ is holomorphic in $U$, $\Psi$ is {\sl of exponential growth of order
1 at infinity in $U$} if for any open interval $J \subset I$, there
exist $\tau >0$ and $C>0$  such that
$$\forall z \in U \cap 0J, \, |\Psi(z)| \leq Ce^{\tau |z|}.$$
\end{defn}

We now introduce the notion of minor. We do that only for a class of formal power
series which will be used in this paper.

\begin{defn}
We consider the formal power series  $\displaystyle \psi(z) = r +
\sum_{n=0}^\infty \frac{\alpha_n}{z^{\frac{n}{m}+1}} \in
\mathbb{C}[[z^{-\frac{1}{m}}]]$, where $m$ is a positive
integer. Then, $r$ is the {\sl   residual coefficient} of $\psi$ and 
$\displaystyle {\widetilde \psi} (\zeta) = 
\sum_{n=0}^\infty \frac{\alpha_n}{\Gamma(\frac{n}{m}+1)}\zeta^\frac{n}{m} \in
\mathbb{C}[[\zeta^{\frac{1}{m}}]]$ is the {\sl minor} of $\psi$.
\end{defn}

In other words, the minor of the formal power series $\psi$ is nothing
but its Borel transform when forgetting its residual coefficient.

This allows to define the  Borel-summability for such formal  series.

\begin{defn}\label{Borel}
The formal power series  $\displaystyle \psi(z) = r +
\sum_{n=0}^\infty \frac{\alpha_n}{z^{\frac{n}{m}+1}} \in
\mathbb{C}[[z^{-\frac{1}{m}}]]$ is {\sl Borel-resummable in the
direction (of argument) $\alpha \in \mathbb{R}$} if:
\begin{enumerate}
\item its minor $\displaystyle {\widetilde \psi} (\zeta)$ defines a
(ramified) analytic function at the origin\footnote{For simplicity, we keep the same
  notation $\displaystyle {\widetilde \psi} (\zeta)$ for the series
  and its sum, and $\zeta$ should be seen as an element of the
  universal covering of $\mathbb{C}^\star$.}, 
\item \label{cond2} there exists an open sector $0I$ with $I$ an open neighbourhood of
$\alpha$ such that $\displaystyle {\widetilde \psi} (\zeta)$ can be analytically
extended in $0I$ and is of exponential growth of order 1 at infinity in $0I$.
\end{enumerate}
The {\sl Borel-sum}  $\mbox{\sc s}_{\alpha} \psi (z)$  with respect to
$z$ in the direction  $\alpha \in \mathbb{R}$ of the
formal series $\psi$ is defined by 
$$\mbox{\sc s}_{\alpha} \psi (z):= r+ \int_0^{\infty e^{i\alpha}} 
e^{-z\zeta} {\widetilde \psi} (\zeta) d\zeta.$$
\end{defn}

In definition \ref{Borel}, when one drops the growth condition at
  infinity (condition \ref{cond2}), then $\displaystyle \psi$ is said to be {\sl Borel
  presummable} in the direction $\alpha$, the summation operator $\mbox{\sc s}_{\alpha}$ being
  replaced by the {\sl presummation} operator which we do not define here, see, e.g., \cite{DP97}.\\

A Borel sum has the following main properties:

\begin{prop}\label{propborel}
If  $\displaystyle \psi(z) = r +
\sum_{n=0}^\infty \frac{\alpha_n}{z^{\frac{n}{m}+1}} \in
\mathbb{C}[[z^{-\frac{1}{m}}]]$ is Borel-resummable in the
direction $\alpha \in \mathbb{R}$, then:
\begin{itemize}
\item its Borel sum  $\mbox{\sc s}_{\alpha} \psi (z)$ is holomorphic
  in a sectorial neighbourhood of infinity $U$ of aperture
  $\displaystyle I= ]-\frac{\pi}{2} -\alpha, \frac{\pi}{2} -\alpha[$.
\item $\mbox{\sc s}_{\alpha} \psi (z)$ is asymptotic  to $\psi(z)$ at
  infinity in $U$. 
More precisely, 
for any strict subinterval $J$  of $I$, there is $C>0$ such that,
for all $n \geq 1$, for all $z \in U \cap 0J$, 
$\displaystyle |\mbox{\sc s}_{\alpha} \psi (z) - r - \sum_{k=0}^{n-1} \frac{\alpha_k}{z^{\frac{k}{m}+1}}|
\leq C^n\Gamma(\frac{n}{m}+1)|z|^{-\frac{n}{m}-1}.$  
\item $\displaystyle \frac{d}{dz} \left(\mbox{\sc s}_{\alpha} \psi (z)\right) =
  \mbox{\sc s}_{\alpha} \left( \frac{d\psi}{dz}(z) \right)$.
\item If two   formal power series  $\displaystyle \psi(z), \, \phi(z)
  \in \mathbb{C}[[z^{-\frac{1}{m}}]]$ are Borel-resummable in the
direction $\alpha \in \mathbb{R}$, then
$\displaystyle \mbox{\sc s}_{\alpha} \left(\psi.\phi \right)(z)  = \mbox{\sc
  s}_{\alpha} (\psi )(z) .  \mbox{\sc
  s}_{\alpha} (\phi )(z).$
\end{itemize}
\end{prop}

\begin{defn}
A formal power series  $\displaystyle \psi(z) = r +
\sum_{n=0}^\infty \frac{\alpha_n}{z^{\frac{n}{m}+1}} \in
\mathbb{C}[[z^{-\frac{1}{m}}]]$ is {\sl resurgent} if its minor 
$\displaystyle {\widetilde \psi} (\zeta)$ defines a
(ramified) analytic function at the origin and is 
{\sl  endlessly continuable}, i.e., for every $L > 0$ there is a
finite  subset  $\Omega_L \subset \mathbb{C}$ such that
${\widetilde \psi}$ can be analytically continued along every path $\lambda$ of length
$< L$  which avoids $\Omega_L$.
\end{defn}

This definition\footnote{Note that Ecalle proposes a more general
 definition.}
 can be extended to an algebra of {\sl resurgent formal
 functions} which we do not precise here. 

\begin{prop}
If two  formal power series  $\displaystyle \psi(z) = r +
\sum_{n=0}^\infty \frac{\alpha_n}{z^{\frac{n}{m}+1}} \in
\mathbb{C}[[z^{-\frac{1}{m}}]]$ and $\displaystyle \phi(z) = s +
\sum_{n=0}^\infty \frac{\beta_n}{z^{\frac{n}{m}+1}} \in
\mathbb{C}[[z^{-\frac{1}{m}}]]$ are resurgent, then the product
$\displaystyle \psi.\phi(z)$ is also resurgent: the minor of
$\displaystyle \psi.\phi(z)$, which is given by
\begin{equation}\label{convo}
\begin{array}{l}
\displaystyle \widetilde{\psi.\phi}(\zeta) = r.{\widetilde \phi}
(\zeta)+s.{\widetilde \psi} (\zeta) + {\widetilde \psi}\ast
{\widetilde \phi} (\zeta),\\
\\
\displaystyle {\widetilde \psi}\ast
{\widetilde \phi} (\zeta) = \int_0^\zeta {\widetilde \psi}(\eta)
{\widetilde \phi}(\zeta -\eta)d\eta \mbox{ (the convolution product),}
\end{array}
\end{equation}
is endlessly continuable.
\end{prop}

For a resurgent formal power series, it may happen that we no longer
can define its Borel (pre)sum in a given
direction $\alpha \in \mathbb{R}$ because the minor can have
singularities along this direction: this is the essence of the  {\sl Stokes
  phenomenon}. 

\begin{figure}[thp]
\begin{center}
\rotatebox{-90}{\includegraphics[height=110mm]{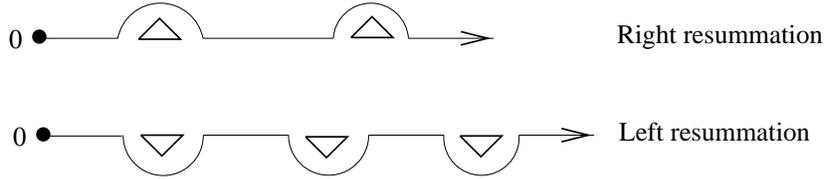}}
\caption{The integration path for {\sl right} (resp. 
{\sl left}) resummation (for $\alpha=0$).
\label{fig:DP13}}
\end{center}
\end{figure}

\begin{defn}\label{rlsum}
We consider a resurgent formal power series  $\displaystyle \psi(z) = r +
\sum_{n=0}^\infty \frac{\alpha_n}{z^{\frac{n}{m}+1}} \in
\mathbb{C}[[z^{-\frac{1}{m}}]]$. Let  $\alpha
\in \mathbb{R}$ be a singular direction for the minor ${\widetilde
  \psi} (\zeta)$. \\
Hypothesis: we assume that there is $\epsilon >0$ such that  ${\widetilde
  \psi} (\zeta)$ can be analytically extended in the open sector
$0]\alpha, \alpha +\epsilon[$ (resp. $0]\alpha -\epsilon, \alpha[$
with an exponential growth of order 1 at infinity. We assume also that
this exponential growth at infinity extends up 
to a path  $[0, \infty e^{i\alpha}+[$ (resp. $[0,
\infty e^{i\alpha}-[$) which circumvents the singularities to the left
(resp. right) along the direction $\alpha$, see figure \ref{fig:DP13}. \\
Then $\psi$ is {\sl right}
(resp. {\sl left}) {\sl  Borel-resummable} in the direction $\alpha$, its
right (resp. left) Borel sum $\mbox{\sc s}_{\alpha+} \psi$ (resp. $\mbox{\sc s}_{\alpha-} \psi$)
being defined by
$$\mbox{\sc s}_{\alpha \pm} \psi (z):= r+ \int_0^{\infty e^{i\alpha} \pm} 
e^{-z\zeta} {\widetilde \psi} (\zeta) d\zeta,$$
for $z$ in a sectorial
neighbourhood of infinity of aperture   
$\displaystyle ]-\frac{\pi}{2} -\alpha, \frac{\pi}{2}
-\alpha[$.
\end{defn}

In definition \ref{rlsum}, it is possible to drop the Hypothesis, 
  replacing right (resp. left) Borel sum by {\sl right} (resp. {\sl
  left}) {\sl Borel presum}, see, e.g., \cite{DP97}. In other words,
  every resurgent formal function is always right and left
  Borel-presummable (in any direction).

We note that proposition \ref{propborel} is still valid for right and
left Borel sum.
Moreover,  when $\psi$ is Borel (pre)summable in the
direction $\alpha \in \mathbb{R}$, then
$$\mbox{\sc s}_{\alpha} \psi (z) = \mbox{\sc s}_{\alpha+} \psi (z) =
\mbox{\sc s}_{\alpha-} \psi (z).$$

In order to be able to compare right and left (pre)summation, one has to
enlarge the set of resurgent formal functions to the set of resurgent symbols. 

\begin{defn}\label{ressymbol}
A {\sl resurgent symbol}\footnote{or {\sl resurgent transseries.}} in the direction $\alpha$ is a formal sum 
$$ \varphi (z) = \sum_{\omega \in \Omega} \psi_\omega \, 
e^{-z\omega}$$
where each $\psi_\omega (z)$ is a  resurgent formal function and $\Omega$,
the {\sl singular support of $\varphi$}, is
a  discrete subset of $[0, \infty e^{i\alpha}[$.\\
The sum and product of two resurgent symbols are defined in an obvious
fashion, so that resurgent symbols in the direction $\alpha$
make up an algebra which we denote
by  ${\mathcal R}_\alpha$.
\end{defn}

The {\sl right} (resp. {\sl left}) (pre)summation operations can be
extended to resurgent symbols in a way so that
$$\mbox{\sc s}_{\alpha+} \, \varphi  
=  \sum_{\omega \in \Omega} (\mbox{\sc s}_{\alpha+} \,  
\psi_\omega) \, e^{-z\omega} \hspace{5mm} \mbox{resp.} \hspace{5mm} \mbox{\sc s}_{\alpha-} \, \varphi 
=  \sum_{\omega \in \Omega} (\mbox{\sc s}_{\alpha-} \,  
\psi_\omega) \, e^{-z\omega}. $$
The construction (which we do not explain here) makes the operations 
$\mbox{\sc s}_{\alpha+}$ and $\mbox{\sc s}_{\alpha-}$\label{ind39}  
{\sl isomorphisms of algebras} and, moreover, $\mbox{\sc s}_{\alpha+}
({\mathcal R}_\alpha) = \mbox{\sc s}_{\alpha-}({\mathcal
  R}_\alpha)$. This  key-result (due to Ecalle) allows to define the
so-called Stokes automorphism, which analyzes the Stokes
phenomenon by  comparing right and left Borel-(pre)summations:

\begin{defn}\label{Stokes}
The {\sl Stokes automorphism in the direction $\alpha$} is defined by:
$${\mathfrak S}_\alpha := \mbox{\sc s}_{\alpha+}^{-1} \circ 
\mbox{\sc s}_{\alpha-} \, : \,
{\mathcal R}_\alpha \rightarrow {\mathcal R}_\alpha .$$
\end{defn}

\begin{figure}[thp]
\begin{center}
\rotatebox{-90}{\includegraphics[height=95mm]{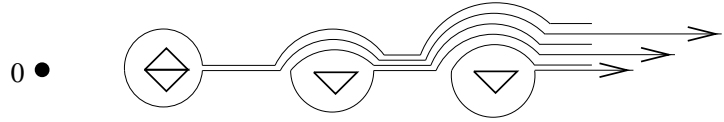}}
\caption{The difference between right and left
summations (Fig. \ref{fig:DP13}) is described as a Laplace integral
along a sum of
contours. 
\label{fig:DP17}}
\end{center}
\end{figure}

\noindent Remark: the action of the Stokes automorphism in a given direction can be understood 
in terms of deformation of the contour of integration in a Laplace
integral, see Figure \ref{fig:DP17}.

It follows from the definitions that the Stokes
 automorphism in the direction $\alpha$ acts  trivially on  exponentials $e^{-z\omega}$, and
 that its  action on a formal resurgent series $\psi$ reads:
$${\mathfrak S}_\alpha \psi = \psi + \sum_\omega  \psi_\omega 
e^{-z\omega},$$
where the sum runs over those singular points of the minor
${\widetilde \psi}$ which have to be avoided when considering left
(pre)summation. The Stokes
 automorphism ${\mathfrak S}_\alpha$ reads as
$${\mathfrak S}_\alpha = 1\! \! 1 + {^+  \! S \! \! \! S}_\alpha $$
where the operator ${^+  \! S \! \! \! S}_\alpha $
 commutes with multiplication
by  exponentials, and  transforms formal resurgent series into
 ``exponentially small resurgent symbols''. 
This implies that the operator
$${\underline {\dot \Delta}}_\alpha  := \ln {\mathfrak S}_\alpha = 
\sum_{n=1}^{\infty} \frac{(-1)^{n-1}}{n} {^+ \! S \! \! \! S}_\alpha^n  $$
is well defined  on ${\mathcal R}_\alpha$. Since ${\mathfrak S}_\alpha$ is an
automorphism of ${\mathcal R}_\alpha$, ${\underline {\dot
    \Delta}}_\alpha$  is a 
derivation of ${\mathcal R}_\alpha$.

\begin{defn}\label{deltal}
 ${\underline {\dot
    \Delta}}_\alpha$ is called the {\sl alien 
derivation in the direction $\alpha$}.
\end{defn}

We note that since $\mbox{\sc s}_{\alpha+}$ and $\mbox{\sc
  s}_{\alpha-}$ commutes with $\frac{d}{dz}$, ${\mathfrak S}_\alpha$
  and therefore ${\underline {\dot
    \Delta}}_\alpha$ commutes with $\frac{d}{dz}$.

\begin{prop}\label{commu}
 ${\underline {\dot
    \Delta}}_\alpha$ commutes with $\frac{d}{dz}$.
\end{prop}

The alien derivation  ${\underline {\dot
    \Delta}}_\alpha$  commutes with
multiplication  by exponentials, and its action on a formal resurgent 
series $\psi$ has the following explicit form:
$${\underline {\dot \Delta}}_\alpha \,  \psi = 
\sum_{\omega \in \Omega ({\widetilde \psi})} 
 \, e^{-z \omega} \, \Delta_\omega \psi,$$
where $\Omega ({\widetilde \psi})$ is a discrete subset of $[0, \infty
e^{i\alpha}[$, the set of so-called {\sl glimpsed singularities} of 
${\widetilde \psi}$, i.e.,  the set of singularities  to be
circumvented when analytically continuing ${\widetilde \psi}$ in the
direction $\alpha$. Since ${\underline {\dot
    \Delta}}_\alpha$ is a derivation, each  $\Delta_\omega$  is also
a derivation. 

\begin{defn}
$\Delta_\omega$ is called the {\sl alien  derivation at  $\omega$}.
\end{defn}

\begin{defn}\label{constres}
A formal resurgent function $\psi$ is a {\sl constant of resurgence}
if for any $\omega$, $\Delta_\omega \psi =0$.
\end{defn}

For instance, if $\psi$ is a convergent power series, then $\psi$ is a
constant of resurgence.\\

Instead of working with the $\Delta_\omega$ 's, it is sometimes
convenient  to work with the so-called {\sl pointed alien derivations}
$${\dot \Delta}_\omega  \, := \,  e^{-z\omega}\,  \Delta_\omega.$$
They have the advantage  of commuting with $\displaystyle
\frac{d}{dz}$ (this is a consequence of  proposition \ref{commu}).

\begin{prop}\label{commute}
The pointed alien derivations ${\dot \Delta}_\omega$ commute with $\frac{d}{dz}$.
\end{prop}

\subsection{Resurgence of solutions of $(\mathfrak{E}_m)$}

We now return to the proof of theorem \ref{Sibuya}. The main idea
will be to consider the asymptotic expansion  $T\Phi_0$ of theorem \ref{Sibuya} as a formal
solution of equation $(\mathfrak{E}_m)$, and to show that $\Phi_0$ can
be   constructed from  $T\Phi_0$
by Borel resummation with respect to an appropriate variable $z =z(x,\a)$ (or
${\widetilde z}$, depending on the parity of $m$) which will be defined
later, uniformly in $\a$ for $\a$ in any given  compact set.

We start with a kind of
preparation theorem, so as to transform equation $(\mathfrak{E}_m)$ into
a normal form. This is based on the Green-Liouville transformation,
but we shall have to dissociate the $m$ odd case from the $m$ even
case for technical reasons. 

In what
follows, $\a$ is assumed to belong to an arbitrary  
given compact set $K \subset \mathbb{C}^m$.

\subsubsection{The $m$ odd case}

We consider the Green-Liouville transformation~:
\begin{equation}\label{Green1}
 \left\lbrace \begin{array}{l}

\displaystyle z=z(x,\a)=\int^x \frac{\sqrt{P_m(t,\a)}}{t}  dt\\
\\
\displaystyle \Psi(z,\a)=\frac{P_m(x,\a)^{\frac{1}{4}}}{\sqrt{x}} \Phi(x,\a)

\end{array} \right. 
\end{equation}
where the Laurent-Puiseux series  expansion in $x$ of $z(x,\a)$,
$$z(x,\a) = \frac{2}{m}x^{\frac{m}{2}}+\sum_{k=1}^{\frac{m-1}{2}}
\frac{b_{\frac{m}{2}-k}(\a)}{\frac{m}{2}-k} 
x^{\frac{m}{2}-k} + O(x^{-\frac{1}{2}}) 
\in x^\frac{m}{2}\mathbb{C}[\a]\{x^{-1}\},$$
coincides with the map $x \mapsto S(x,\a)$ defined in
theorem \ref{Sibuya} modulo an analytic (multivalued) function 
vanishing at infinity. We note that, if $\displaystyle 
\begin{array}{l}
X\\
\downarrow \pi\\
\mathbb{C}
\end{array}
$ ({\sl resp. } $\displaystyle 
\begin{array}{l}
Z\\
\downarrow \widetilde{\pi}\\
\mathbb{C}
\end{array}
$)
is the (ramified) Riemann surface of $x^{1/2}$  ({\sl resp. }
$z^{1/m}$), then there exists a
compact set $B$ (depending on $K$) ({\sl resp. } $\widetilde{B}$) 
such that the map $(x,\a) \in
\pi^{-1} (\mathbb{C} \backslash B) \times K \mapsto (z(x,\a),\a)\in
\widetilde{\pi}^{-1} (\mathbb{C} \backslash \widetilde{B}) \times K$ is
bi-holomorphic. 

Remark here that by quasi-homogeneity of
$P_m(x,\a)$,
\begin{equation}\label{qhz1}
z(\omega x,\omega .\a)=
\omega^\frac{m}{2}z(x, \a).
\end{equation}
The transformation (\ref{Green1})
converts $(\mathfrak{E}_m)$ into the following equation, 
\begin{equation}\label{mathfrakF}
 - \frac{d ^2}{d z^2}\Psi + (1-F(z,\a))\Psi=0,
\end{equation}
which is our prepared normal form.

It is straightforward to see that 
 studying $(\mathfrak{E}_m)$ at infinity in the variable $x$
amounts to studying (\ref{mathfrakF}) at infinity in the variable
 $z$. The inverse map $x: (z,\a) \mapsto x(z,\a)$ can be identified with
  its Laurent-Puiseux series expansion 
\begin{equation}\label{inverse1}
 \displaystyle x(z,\a) =
 \left(\frac{m}{2}z\right)^\frac{2}{m} +O(1) \in 
z^\frac{2}{m}\mathbb{C}[\a]\{z^{-\frac{2}{m}}\},
\end{equation}
and from (\ref{qhz1}),
\begin{equation}\label{qhx1}
x(\omega^\frac{m}{2} z,\omega .\a)=
\omega x(z, \a).
\end{equation}
In (\ref{mathfrakF}), $F(z,\a)$ is defined by:
\begin{equation}\label{cestF} 
\displaystyle F(z,\a)=-\frac{xP_m'(x,\a)+P_m(x,\a)}{4P_m(x,\a)^2}-x^2\big(\frac{1}{4}
 \frac{P_m''(x,\a)}{P_m(x,\a)^2}-\frac{5}{16}
\frac{(P_m'(x,\a))^2}{P_m(x,\a)^3}\big) \, |_{\displaystyle x=x(z,\a)}.
\end{equation}
One infers from (\ref{inverse1}) that
\begin{equation}\label{estimateF}
\displaystyle 
F(z,\a) = \frac{m^2-4}{4m^2z^2}(1+O(z^{-\frac{2}{m}}))
 \in \frac{1}{z^2} \mathbb{C}[\a]\{z^{-\frac{2}{m}}\}
\end{equation}
is an analytic 
function at infinity in $z^{-1/m}$, uniformly in $\a \in K$.
It is easy to show the existence of a unique formal
solution $\Psi_+(z,\a)$ of (\ref{mathfrakF}) satisfying 
\begin{equation}\label{psiplus}
\Psi_+(z,\a) = e^{-z} \psi_+(z,\a),
\end{equation}
where
\begin{equation}\label{psiplusbis}
\psi_+(z,\a) = 1 + \sum_{n=0}^\infty \frac{\alpha_n (\a)}{z^{\frac{n}{m}+1}}\in
\mathbb{C}[\a][[z^{-\frac{1}{m}}]]
\end{equation}
 with $1$ for its  residual
coefficient. Moreover the formal power series
expansion $\psi_+(z,\a)$ satisfies the 
following ordinary differential equation:
\begin{equation}\label{nonconvequat}
 - \frac{d ^2}{d z^2}\psi_+ + 2\frac{d}{d z}\psi_+ -F(z,\a)\psi_+=0.
\end{equation}

In addition, from the quasi-homogeneity property of $P_m$, from
(\ref{qhx1}) et (\ref{cestF}), one sees that
\begin{equation}\label{qhF1}
F(\omega^\frac{m}{2} z , \omega . \a)= F(z, \a).
\end{equation}
Defining
\begin{equation}\label{cestphiminus}
\psi_-(z, \a) =  
\psi_+(\omega^\frac{m}{2}z, \omega . \a),
\end{equation}
one deduces the existence of a unique formal 
solution  $\Psi_-(z,\a)$ of (\ref{mathfrakF}) such that 
$\Psi_-(z,\a) = e^{z} \psi_-(z,\a)$ with  $\psi_- \in
\mathbb{C}[\a][[z^{-\frac{1}{m}}]]$ with a residual coefficient equal to
$1$. Note that $\Psi_+,\Psi_-$ are  linearly independent, 
 so that $(\Psi_+, \Psi_-)$ provides a fundamental system of
formal solutions for the linear second order equation
(\ref{mathfrakF}).

The formal series $\psi_\pm$ enjoys the following properties:

\begin{prop}\label{propborel}
The formal power series expansion 
$\displaystyle \psi_+(z,\a) = 1 + 
\sum_{n=0}^\infty \frac{\alpha_n (\a)}{z^{\frac{n}{m}+1}}\in
\mathbb{C}[\a][[z^{-\frac{1}{m}}]]$ ({\sl resp.} $\displaystyle
\psi_-(z,\a) = \psi_+(\omega^\frac{m}{2}z, \omega . \a)$)
  is Borel-resummable with respect
to $z$, uniformly in $\a$ for $\a$ in any compact set of
$\mathbb{C}^m$, for every direction of summation except those of
argument  $\pi \mod (2\pi)$ ({\sl resp.} $0 \mod (2\pi)$).
\end{prop}

\begin{proof}
We have to  analyze the analytic properties of the 
 minor 
\begin{equation}\label{minor}
\widetilde \psi_+(\zeta, \a) = 
\sum_{n=0}^\infty \frac{\alpha_n (\a)}{\Gamma(\frac{n}{m}+1)}
\zeta^{\frac{n}{m}} \in
\mathbb{C}[\a][[\zeta^{\frac{1}{m}}]].
\end{equation}
of $\psi_+(z,\a)$.
To proceed, we go back to equation (\ref{nonconvequat}). Instead of
 considering this differential equation, 
 we shall rather introduce its deformation,
\begin{equation}\label{nonconvequatbis}
 - \frac{d ^2}{d z^2}\psi + 2\frac{d}{d z}\psi -F(z,\a)
   +\varepsilon F(z,\a)(1- \psi)=0,
\end{equation}
where $\varepsilon$ can be thought of as a parameter of
perturbation. The introduction of this parameter will help us to
rewrite $\psi_+$ and its minor $\widetilde \psi_+$ into an analyzable 
form, since (\ref{nonconvequatbis}) reduces to
(\ref{nonconvequat}) when $\varepsilon=1$. We now look for a formal
solution of (\ref{nonconvequatbis}) in the form of a normalized 
series expansion with respect to  $\varepsilon$:
\begin{equation}\label{perturb}
\psi(z, \a, \varepsilon) = 1+ \sum_{n=0}^\infty \psi_n (z,
\a)\varepsilon^n, \hspace{5mm} \psi_n \in 
\frac{1}{z}\mathbb{C}[\a][[z^{-\frac{1}{m}}]].
\end{equation}
Plugging (\ref{perturb}) into (\ref{nonconvequatbis}) and identifying
the powers of $\varepsilon$, one gets:
\begin{equation}\label{rec1}
\left\{
\begin{array}{l}
\displaystyle  - \frac{d ^2}{d z^2}\psi_0 + 2\frac{d}{d z}\psi_0
=F(z,\a)\\
\\
\displaystyle  - \frac{d ^2}{d z^2}\psi_{n+1} + 2\frac{d}{d z}\psi_{n+1}
=F(z,\a)\psi_{n}, \, \, n \geq 0 .
\end{array}
\right.
\end{equation}
This translates into the fact that the minors ${\widetilde
  \psi_n}(\zeta,\a)$ of  the $\psi_n(z,\a)$ have to satisfy the
  following convolution equations,
\begin{equation}\label{rec2}
\left\{
\begin{array}{l}
\displaystyle -\zeta(2+\zeta){\widetilde
  \psi_0} = {\widetilde F}\\
\\
\displaystyle  -\zeta(\zeta+2) {\widetilde \psi}_{n+1} =
{\widetilde \psi}_{n} \ast  {\widetilde F} , \, \, n \geq 0 ,\\
\end{array}
\right.
\end{equation}
where ${\widetilde F}(\zeta, \a)$ is the minor of $F(z,\a)$, while
$\ast$ stands for the convolution product (cf. (\ref{convo})). \\
We have now to analyze the analytic properties of
\begin{equation}\label{perte}
{\widetilde \psi} (\zeta, \a, \varepsilon) = \sum_{n=0}^\infty 
{\widetilde \psi_n} (\zeta, \a) \varepsilon^n .
\end{equation}
The key-point of the analysis will come from the properties of $F$. Writing
\begin{equation}\label{cestFbis}
\displaystyle F(z,\a)  = 
\sum_{n=0}^\infty \frac{f_n(\a)}{z^{\frac{2n}{m}+2}} \in 
\frac{1}{z^2} \mathbb{C}[\a]\{z^{-\frac{2}{m}}\},
\end{equation}
we know that
\begin{equation}\label{cestG}
\displaystyle G(z)  = 
\sum_{n=0}^\infty \frac{g_n}{z^{\frac{2n}{m}+2}} \hspace{5mm}
\mbox{with} \hspace{5mm} g_n = \sup_{\a \in K} |f_n(\a)|,   
\end{equation}
is an analytic function at infinity in $z^{-1/m}$. Therefore, 
 its minor
\begin{equation}\label{cestGtilde}
\displaystyle \widetilde  G (\zeta)  = 
\sum_{n=0}^\infty \frac{g_n}{\Gamma(\frac{2n}{m}+2)} 
\zeta^{\frac{2n}{m}+1} \in  \zeta\mathbb{C}\{\zeta^{\frac{2}{m}}\}
\end{equation}
 is an
entire function in $\zeta^{1/m}$  (with an exponential growth at infinity of
order at most $1$). Thus, 
 if $\mathfrak{C}_m$
denotes the Riemann surface of $\zeta^{1/m}$, then 
$${\widetilde F}  (\zeta, \a)  = 
\sum_{n=0}^\infty \frac{f_n(\a)}{\Gamma(\frac{2n}{m}+2)} 
\zeta^{\frac{2n}{m}+1} \in
\zeta\mathbb{C}[\a]\{\zeta^{\frac{2}{m}}\}$$
 is
a holomorphic function in $(\zeta, \a) \in \mathfrak{C}_m \times K$
such that 
\begin{equation}\label{major}
\forall (\zeta, \a) \in \mathfrak{C}_m \times K, \, \, 
|\widetilde F(\zeta, \a) | \leq G(|\zeta|) .
\end{equation}
Using the fact that ${\widetilde F}$ is a  holomorphic function in
$(\zeta, \a) \in \mathfrak{C}_m \times K$  such that $F(\zeta, \a) =
O(\zeta)$ uniformly in $\a \in K$, and from the  properties
of the convolution product, one easily deduces
from (\ref{rec2}) that each $\widetilde \psi_n$ belongs to the space
$\mathbb{C}[\a]\{\zeta^{\frac{1}{m}}\}$ and extends analytically to
$\widetilde{\mathbb{C} \backslash \{0,-2\}} \times K$, where 
$\widetilde{\mathbb{C} \backslash \{0,-2\}}$ is the universal covering
of $\mathbb{C} \backslash \{0,-2\}$.\\
For $\rho >0$, we now define  the star-shape domain
\begin{equation}
\Omega_m (\rho) = \{\zeta \in \mathfrak{C}_m, \, \, |\dot \zeta
+2|>\rho, \, \, [0, \zeta] \in \Omega_m (\rho)\}
\subset \mathfrak{C}_m
\end{equation}
where $\dot \zeta$ is the projection of $\zeta$ by the natural mapping 
$\mathfrak{C}_m \rightarrow \mathbb{C}$.  We also introduce
the sequence of analytic functions $h_n(\zeta)$ defined for $\zeta \in
\mathfrak{C}_m$ by:
\begin{equation}\label{rec3}
\left\{
\begin{array}{l}
\displaystyle \zeta \rho {\widetilde
  h_0} = {\widetilde G}\\
\\
\displaystyle  \zeta \rho {\widetilde h}_{n+1} =
{\widetilde h}_{n} \ast  {\widetilde G} , \, \, n \geq 0\\
\end{array}
\right..
\end{equation}
Comparing (\ref{rec3}) with (\ref{rec2}), and using (\ref{major}), one gets
\begin{equation}\label{major2}
\forall (\zeta, \a) \in \Omega_m(\rho) \times K, \, 
\forall n \in \mathbb{N}, 
\, \, 
|{\widetilde \psi}_{n}(\zeta, \a) | \leq {\widetilde h}_{n} (|\zeta|).
\end{equation}
This can be shown by an easy recursion. We just detail here the $n=0$
and $n=1$ cases.\\
For all $(\zeta, \a) \in \left(\Omega_m(\rho)\backslash \{0\}
\right)\times K$ we first have:
$$|\widetilde \psi _0(\zeta, \a)| = 
\frac{|\widetilde F (\zeta, \a)|}{|\zeta||\zeta +
  2|}\leq \frac{\widetilde G (|\zeta|)}{|\zeta|\rho} = \widetilde h_0
(|\zeta|),$$
and this inequality extends to $\zeta=0$ by continuity.
This proves (\ref{major2}) for $n=0$.\\
We thus deduce that, for all $(\zeta, \a) \in \left(\Omega_m(\rho)\backslash \{0\} \right)\times K$:
$$|\widetilde \psi _1(\zeta, \a)| = \frac{|\widetilde F  * \widetilde \psi
  _0(\zeta, \a)|}{|\zeta||\zeta +  2|} \leq
 \frac{|\int_0^{\zeta}\widetilde F (\eta, \a)
 \widetilde \psi  _0(\zeta - \eta, \a) d\eta |}{|\zeta|\rho}.$$
Writing $\displaystyle \zeta = |\zeta|e^{i\theta}$ and making the
  change of variable $\eta = te ^{i\theta}$, we get:
$$\left |\int_0^{\zeta}\widetilde F (\eta, \a)
 \widetilde \psi  _0(\zeta - \eta, \a) d\eta \right | = \left | \int_0^{|\zeta|}\widetilde F (te ^{i\theta}, \a)
 \widetilde \psi  _0((|\zeta| - t)e ^{i\theta}, \a) dt  \right |$$
$$\le \int_0^{|\zeta|}|\widetilde F (te ^{i\theta}, \a)|.
 |\widetilde \psi  _0((|\zeta| - t)e ^{i\theta}, \a)|dt \le \int_0^{|\zeta|}\widetilde G (t)
 \widetilde h  _0(|\zeta| - t) dt = \widetilde G * \widetilde h  _0(|\zeta|). $$
Therefore, for all $(\zeta, \a) \in \left(\Omega_m(\rho)\backslash \{0\} \right)\times K$,
 $$|\widetilde \psi _1(\zeta, \a)| \le \frac{\widetilde G * \widetilde h
  _0(|\zeta|)}{|\zeta|\rho} = \widetilde h  _1(|\zeta|).$$
This gives (\ref{major2}) for $n=1$ by an argument of continuity.\\
Now, 
$\displaystyle 
{\widetilde h} (\zeta, \varepsilon) = \sum_{n=0}^\infty 
{\widetilde h_n} (\zeta) \varepsilon^n$ is nothing but the minor of
the series expansion $\displaystyle 
{ h} (z, \varepsilon) = \sum_{n=0}^\infty 
{h_n} (z) \varepsilon^n$, where the $h_n$'s are defined recursively by:
\begin{equation}\label{rec4}
\left\{
\begin{array}{l}
\displaystyle - \rho \frac{d}{dz}{ h_0} = {G}\\
\\
\displaystyle  - \rho \frac{d}{dz}{ h}_{n+1} =
{h}_{n}  {G} , \, \, n \geq 0 .\\
\end{array}
\right.
\end{equation}
This means that $h$ satisfies the following ordinary differential
equation:
\begin{equation}\label{forh}
- \rho \frac{d}{dz} h = \varepsilon G(z)h +G(z).
\end{equation}
From (\ref{cestG}), we see that $G$ is integrable at infinity, so that
the function
\begin{equation}\label{solh}
\displaystyle (z, \varepsilon) \mapsto  \frac{e^{-\frac{\varepsilon}{\rho}
    \int_{+\infty}^{z}G(z') dz'} - 1}{\varepsilon}
\end{equation}
is a solution of equation (\ref{forh}) which is holomorphic 
for $z$ in a neighbourhood of infinity of  $\mathfrak{C}_m$
 and
$\varepsilon \in D(0,R)$, $R>1$. Moreover, its Taylor series expansion
at $\varepsilon=0$ is exactly $\displaystyle 
{ h} (z, \varepsilon) = \sum_{n=0}^\infty 
{h_n} (z) \varepsilon^n$.\\
In return, this proves that $\widetilde h (\zeta, \varepsilon)$
defines a holomorphic function in  $(\zeta,\varepsilon) \in 
\mathfrak{C}_m \times D(0,R)$, with an exponential growth of order not
greater than $1$ at infinity in $\zeta$, uniformly in $\varepsilon
\in D(0,R)$: there exist $A, B \in ]0, +\infty[$ such that
$$\forall (\zeta,\varepsilon) \in 
\mathfrak{C}_m \times D(0,R), \, \, |h(\zeta, \varepsilon) | \leq 
Ae^{B |\zeta|}.$$
This last result, together with (\ref{major2}), shows that the series expansion
${\widetilde \psi}_{n}(\zeta, \a, \varepsilon)$ converges uniformly
for $\zeta$ in every compact set of $\Omega_m(\rho)$, $\a \in K$ and
$\varepsilon \in D(0,R)$, and moreover,
$$\forall (\zeta, \a, \varepsilon) \in \Omega_m(\rho) \times K \times
D(0,R), \, \, 
|{\widetilde \psi} (\zeta, \a,  \varepsilon) | 
\leq {\widetilde h} (|\zeta|, |\varepsilon|) \leq Ae^{B |\zeta|}.$$
Putting $\varepsilon =1$, we deduce the same result for 
${\widetilde \psi}_+ (\zeta, \a)$: holomorphy in $\Omega_m(\rho)
\times K$, exponential growth of order not
greater than $1$ at infinity in $\zeta$, uniformly in $\a \in
K$. \\
Since  $\rho>0$ can be chosen arbitrarily small, we have shown that, 
except for the directions of argument  $\alpha = \pi \mod (2\pi)$, 
there is no singular point on the half line 
$\arg \zeta =\alpha$ and,  ${\widetilde \psi}_+ (\zeta, \a)$ having
an exponential growth of order not greater than $1$ at infinity in
$\zeta$, uniformly in $\a \in
K$, we deduce that $\psi_+(z,\a)$  is Borel-resummable with respect
to $z$, uniformly in $\a$ for $\a$ in any compact set of
$\mathbb{C}^m$, for every direction of summation except those of
argument  $\pi \mod (2\pi)$. 

Thanks to   (\ref{cestphiminus}), an analogous result can be obtained
 ${\psi_-} (z, \a)$.  This yields proposition \ref{propborel}.
\end{proof}

Proposition \ref{propborel} is enough to prove theorem
\ref{Sibuya}. Let us define $\phi_0(x,\a) \in 
\mathbb{C}[\a][[x^{-\frac{1}{2}}]]$ by the following formula~:
$$ \displaystyle
x^{r(\a)}e^{-S(x,\a)}\phi_0(x,\a)=
\frac{\sqrt{x}}{P_m(x,\a)^{\frac{1}{4}}}e^{-z}\psi_+(z,\a)\, 
|_{\displaystyle z=z(x,\a)}.$$
Due to the very definition of $\psi_+$, the left-hand side of this
equality is a formal solution of equation $(\mathfrak{E}_m)$. \\
We know from proposition \ref{propborel}
that $\psi_+(z,\a)$ is Borel-resummable for the direction of argument
 $0$. For $\a$ in any given compact set $K$ of $\mathbb{C}^m$, 
this allows us to define the function
$$\Phi_0 (x,\a) = \frac{\sqrt{x}}{P_m(x,\a)^{\frac{1}{4}}}e^{-z}\, 
\mbox{\sc s}_{0} \psi_+(z,\a)\, 
|_{\displaystyle z=z(x,\a)},$$
which is an analytic solution of $(\mathfrak{E}_m)$ for $z$ in a sectorial neighbourhood of
infinity of aperture $]-\frac{\pi}{2} , \frac{\pi}{2}[$ and $\a$ in $K$. 
Note that the size of the sectorial
neighbourhood may depend on $K$. 
 By the inverse map $z \leftrightarrow x$ (given by  
(\ref{inverse1})), this corresponds to a $x$-sectorial neighbourhood of
infinity    of aperture
$]-\frac{\pi}{m} , \frac{\pi}{m}[$.
From proposition 
\ref{propborel} again,  $\Phi_0$ can be analytically
extended by varying the direction of summation on $]-\pi, \pi[$. 
This shows that $\Phi_0$ is holomorphic in a 
$x$-sectorial neighbourhood  of infinity  $\Sigma_0'$  of aperture
$]-\frac{3\pi}{m} , \frac{3\pi}{m}[$ and,  by construction, $\Phi_0$
is asymptotic to $\displaystyle x^{r(\a)}e^{-S(x,\a)}\phi_0(x,\a)$
at infinity in $\Sigma_0'$, uniformly in $\a \in K$. \\
The uniqueness of $\Phi_0$ follows from the Watson theorem
\cite{Mal95}.\\
Also, since  for any
strict sub-sector $\Sigma$ of $\Sigma_0$ the set $\Sigma \backslash
\Sigma \cap \Sigma_0'$ is bounded, all we have to do now to get 
part $(1)$ of theorem \ref{Sibuya} is to show that $\Phi_0$ extends
analytically in $x \in \Sigma_0$. This is a consequence of  the 
Cauchy-Kovalevskaya theorem: take a point $x_0$ in $\Sigma_0'$
and consider the datum $(\Phi_0(x_0,\a), \Phi_0'(x_0,\a))$. Then $\Phi_0$
is uniquely defined by this Cauchy datum, which is holomorphic in $\a
\in K$. Since the linear differential equation  $(\mathfrak{E}_m)$ is
holomorphic in $(x,\a) \in \mathbb{C}^\star \times \mathbb{C}$, we
conclude that $\Phi_0$ extends analytically to
$\widetilde{\mathbb{C}^\star} \times K$, where
$\widetilde{\mathbb{C}^\star}$ stands for the universal covering of 
$\mathbb{C}^\star$. We end by noticing  that $K$ can be chosen
arbitrarily. This also shows part (2) of
theorem \ref{Sibuya}.
\\
Part $(3)$ of  theorem \ref{Sibuya} follows from the fact that the Borel
resummation  with respect to $z$ commutes with the
derivative $\displaystyle \frac{d}{dz}$.

Note that besides proving theorem \ref{Sibuya}, we have obtained the
following interesting result:

\begin{prop}\label{PropBorelPhi0}
When $m$ is odd, 
the analytic function $\Phi_0$ of theorem \ref{Sibuya} is given by
\begin{equation}\label{Phi0Borel1}
\Phi_0 (x,\a) = \frac{\sqrt{x}}{P_m(x,\a)^{\frac{1}{4}}}e^{-z}\, 
\mbox{\sc s}_{\alpha} \psi_+(z,\a)\, 
|_{\displaystyle z=z(x,\a)},
\end{equation}
for $x$ in a sectorial neighbourhood of
infinity    of aperture
$]-\frac{\pi}{m} -\frac{2\alpha}{m},
\frac{\pi}{m}-\frac{2\alpha}{m}[$, uniformly in $\a$ for $\a$ in any
compact set of $\mathbb{C}^m$,
where the direction of Borel resummation $\alpha$ runs through $]-\pi, +\pi[$.
\end{prop}

The arguments used to prove proposition \ref{propborel} can be
extended to analyze the whole analytic structure of the minor 
$ {\widetilde \psi}_+ (\zeta, \a) $ of  ${\psi}_+ (z, \a)$. Since the
techniques involved are the same as those used in \cite{L95} and
\cite{GS001}, we just give the final 
result\footnote{For the particular reader, this part is 
detailed in Rasoamanana \cite{RasoaThesis}}: 

\begin{prop}\label{resurgenceprop}
The minor $ {\widetilde \psi}_+ (\zeta, \a) \in 
\mathbb{C}[\a]\{\zeta^{\frac{1}{m}}\}$ ({\sl resp.}  
$ {\widetilde \psi}_- (\zeta, \a) \in 
\mathbb{C}[\a]\{\zeta^{\frac{1}{m}}\}$)  of $ {\psi}_+ (z, \a)$ ({\sl
  resp.} 
$ {\psi}_- (z, \a)$) can be extended analytically to $(\zeta, \a) \in  
\widetilde{\mathbb{C} \backslash \{0,-2\}} \times \mathbb{C}^m$ 
({\sl  resp.} $(\zeta, \a) \in  
\widetilde{\mathbb{C} \backslash \{0,+2\}} \times \mathbb{C}^m$), where 
$\widetilde{\mathbb{C} \backslash \{0,\pm 2\}}$ is the universal covering
of $\mathbb{C} \backslash \{0,\pm 2\}$. Moreover, $ {\widetilde
  \psi}_\pm$ has an exponential growth of order not greater than 1 at
infinity in $\zeta$, uniformly in $\a$ for $\a$ in any given compact
set of  $\mathbb{C}^m$.
\end{prop}

One can make things more precise concerning the resurgent
structure, that is the behavior of $ {\widetilde \psi}_+$ and $
{\widetilde \psi}_-$ near their singular points. To do that, we shall
use the alien derivations. \\
We would like to
compute  $\Delta_\tau \psi_+$, where $\Delta_\tau$ stands for 
 the alien derivation at $\tau$. From
proposition \ref{resurgenceprop} we know that the singular points
of the minor of $\psi_+$ lie above $-2$ and $0$. However, since 
$\psi_+$ belongs to $\mathbb{C}[\a][[z^{-\frac{1}{m}}]]$, 
the non vanishing  $\Delta_\tau \psi_+$ can
be indexed by the elements $\tau$ above $-2$
 on the Riemann surface $\mathfrak{C}_m$ of
$z^{\frac{1}{m}}$. \\
We now use one of the fundamental properties of the alien derivations:
the pointed alien derivation $\dot \Delta_\tau = e^{-\tau z}
\Delta_\tau$ commutes with $\frac{d}{dz}$ (proposition \ref{commute}).
Using the fact that the resurgent symbol (definition \ref{ressymbol})
 $\Psi_+(z,\a) =
e^{-z} \psi_+(z,\a)$ is solution of equation (\ref{mathfrakF}) and
that $F$ is a constant of resurgence (definition \ref{constres}), we
obtain:
$$ - \frac{d ^2}{d z^2}\left(\dot \Delta_\tau \Psi_+ \right) + 
(1-F(z,\a))\left(\dot \Delta_\tau \Psi_+ \right)=0.
$$
This means that  $\dot \Delta_\tau \Psi_+$ satisfies the same equation 
(\ref{mathfrakF}). Since $(\Psi_+, \Psi_-)$ is a fundamental
 system of formal solutions for (\ref{mathfrakF}), we can conclude
 that 
$\dot \Delta_\tau \Psi_+$ has to be proportional to the resurgent symbol $\Psi_-$ by an
argument of singular support (definition \ref{ressymbol}): the
singular support of $\Psi_+$ (resp. $\Psi_-$) reduces to $\{+1\}$
(resp. $\{-1\}$), whereas,  by definition, the resurgent symbol $\dot \Delta_\tau \Psi_+ =
e^{+z}\times(\mbox{a formal resurgent function})$ has $\{-1\}$ for its
singular support. We deduce that there is
$\delta_\tau(\a)$ such that $\dot \Delta_\tau \Psi_+ (z,\a) = \delta_\tau(\a)
\Psi_- (z,\a)$, i.e., $\Delta_\tau \psi_+ (z,\a) = \delta_\tau(\a)
\psi_- (z,\a)$. 
Similarly, one obtains the existence of $\delta_\tau(\a)$ such that 
 $\Delta_\tau \psi_- (z,\a) = \delta_\tau(\a)
\psi_+ (z,\a)$, where $\tau$ is above $+2$
 on the Riemann surface $\mathfrak{C}_m$.

 The coefficients   $\delta_\tau(\a)$ are entire
functions of $\a$: this stems directly from the regularity in 
$\a$ of the formal
series $\psi_+, \psi_-$, 
and from the fact that the location of the singular points of the minors 
does not depend on $\a$, 
 so that the Stokes automorphism (in any direction) commutes with the analytic
continuation in $\a$. However, this will be a consequence of theorem
\ref{relSkCk} which will be discussed in a moment.

To sum up:

\begin{thm}\label{Summability1}
For $m$ odd, there exists a unique formal power series $\psi_+(z,\a) \in
\mathbb{C}[\a][[z^{-\frac{1}{m}}]]$ ({\sl resp.}  $\psi_-(z,\a) \in
\mathbb{C}[\a][[z^{-\frac{1}{m}}]]$) whose residual coefficient is 1,
such that $e^{- z} \psi_+(z,\a)$ ({\sl resp.} $e^{+ z} \psi_-(z,\a)$) is
solution of equation  (\ref{mathfrakF}), and moreover:
\begin{equation}\label{psimoins}
\psi_-(z, \a) =  
\psi_+(\omega^\frac{m}{2}z, \omega . \a).
\end{equation}
These formal power series $\psi_{\pm}$ are resurgent in $z$
with holomorphic dependence in $\a$, and
Borel-resummable\footnote{Except of course for the singular directions
  which are described by the resurgence structure.} in $z$, 
uniformly with respect to $\a$ in any compact set.\\
Their resurgent structure is given by:
$$ \left\lbrace \begin{array}{l}
\displaystyle \Delta_{2e^{ki\pi}}\psi_-(z, \a)=S_k(\a) \psi_+(z,\a) \quad 
\text{for} \quad   k \in 2\mathbb{Z}\\
\displaystyle \Delta_{2e^{ki\pi}}\psi_+(z,\a)=S_k(\a) \psi_-(z,\a) \quad 
\text{for} \quad   k-1 \in 2\mathbb{Z}\\
\displaystyle \Delta_\tau \psi_{\pm}=0 \quad  \quad \text{otherwise},
\end{array} \right.$$ 
where $\Delta_\tau$ is the alien derivation at $\tau$. The 
coefficients  $S_k(\a)$ are  entire functions in $\a$ and, 
for all $k \in
\mathbb{Z}$, $S_k = S_{k \mod 2m}$.
\end{thm}

\begin{defn}
The coefficients $S_k (\a)$, $k \in \mathbb{Z}$, 
 are called {\sl the Stokes multipliers}.
\end{defn}

\subsubsection{The $m$ even case}

The fundamental difference with the previous $m$ odd case is now the existence
of the term $b_0(\a)\ln(x)$ in the asymptotic expansion of $z(x,\a)$ 
(defined by
(\ref{Green1})) at infinity in $x$. This is
why it is worth considering the following new Green-Liouville transformation 
\begin{equation}\label{Green2}
 \left\lbrace \begin{array}{l}
\displaystyle \widetilde{z}=\widetilde{z}(x,\a)=
\int^x \frac{\sqrt{P_m(t,\a)}-b_0(\a)}{t} dt\\
\\
\displaystyle \Psi(\widetilde{z},\a)=
\frac{\sqrt{\sqrt{P_m(x,\a)} -b_0(\a)}}{\sqrt{x}} \Phi(x,\a)
\end{array} \right. 
\end{equation}
so that the Laurent-Puiseux series expansion of $x \mapsto
\widetilde{z}(x,\a)$
coincides with the map $x \mapsto S(x,\a)$ defined in
theorem \ref{Sibuya} modulo an analytic function 
vanishing at infinity. The quasi-homogeneity properties (\ref{qhz1})
and (\ref{qhx1}) are still valid for the maps $(x,\a) \mapsto
\widetilde{z}(x,\a)$ and  $(\widetilde{z}, \a) \mapsto 
x(\widetilde{z},\a)$ respectively.

Equation $(\mathfrak{E}_m)$ is converted into the prepared equation~: 
\begin{equation}\label{mathfrakFprime}
 - \frac{d ^2}{d \widetilde{z}^2}\Psi+
\left(1+\frac{4b_0(\a)}{m\widetilde{z}}
 -H(\widetilde{z},\a)
\right) \Psi=0 
\end{equation}
with
\begin{equation}
 \left\lbrace 
\begin{array}{l}
\displaystyle H(\widetilde{z},a)= 
1+\frac{4b_0(\a)}{m\widetilde{z}}
- \frac{P_m(x,\a)}{(\sqrt{P_m(x,\a)}-b_0(\a))^2}\\
\displaystyle
 -\Big(\frac{xP_m'(x,\a)+P_m(x,\a)-b_0(\a)
\sqrt{P_m(x,\a)}}{4\sqrt{P_m(x,\a)}(\sqrt{P_m(x,\a)}-b_0(\a))^3}\\
\displaystyle \hspace{30mm} +x^2\big(\frac{1}{4}
 \frac{P_m''(x,\a)}{\sqrt{P_m(x,\a)} (\sqrt{P_m(x,\a)}-b_0(\a))^3} \\
\qquad \displaystyle -\frac{1}{16} 
\frac{(P_m'(x,\a))^2(5\sqrt{P_m(x,\a)}-
2b_0(\a))}{(P_m(x,\a)-b_0(\a)\sqrt{P_m(x,\a)})^3(\sqrt{P_m(x,\a)}-
b_0(\a))}\big)\Big)\, 
|_{\displaystyle \widetilde{z}=\widetilde{z}(x,\a)}
\end{array} \right.
\end{equation}
and 
\begin{equation}
 H(\widetilde{z},\a)=
\frac{m^2-4}{4m^2\widetilde{z}^2}(1+O(\tilde{z}^{-\frac{2}{m}})) \in 
\frac{1}{\widetilde{z}^2}\mathbb{C}[\a]\{\widetilde{z}^{-\frac{2}{m}}\}.
\end{equation}
Furthermore, $H$ satisfies the quasi-homogeneity property (\ref{qhF1}). 

One easily proves the existence of a unique formal
solution 
$$\displaystyle \Psi_+(\widetilde{z},\a) = 
e^{-\widetilde{z}}\psi_+(\widetilde{z},\a)$$ 
of (\ref{mathfrakFprime}) satisfying 
$$\displaystyle  \psi_+(\widetilde{z},\a) = 
\widetilde{z}^{-\frac{2b_0(\a)}{m}} \mu_+ (\widetilde{z},\a)$$
 where
 $\mu_+ \in
\mathbb{C}[\a][[\widetilde{z}^{-\frac{2}{m}}]]$ with  residual
coefficient $1$. By quasi-homogeneity, one deduces the existence of another
formal solution 
$$\displaystyle \Psi_-(\widetilde{z},\a) = 
e^{+\widetilde{z}}\psi_-(\widetilde{z},\a) = 
e^{+\widetilde{z}}\widetilde{z}^{+\frac{2b_0(\a)}{m}} 
\mu_- (\widetilde{z},\a)$$
 of (\ref{mathfrakFprime}) such
that 
$$
\psi_-(\widetilde{z}, \a) =  
\psi_+(\omega^\frac{m}{2}\widetilde{z}, \omega .\a).$$

From now on the analysis is exactly the same as in the $m$ odd
case and yields the following results:

\begin{thm}\label{Summability2}
For $m$ even, there exists a unique formal series
 $\displaystyle  \psi_+(\widetilde{z},\a) = 
\widetilde{z}^{-\frac{2b_0(\a)}{m}} \mu_+ (\widetilde{z},\a)$
where $\mu_+ \in
\mathbb{C}[\a][[\widetilde{z}^{-\frac{2}{m}}]]$ with  residual
coefficient $1$
 ({\sl resp.}  $\displaystyle \psi_-(\widetilde{z},\a) = 
\widetilde{z}^{+\frac{2b_0(\a)}{m}} 
\mu_- (\widetilde{z},\a)$, $\mu_- \in
\mathbb{C}[\a][[\widetilde{z}^{-\frac{2}{m}}]]$ with  residual
coefficient $1$),
such that $e^{- \widetilde{z}} \psi_+(\widetilde{z},\a)$ 
({\sl resp.} $e^{+ \widetilde{z}} \psi_-(\widetilde{z},\a)$) is
solution of equation  (\ref{mathfrakFprime}). Moreover,
\begin{equation}\label{psimoinsprime}
\psi_-(\widetilde{z}, \a) =  
\psi_+(\omega^\frac{m}{2}\widetilde{z}, \omega .\a).
\end{equation}
The formal power series $\psi_{\pm}$ are resurgent in 
$\widetilde{z}$
with holomorphic dependence in $\a$, and Borel-resummable in $\widetilde{z}$, 
uniformly with respect to $\a$ in any compact set.\\
There exists a set of  entire functions 
$S_k(\a)$, the
{\sl  Stokes multipliers},   such that~:
$$ \left\lbrace \begin{array}{l}
\displaystyle \Delta_{2e^{ki\pi}}\psi_-(\widetilde{z},\a)=
S_k(\a) \psi_+(\widetilde{z},\a) \quad 
\text{for} \quad   k \in 2\mathbb{Z}\\
\displaystyle \Delta_{2e^{ki\pi}}\psi_+(\widetilde{z},\a)=
S_k(\a) \psi_-(\widetilde{z},\a) \quad 
\text{for} \quad   k-1 \in 2\mathbb{Z}\\
\displaystyle \Delta_\tau \psi_{\pm}=0 \quad  \quad \text{otherwise},
\end{array} \right.$$ 
where $\Delta_\tau$ is the alien derivation at $\tau$. 
\end{thm}

In this theorem, due to  the fact that the formal solutions  
$\mu_+$ and  $\mu_-$ belong to
$\mathbb{C}[\a][[z^{-\frac{2}{m}}]]$, the alien derivations need only
to be  indexed by elements on the Riemann surface of
$z^{\frac{2}{m}}$. Thus {\sl a priori} only $m$ Stokes multipliers
govern the resurgence structure. Nevertheless, it is better to
describe the resurgence structure in terms of $\psi_+$ and  $\psi_-$,
which have to be thought of as formal functions on the universal covering of
$\mathbb{C}^\star$. \\

Returning to theorem \ref{Sibuya}, we finally get the desired result:

\begin{prop}\label{PropBorelPhi0bis}
When $m$ is even, 
the analytic function $\Phi_0$ of theorem \ref{Sibuya} is given by
\begin{equation}\label{Phi0Borel2}
\Phi_0 (x,\a) = \frac{\sqrt{x}}{(\sqrt{P_m(x,\a)} - b_0(\a))^\frac{1}{2}}
e^{-\widetilde{z}}\, 
\mbox{\sc s}_{\alpha} \psi_+(\widetilde{z},\a)\, 
|_{\displaystyle \widetilde{z}=\widetilde{z}(x,\a)},
\end{equation}
for $x$ in a sectorial neighbourhood of
infinity    of aperture
$]-\frac{\pi}{m} -\frac{2\alpha}{m},
\frac{\pi}{m}-\frac{2\alpha}{m}[$, uniformly in $\a$ for $\a$ in any
compact set of $\mathbb{C}^m$,
where the direction of Borel resummation $\alpha$ runs through $]-\pi, +\pi[$.
\end{prop}

\subsection{Some properties of the Stokes multipliers} 

The quasi homogeneity induces some interesting
properties of the Stokes multipliers.

To simplify, we assume from now on that:

\begin{notation}
\begin{equation}\label{newomega}
\omega=e^{\frac{2i\pi}{m}}.
\end{equation}
\end{notation}

We recall the following easy result in resurgence theory, cf.
\cite{Ec81-1}.

\begin{lem}\label{nonsense}
Let $\psi_1(y)$ be a formal resurgent function and let $\nu$ be a
nonzero complex number. 
Setting $y=\nu t$ and $\psi_2(t)=\psi_1(y)$, we have the following
equality~:
$$ \displaystyle \Delta^t_{\nu \tau} \psi_2 = \Delta^y_{\tau} \psi_1.$$
where $\Delta_\tau^x$ denotes the alien derivation at $\tau$ with
respect to the variable $x$.
\end{lem}

\begin{prop}\label{qhSk}
With the notations of theorem \ref{Summability1} and 
\ref{Summability2}, we have, for all $k \in \mathbb{Z}$:
$$ \displaystyle S_k(\a)=
S_0(\omega^k.\a).$$
\end{prop}

\begin{proof}
In theorems \ref{Summability1} and \ref{Summability2}
we introduce $t= z$ for $m$ odd, $t=\widetilde z$ for $m$ even.
From (\ref{psimoins}) and (\ref{psimoinsprime}), we get
$$ \left \lbrace \begin{array}{l}
\displaystyle \psi_+(e^{+ i\pi}t,\omega.\a)=\psi_-(t,\a) \\
\\
\displaystyle \psi_+(t,\a)=
\psi_-(e^{-i\pi}t, \omega^{-1}.\a).
\end{array} \right. $$
Using lemma \ref{nonsense} with $y=e^{i\pi} t$, we deduce
$$ \displaystyle \Delta^y_{2e^{i0}} \psi_-(y,\omega.\a) = 
 \Delta^t_{2e^{i\pi}} \psi_+(t,\a).$$
Now, by definition of $S_0$ and $S_1$ 
(theorems \ref{Summability1} and \ref{Summability2}), we have the equalities~:
$$ \left \lbrace \begin{array}{l}
\displaystyle \Delta^y_{2e^{i0}} \psi_-(y, \omega.\a)=
S_0(\omega.\a) \psi_+(y,\omega.\a) \\
\\
\displaystyle
\Delta^z_{2e^{i\pi}}\psi_+(t,\a)=S_1(\a)
\psi_-(t,\a).
\end{array} \right. $$
Finally, we obtain~:
$$ \displaystyle S_1(\a)=S_0(\omega.\a).$$
We end the proof by an easy induction argument.
\end{proof}

Since $\omega^m = e ^{2i\pi}$, we get:

\begin{cor}
For all $k \in \mathbb{Z}$, 
$$S_k = S_{k \mod m}.$$
\end{cor}

\subsection{Stokes-Sibuya coefficients and  Stokes
  multipliers}

To describe the  connection formulas, we have now two sets of Stokes
coefficients  at our disposal. One is made up of  the Stokes-Sibuya
coefficients $C_k (\a)$,  the other is made up of the Stokes multipliers
$S_k (\a)$. The following proposition clarifies 
the relations between these two fundamental data.

\begin{thm}\label{relSkCk}
We consider the Stokes-Sibuya  coefficients $C_k$ given by theorem
\ref{SibuyaCoeff} and the Stokes multipliers described by theorems
\ref{Summability1} and \ref{Summability2}. Then, for all $k \in \mathbb{Z}$,
\begin{equation}\label{SkCk}
  S_k(\a)=\omega^{r(\omega^k.\a)}C_k(\a)
\end{equation}
where $\omega$ is given by (\ref{newomega}).
\end{thm}

In this theorem, $r(\a)$ has been defined in theorem
\ref{Sibuya}. In particular, when $m$ is odd, then 
$\displaystyle r(\a)=\frac{1}{2}-\frac{m}{4}$ does not depend on $\a$,
so that (\ref{SkCk}) simply reads
$$\mbox{for } m \mbox{ odd }, 
 S_k(\a)=\omega^{r}C_k(\a).$$

\begin{proof}

To simplify, we give the proof for $m$ odd only, so that 
$\displaystyle r=\frac{1}{2}-\frac{m}{4}$.

By proposition \ref{qhSk} and formula (\ref{qhCk}) of theorem
\ref{SibuyaCoeff}, it is sufficient to show (\ref{SkCk})  for $k=0$. 
By proposition \ref{PropBorelPhi0}, 
 $\Phi_0$ of theorem \ref{Sibuya}  can be defined by 
\begin{equation}\label{firstpoint}
\Phi_0 (x,\a) = \frac{\sqrt{x}}{P_m(x,\a)^{\frac{1}{4}}}e^{-z}\, 
\mbox{\sc s}_{0} \psi_+(z,\a)\, 
|_{\displaystyle z=z(x,\a)},
\end{equation}
for $z$ in a sectorial neighbourhood of
infinity of aperture $]-\frac{\pi}{2} , \frac{\pi}{2}[$, which
corresponds to $x$ in a sectorial neighbourhood of
infinity    of aperture
$]-\frac{\pi}{m} , \frac{\pi}{m}[$.\\
Now, by theorem \ref{SibuyaCoeff},  $\Phi_1$ is defined by 
$$
 \displaystyle \Phi_1(x,\a)=\Phi_0(\omega x, \omega.\a). 
$$
Using (\ref{firstpoint}), we get 
the following representation for $\Phi_1$~:
\begin{equation}\label{secondpoint} 
\Phi_1(x,\a)=\frac{\sqrt{x\omega}}{P_m(x\omega,\omega.
  \a)^{\frac{1}{4}}}e^{-z} \mbox{\sc s}_{0} \psi_+(z,\omega. \a)\, 
|_{\displaystyle z=z(\omega x,\omega.\a)},
\end{equation}
for $\omega x$ in a sectorial neighbourhood of
infinity    of aperture
$]-\frac{\pi}{m} , \frac{\pi}{m}[$, i.e., 
$x$ in a sectorial neighbourhood of
infinity    of aperture
$]-\frac{3\pi}{m} , -\frac{\pi}{m}[$, so that $z$ belongs to  
a sectorial neighbourhood of
infinity    of aperture $]-\frac{3\pi}{2}, -\frac{\pi}{2}[$.
By quasi-homogeneity considerations, we have seen that 
$\displaystyle z(\omega x,\omega .\a)=e^{i\pi}z(x,\a),$
(cf. formula (\ref{qhz1})),
so that,  by formula (\ref{psimoins}) of theorem \ref{Summability1},
$$\displaystyle \psi_+(z(\omega x,\omega .\a), \omega .\a) =\psi_-(z,\a).$$
Also,
$$\displaystyle P_m(\omega x,\omega .\a)=\omega^m P_m(x,\a).$$
This means that (\ref{secondpoint}) can be written as
$$ \displaystyle \Phi_1(x,\a)=\omega^{r}
\frac{\sqrt{x}}{P_m(x,\a)^{\frac{1}{4}}}e^z
 \mbox{\sc s}_{\pi}\psi_-(z,\a)\,  |_{\displaystyle z=z(x,\a)} $$
for $x$ ({\sl resp.} $z$) in a  sectorial neighbourhood of
infinity    of aperture  $]-\frac{3\pi}{m} , -\frac{\pi}{m}[$ 
({\sl resp.} $]-\frac{3\pi}{2} , -\frac{\pi}{2}[$).
As for $\Phi_1$, we have the following representation for $\Phi_{-1}$~:
$$ \displaystyle \Phi_{-1}(x,\a)=\omega^{-r}
\frac{\sqrt{x}}{P_m(x,\a)^{\frac{1}{4}}}e^z
 \mbox{\sc s}_{-\pi}\psi_-(z,\a)\,  |_{\displaystyle z=z(x,\a)}$$
for $x$ ({\sl resp.} $z$) in a  sectorial neighbourhood of
infinity    of aperture  $]\frac{\pi}{m} , \frac{3\pi}{m}[$ 
({\sl resp.} $]\frac{\pi}{2} , \frac{3\pi}{2}[$).

To compare $\Phi_{-1}$, $\Phi_{0}$, and $\Phi_{1}$, we rotate the
directions of resummation so as to sum along the direction $0$. 
 Since $\displaystyle \Delta_{2e^{i0}}
\psi_-(z,\a)=S_0(\a) \psi_+(z,\a)$ (cf. theorem \ref{Summability1}),
$\psi_-$ is not Borel-resummable in the direction $0$ if $S_0(\a) \neq
0$, but only right
or left Borel-resummable. In other words, we
have to take into account a Stokes phenomenon. Since 
${\underline {\dot \Delta}}_0 \psi_-(z,\a)$ (definition \ref{deltal})  reduces to
$\displaystyle {\dot \Delta}_{2e^{i0}} \psi_-(z,\a)$, one gets  
$\displaystyle {\mathfrak S}_0  \psi_-(z,\a) = \psi_-(z,\a) +
{\dot \Delta}_{2e^{i0}} \psi_-(z,\a)$,
where ${\mathfrak S}_0$ is the Stokes automorphism in the direction
$0$ (definition \ref{Stokes}). Therefore, $\displaystyle \mbox{\sc s}_{0-} \psi_-(z,\a) =
\mbox{\sc s}_{0+} \left[ \psi_-(z,\a) + e^{-2z}S_0(\a)
  \psi_+(z,\a)\right]$, where  $\mbox{\sc s}_{0+}$ (resp. $\mbox{\sc s}_{0-}$) is the
right (resp. left) Borel-resummation in the direction $0$.

We obtain:
\begin{equation}\label{thirdpoint}
 \left \lbrace \begin{array}{l}
\displaystyle \Phi_0(x,\a)=
\frac{\sqrt{x}}{P_m(x,\a)^{\frac{1}{4}}}e^{-z} 
\mbox{\sc s}_0\psi_+(z,\a)\,  |_{\displaystyle z=z(x,\a)}\\
\\
\displaystyle \Phi_1(x,\a)=
\omega^{r}\frac{\sqrt{x}}{P_m(x,\a)^{\frac{1}{4}}}e^{+z} 
\mbox{\sc s}_{0+}\psi_-(z,\a)\,  |_{\displaystyle z=z(x,\a)}\\
\\
\displaystyle
\Phi_{-1}(x,\a)=\omega^{-r}\frac{\sqrt{x}}{P_m(x, \a)^{\frac{1}{4}}}
\mbox{\sc s}_{0+}[e^z\psi_-(z,\a)+e^{-z}S_0(\a)\psi_+(z,\a)]\,  
|_{\displaystyle z=z(x,a)}.
\end{array} \right. 
\end{equation}

By theorem \ref{SibuyaCoeff}, we have 
 the connection formula $\Phi_{-1}=C_0(\a) \Phi_0
+\widetilde{C}_0(\a) \Phi_1$; in this equality, 
replacing $\Phi_{-1}$, $\Phi_0$, $\Phi_1$ by the right-hand sides
 of (\ref{thirdpoint}) and equating the coefficients of  $e^{-z}\mbox{\sc
 s}_0\psi_+$ and $e^{+z} \mbox{\sc s}_{0^+}\psi_-$, we finally get:
$$ \left \lbrace \begin{array}{l}
\displaystyle S_0(\a)=\omega^{r} C_0(\a)\\ 
\\
\displaystyle
\widetilde{C}_0(\a)=\omega^{-2r}=\omega^{m-2+2r}.
\end{array} \right. $$
\end{proof}   


\section{Solutions of $(\mathfrak{E}_m)$ in the neighbourhood
  of the origin: Fuchs theory}\label{sec3}

In order to get more information about the Stokes-Sibuya
coefficients $C_k$ (or about the Stokes multipliers  $S_k$,
 since this  is equivalent, by theorem \ref{relSkCk}), we
have to pick up the necessary information coming from the other
singular point of $(\mathfrak{E}_m)$, i.e., the origin.

Since the origin is a regular singular  point of $(\mathfrak{E}_m)$, the 
classical Fuchs theory allows to describe ``canonical'' systems of
solutions of $(\mathfrak{E}_m)$ near the origin (see, e.g.,
\cite{Rein82, W65}).
The characteristic equation is  $s(s-1)-a_m=0$, so that
$\displaystyle
 \frac{1\pm p}{2}$  are the characteristic values, with 
 $\displaystyle 
p = (1+4a_m)^{\frac{1}{2}}$.

\begin{notation}
In what follows,  $\displaystyle 
p = (1+4a_m)^{\frac{1}{2}}$ and  $\displaystyle s(p)=
\frac{1+p}{2}$. \\
We note $\ap := (a_1, \cdots , a_{m-1})$, and for all $\tau \in
\mathbb{C}$,
$$\tau. \ap := (\tau a_1, \cdots ,\tau^{m-1} a_{m-1}).$$
\end{notation}

As is well-known, we have to distinguish between the following three cases~:
$p\notin \mathbb{Z}$, $p \in
\mathbb{Z}^\star$ and $p=0$. Since we have the freedom for choosing
the determination of the square root $(1+4a_m)^{\frac{1}{2}}$, we
can avoid the case where $p \in - \mathbb{N}^\star$ in the following
theorem.

\begin{thm}\label{Fuchs0}
There exist 
two unique linearly independent solutions $f_1$,
$f_2$ of $(\mathfrak{E}_m)$
such that
$$ \left\lbrace \begin{array}{l}
\displaystyle f_1(x, \ap ,p)=
x^{s(p)} g_1(x, \ap ,p) = 
x^{s(p)}\left(1+\sum_{k=1}^\infty A_k(\ap, p)x^k\right)\\
\\
\displaystyle f_2(x, \ap ,p)=
\lambda(\ap, p) f_1(x, \ap ,p) \ln (x) + x^{s(-p)} g_2(x, \ap ,p)
\end{array} \right.$$ 
where $g_1$, $g_2$ are entire functions in 
$x$ and $\ap$, while $\lambda$ is entire  in $\ap$.
Moreover, $g_1$ is meromorphic in $p$ 
with at most simple poles when $-p \in \mathbb{N}^\star $. 
Precisely, for all $k \in
\mathbb{N}^\star$, 
$\displaystyle A_k(\ap, p) \prod_{l=1}^{k}(p+l) 
\in \mathbb{C}[\ap, p]$.
\begin{enumerate}
\item When $p \notin \mathbb{Z}$, then $\lambda(\ap, p)=0$ and 
$g_2(x, \ap ,p)=g_1(x, \ap ,-p)$. 
\item When $p \in \mathbb{N}^\star$, then
$$g_2(x,\ap, p) = \left(1+\sum_{k=1}^\infty B_k(\ap, p)x^k\right)  
\hspace{5mm} \mbox{with} \hspace{5mm} 
B_p =0.
$$
Moreover,  for all $k \in \mathbb{N}^\star$, 
$\displaystyle \lambda(\ap ,p), \, \, B_k (\ap ,p)  \in \mathbb{C}[\ap]$,
and 
$$\lambda(\omega .\ap, p) = \omega^{-p} \lambda (\ap, p).$$
\item When $p=0$, then $\lambda(\ap ,p) =1$ and
$$g_2(x, \ap, p) = \sum_{k=1}^\infty B_k(\ap, p)x^k
$$
with, for all $k \in \mathbb{N}^\star$, 
$\displaystyle  B_k(\ap, p)  \in \mathbb{C}[\ap]$.
\end{enumerate}
\end{thm}

\begin{rem}\label{rem0bis}
When $-p  \in \mathbb{N}^\star$, 
just change $p$ into $-p$ in theorem \ref{Fuchs0}, which corresponds
to choosing the other root for $\displaystyle (1+4a_m)^{\frac{1}{2}}$,
or equivalently, which corresponds to $a_m$ making a
loop around $-1/4$.
\end{rem}

\begin{rem}\label{rem40}
In the special case when
$\ap= 0$, the function $g_1$ is
meromorphic in $p$ 
with at most simple poles when $-p \in m\mathbb{N}^\star$.
\end{rem}

The existence and unicity of $f_1$ and $f_2$ follow from the Fuchs 
theory, and the chosen normalization for $g_1$ and $g_2$.
The properties of the coefficients $A_k$, $B_k$ and $\lambda$ 
can be proved by induction, and this induces 
the analytic properties of $g_1$, $g_2$. The quasi-homogeneity
property of $\lambda$ is a consequence of 
the quasi-homogeneity of equation
$(\mathfrak{E}_m)$. 

The following result can be shown also by induction, see Rasoamanana 
\cite{RasoaThesis}: 

\begin{prop}\label{Calambda}
We consider $p \in  \mathbb{N}^\star$ and we note
 $p= km + r$, $0 \leq r \leq  m-1$ its Euclidian division by $m$. We
 introduce $\displaystyle \epsilon (r) = \left\{
\begin{array}{l} 
1 \mbox{ if } r \neq 0\\
0 \mbox{ if } r = 0
\end{array}
\right.$. Then
$$
\begin{array}{c}
\displaystyle \lambda(\ap, p) = \\
\\
\displaystyle \frac{1}{p} \sum_{l=k+\epsilon (r)}^p 
\renewcommand{\arraystretch}{0.5}
\begin{array}[t]{c}
\displaystyle \sum \\
{\scriptstyle i_1+ \cdots +i_l =p }\\
{\scriptstyle 1 \leq i_j \leq m}
\end{array}
\renewcommand{\arraystretch}{1}
\frac{a_{m-i_1} \cdots a_{m-i_l}}{i_1 (i_1-p) \cdots (i_1 +\cdots
  +i_{l-1}) (i_1 +\cdots +i_{l-1}-p)} 
\end{array}$$
with the convention $a_0 =1$.
\end{prop}

\begin{rem}\label{rem50}
For $p\in \mathbb{N}^\star$,
\begin{itemize}
\item in the special case when
$\ap=0$, then 
$$\lambda (0, p)\, |_{p \neq 0 \mod m} =0$$
while, for $k \in \mathbb{N}^\star$,
$$\lambda (0, p)\, |_{p = km} =
\frac{(-1)^{k+1}}{m^{2k-1}k\Gamma(k)^2}.$$
\item When $m=2$, then
$$
\left\{
\begin{array}{l}
\displaystyle 
\lambda (\ap, p) \, |_{p = 0 \mod 2} =-
\frac{1}{p\Gamma(p)^2}\prod_{k=1}^{p/2} (a_1 +2k-1)(a_1
-2k+1)\\
\\
\displaystyle 
\lambda (\ap, p) \, |_{p = 1 \mod 2} =
\frac{1}{p\Gamma(p)^2}a_1\prod_{k=1}^{(p-1)/2} (a_1 +2k)(a_1
-2k)\\
\end{array}
\right. .
$$
\end{itemize}
\end{rem}

From the uniqueness of $f_1$ and $f_2$ in theorem \ref{Fuchs0}, and
from the quasi-homogeneity of equation $(\mathfrak{E}_m)$, we easily
obtain:

\begin{cor}\label{CorFuchs1}
We consider the fundamental system of solutions $(f_1, f_2)$ of theorem
\ref{Fuchs0}. Then,
\begin{equation}\label{lesMkgg}
\left(\begin{array}{c} 
f_1 \\ 
f_2
\end{array} 
\right) (\omega x,\omega .\ap, p)= \mathfrak{N} (\ap, p)
\left(
\begin{array}{c} 
f_1 \\
f_2 
\end{array} 
\right) (x,\ap, p)
\end{equation}
with
\begin{equation}\label{lesNkgg}
\mathfrak{N} (\ap, p)=
\left(\begin{array}{cc} 
\omega^{s(p)} & 0 \\ 
\frac{2i \pi}{m} \lambda (\ap, p)\omega^{s(-p)}  & \omega^{s(-p)}
\end{array} 
\right).
\end{equation}
 Moreover
\begin{equation}\label{Mono1}
\left(\begin{array}{c} 
f_1 \\ 
f_2
\end{array} 
\right) (\omega^m x, \ap, p)= \mathfrak{M} (\ap, p) 
\left(
\begin{array}{c} 
f_1 \\
f_2 
\end{array} 
\right) (x,a)
\end{equation}
where
\begin{equation}\label{Mono2}
\mathfrak{M} (\ap, p)  = 
\left(\begin{array}{cc} 
e^{2i\pi s(p)} & 0 \\ 
2i\pi \lambda(\ap, p)e^{2i\pi s(-p)}  & e^{2i\pi s(-p)}
\end{array} 
\right)
\end{equation}
is the monodromy matrix at the origin.
\end{cor}


\section{The $0\infty$ connection matrices}\label{sec4}

In section \ref{sec2}, we described a set of fundamental systems of
solutions $(\Phi_{k-1},\Phi_k)$ of $(\mathfrak{E}_m)$, where 
$k \in \mathbb{Z}$. In section \ref{sec3}, we have obtained another 
fundamental system of solutions  $(f_1,f_2)$. To compare these
fundamental systems, we introduce, for all $k \in \mathbb{Z}$:

\begin{equation}\label{lesMk}
\left(\begin{array}{c} 
\Phi_{k-1} \\ 
\Phi_k 
\end{array} 
\right) (x, \a)= M_k (\ap, p)
\left(
\begin{array}{c} 
f_1 \\
f_2 
\end{array} 
\right) (x, \ap, p)
\end{equation}
where the matrices  $M_k(\ap, p)$ are
invertible.

\begin{defn}
The matrices  $M_k(\ap, p)$ are
called the {\sl $0\infty$-connection matrices}.
\end{defn}

We now give some properties of the $M_k$. 
These properties depend essentially on  $p$, as does
the fundamental system $(f_1, f_2)$.

We start with an obvious result.

\begin{prop}\label{easyprop}
For every $k \in \mathbb{Z}$,
\begin{equation}
M_{k+1}(\ap, p ) = M_k(\omega .\ap, p ) \mathfrak{N}(\ap, p)
\end{equation}
where the invertible matrix $\mathfrak{N}(\ap, p)$ is given by
(\ref{lesNkgg}).
\end{prop}

\begin{proof}
By theorem \ref{SibuyaCoeff}, we write 
$$\displaystyle 
\left(\begin{array}{c} 
\Phi_{k} \\ 
\Phi_{k+1} 
\end{array} 
\right) (x,\a)=
\left(\begin{array}{c} 
\Phi_{k-1} \\ 
\Phi_k 
\end{array} 
\right)(\omega x, \omega .\a)  
 = M_k(\omega .\ap, p)
\left(
\begin{array}{c} 
f_1 \\
f_2 
\end{array} 
\right)(\omega x, \omega . \ap, p).$$
 Since, by definition of $\mathfrak{N}$,
$\displaystyle \left(
\begin{array}{c} 
f_1 \\
f_2 
\end{array} 
\right)(\omega x, \omega .\ap, p) = \mathfrak{N}(\ap, p) \left(
\begin{array}{c} 
f_1 \\
f_2 
\end{array} 
\right)(x, \ap, p)$, 
we can conclude because $(f_1, f_2)$ is a fundamental system.
\end{proof}

\begin{thm}\label{0infty0}
a) For every $k \in \mathbb{Z}$,
\begin{equation}\label{det1}
\displaystyle \det M_k(\ap, p) = 
\left\{
\begin{array}{l}
\displaystyle 2(-1)^{k}\frac{
\omega^{(k-1)(1-\frac{m}{2})+r(\omega^{k}.\a)}
}{p} \hspace{2mm} \mbox{for}\hspace{2mm} p \neq 0 \\
\\
\displaystyle 2(-1)^{k-1}
\omega^{(k-1)(1-\frac{m}{2})+r(\omega^{k}.\a)}
\hspace{2mm} \mbox{for}\hspace{2mm} p = 0.
\end{array}
\right.
\end{equation}

b) For every $k \in \mathbb{Z}$, the matrix
$M_k (\ap, p)$ is entire in $\ap$. More precisely,   
\begin{equation}\label{Mk2}
\begin{array}{c}
 \displaystyle M_k (\ap, p) =\\
 \displaystyle\left(
\begin{array}{cc} 
L_k(\ap, p) & \widetilde{L}_k (\ap, p)\\
\omega^{s(p)}L_k(\omega. \ap, p) + \frac{2i\pi}{m}\lambda(\ap, p)
\omega^{s(-p)}\widetilde{L}_k(\omega. \ap, p)  & 
\omega^{s(-p)} \widetilde{L}_k(\omega. \ap, p)
\end{array} 
\right) 
\end{array}
\end{equation} 
where $L_k(\ap, p)$ and  $\widetilde{L}_k (\ap, p)$ are entire in
$\ap$.

c) For every $k \in \mathbb{Z}$,  the matrix
$M_k (\ap, p)$ is holomorphic in $p \notin \mathbb{Z}$, and
$$\forall p \notin \mathbb{Z}, \forall \ap \in \mathbb{C}^{m-1}, 
\, \, \widetilde{L}_k (\ap, p) = 
L_k(\ap, -p).$$
Moreover,  $\widetilde{L}_k$ extends analytically at $p \in
\mathbb{N}^\star$.

d) For every $k \in \mathbb{Z}$,
\begin{equation}\label{Mk20k}
M_k (\ap, p) = M_0 (\omega^k . \ap, p) 
\left(\begin{array}{cc} 
\omega^{ks(p)} & 0 \\ 
\frac{2i \pi k}{m} \lambda (\ap, p)\omega^{k s(-p)}  & \omega^{ks(-p)}
\end{array} 
\right).
\end{equation}
In particular,
\begin{equation}\label{Mk20k1}
M_m (\ap, p) = M_0 (\ap, p) \mathfrak{M} (\ap, p). 
\end{equation}
\end{thm}

\begin{proof} We only detail the proof for $p \notin \mathbb{Z}$.

a) We deduce from (\ref{lesMk}) that
$$W(\Phi_{k-1},  \Phi_{k}) = \det(M_k) W(f_1,  f_2)$$ 
where $W(., .)$ is the Wronskian. 
From lemma \ref{lemme2}, we know that
$$W(\Phi_{k-1},\Phi_{k}) =  2(-1)^{k-1}
\omega^{(k-1)(1-\frac{m}{2})+r(\omega^{k}.\a)},$$ 
while, by theorem \ref{Fuchs0}, taking the limit $x \rightarrow 0$ and
 using the fact that the wronskian is $x$-independent, one easily gets
$$W(f_1,  f_2) = s(-p)-s(p) = -p.$$

b) From (\ref{lesMk}), we have  $\forall k \in \mathbb{Z}$,
$\displaystyle 
\left(\begin{array}{c} 
\Phi_{k-1} \\ 
\Phi_k
\end{array} 
\right) (x,\a)= M_k( \ap, p)
\left(
\begin{array}{c} 
f_1 \\
f_2 
\end{array} 
\right)(x,\ap, p)$
with 
$\displaystyle
M_k(\ap, p)=
\left(
\begin{array}{cc} 
\beta_{k1}(\ap, p) & \beta_{k2}(\ap, p)\\
\beta_{k3}(\ap, p) & \beta_{k4}(\ap, p)
\end{array} 
\right) 
$
so that, in particular,
\begin{equation}\label{pourphi1}
\Phi_k(x,\a) = \beta_{k3}(\ap, p) f_1(x,\ap, p) + 
\beta_{k4}(\ap, p)f_2(x,\ap, p).
\end{equation}
Then
$$
\left(\begin{array}{c} 
\Phi_{k-1} \\ 
\Phi_k 
\end{array} 
\right) (\omega x,\omega.\a)= M_k(\omega.\ap, p)
\left(
\begin{array}{c} 
f_1 \\
f_2 
\end{array} 
\right)(\omega x,\omega. \ap, p).
$$
 By proposition  \ref{easyprop}
 and  corollary \ref{CorFuchs1}, we get:  
$$
\left(\begin{array}{c} 
\Phi_{k} \\ 
\Phi_{k+1}
\end{array} 
\right) ( x, \a)= M_k(\omega. \ap, p)
\left(
\begin{array}{cc} 
\omega^{s(p)} & 0 \\
0 & \omega^{s(-p)}
\end{array} 
\right)
\left(
\begin{array}{c} 
f_1 \\
f_2 
\end{array} 
\right)(x,\ap, p)
$$
so that
\begin{equation}\label{pourphi1encore}
\Phi_k(x,\a) = \omega^{s(p)} \beta_{k1}(\omega .\ap, p) f_1(x,\ap, p) + 
 \omega^{s(-p)} \beta_{k2}(\omega .\ap, p)f_2(x,\ap, p).
\end{equation}
Comparing (\ref{pourphi1}) and (\ref{pourphi1encore}), we obtain the
announced form for $M_k$ with 
$\beta_{k1} = L_k$ and $\beta_{k2} = \widetilde L_k$,
since $(f_1, f_2)$ is a fundamental system. 

b) and c) We have 
$$ M_k = -\frac{1}{p}
\left(\begin{array}{cc} 
\Phi_{k-1} & \Phi_{k-1}' \\ 
\Phi_k & \Phi_k' 
\end{array} 
\right)
\left(\begin{array}{cc} 
f_2' & -f_1' \\ 
-f_2 & f_1 
\end{array} 
\right)
$$
so that the analytic properties of $M_k$ 
easily follow from the analytic properties  of the
$\Phi_k$'s (theorem \ref{SibuyaCoeff}) and of $f_1, f_2$ 
(theorem \ref{Fuchs0}).

d) The given statement follows from 
 proposition \ref{easyprop}, by induction, inferring from (\ref{lesNkgg})
 that $\mathfrak{N} (\omega. \ap, p) =
\mathfrak{N} ( \ap, p) $.
\end{proof}

In addition to theorem \ref{0infty0}, it is easy to show the following
proposition  (the special case where $\ap = 0$ follows from  remark
\ref{rem40}): 

\begin{prop}\label{rem41}
The restriction to $p \notin \mathbb{Z}$ of the function
$L_k(\ap, p)$ ({\sl resp.} $\widetilde L_k(\ap, p)$) 
 has a meromorphic continuation  in $p$, with at
most simple poles when $p \in \mathbb{N}$ ({\sl resp.} $-p \in \mathbb{N}$).\\
In the special case where
$\ap = 0$,  the
restriction to $p \notin \mathbb{Z}$ of the function
$L_k(\ap, p)$ 
({\sl resp.} $\widetilde{L}_k (\ap, p)$)  has a meromorphic continuation 
 in $p$, with at most simple
poles at $p \in m\mathbb{N}$ ({\sl resp.} $-p \in m\mathbb{N}$). 
\end{prop}


\section{Monodromy, Stokes-Sibuya and  
$0\infty$ connection matrices}\label{sec5}

We collect here the different results we have got on the
monodromy, the Stokes-Sibuya and the $0\infty$ connection matrices to
obtain a set of functional relations.

\subsection{First functional equation}

From the very definition (\ref{lesMk}) of the $0\infty$
connection matrices $M_k$ and the fundamental property (\ref{lesCk}) 
of the Stokes-Sibuya connection matrices, we have for all $k \in \mathbb{Z}$:
$$\displaystyle 
M_k \left(
\begin{array}{c} 
f_1 \\
f_2 
\end{array} 
\right) = 
\left(\begin{array}{c} 
\Phi_{k-1} \\ 
\Phi_{k} 
\end{array} 
\right)= 
\mathfrak{S}_k 
\left(\begin{array}{c} 
\Phi_{k} \\ 
\Phi_{k+1} 
\end{array} 
\right) 
 = \mathfrak{S}_k  M_{k+1} 
\left(
\begin{array}{c} 
f_1 \\
f_2 
\end{array} 
\right).$$ 
Since $(f_1, f_2)$ is a fundamental
system, we  thus have the following proposition:

\begin{prop}\label{propSetM}
For all $k \in \mathbb{Z}$,
\begin{equation}\label{SetM}
\displaystyle \mathfrak{S}_k (\a) = M_k (\ap, p)
  M_{k+1}^{-1}(\ap, p).
\end{equation}
\end{prop}

Using (\ref{SetM}), we see that
$$
\mathfrak{S}_0 (\a) \mathfrak{S}_1 (\a) \cdots 
\mathfrak{S}_{m-1} (\a) = M_0 (\ap, p)  M_{m}^{-1}(\ap, p).$$
Using (\ref{Mk20k1}), we obtain the following theorem:

\begin{thm}\label{FUNCT}
The Stokes-Sibuya connection matrices satisfy the following functional
relation: 
\begin{equation}\label{functional}
\displaystyle \mathfrak{S}_0 (\a) \mathfrak{S}_1 (\a) \cdots 
\mathfrak{S}_{m-1} (\a) = 
M_0 (\ap, p) \mathfrak{M}^{-1}(\ap, p) M_{0}^{-1}(\ap, p).
\end{equation}
\end{thm}

This functional relation is equivalent to formula
(\ref{Monoinfty3}) of theorem \ref{Monoinfty}. But this new
formulation is interesting  
thanks to the following two corollaries. 

\begin{cor}\label{CORFUNCT}
We have
$$Tr \left( \mathfrak{S}_0 (\a) \mathfrak{S}_1 (\a) \cdots 
\mathfrak{S}_{m-1} (\a) \right) =  
-2\cos(\pi p)$$
where $Tr$ is the Trace.
\end{cor}

\begin{proof}
This follows from the fact that 
$\displaystyle  Tr \left( M_0 (\ap, p) \mathfrak{M}^{-1}(\ap, p)
  M_{0}^{-1}(\ap, p) \right) = Tr \left( \mathfrak{M}^{-1}(\ap, p)
\right) = -2\cos(\pi p)$.
\end{proof}

We have also the following result:

\begin{cor}\label{CORFUNCT2}
We assume that $\displaystyle p \in
\mathbb{N}^\star$. Then, with the notations of theorem  
\ref{Fuchs0},
$$
\mathfrak{S}_0 (\a) \mathfrak{S}_1 (\a) \cdots 
\mathfrak{S}_{m-1} (\a) \,|_{\lambda(\ap, p)=0} = 
(-1)^{p+1} \left(\begin{array}{cc} 
1 & 0 \\ 
0 & 1
\end{array} 
\right).
$$
\end{cor}

\begin{proof}
From corollary \ref{CorFuchs1}, we know that $\displaystyle 
\mathfrak{M} (\ap, p)  = e^{2i\pi s(p)}
\left(\begin{array}{cc} 
1 & 0 \\ 
2i \pi \lambda (\ap, p) & 1
\end{array} 
\right)$ with $2s(p)=1+p$,  so that   
$$\begin{array}{c}
\displaystyle  M_0 (\ap, p) \mathfrak{M}^{-1}(\ap, p)
  M_{0}^{-1}(\ap, p)\,|_{\lambda(\ap, p)=0} \\
\\
\displaystyle  =\mathfrak{M}^{-1}(\ap, p)\,|_{\lambda(\ap, p)=0} = 
(-1)^{p+1} \left(\begin{array}{cc} 
1 & 0 \\ 
0 & 1
\end{array} 
\right).
\end{array}$$
\end{proof}

\subsection{Second functional equation}

\begin{thm}\label{genial1}
We use the notations of theorem \ref{0infty0}.
\begin{itemize}
\item 1. We assume  $p \notin \mathbb{Z}$. 
We assume furthermore  that $\a$ is chosen so that,
 for all $ k =0, \cdots, m-1$,
$\widetilde{L}_0 (\omega^k. \ap, p) \neq 0$\footnote{Note that $\widetilde{L}_0
  (\ap, p)$ cannot be identically zero; therefore, this is a generic
  hypothesis on $\a$.}. Then
\begin{equation}\label{genial2}
\frac{L_0(\ap, p)}{\widetilde{L}_0(\ap, p)} = -i
\frac{\omega^{-\frac{3}{2}}\omega^{-(m+1)\frac{p}{2}}}{p \sin(\pi p) } 
\sum_{k=0}^{m-1}  \frac{\omega^{r(\omega^k .\a) +(k+1)p}}
{ \widetilde{L}_0(\omega^{k}.\ap, p) \widetilde{L}_0 (\omega^{k+1}.\ap, p)}.
\end{equation}
\item 2.  We assume  $p \in \mathbb{N}^\star$. 
Assuming also that $\a$ is chosen so that,
 for all $ k =0, \cdots, m-1$,
$\widetilde{L}_0 (\omega^k. \ap, p) \neq 0$, then
\begin{equation}\label{genial2bis}
i\pi p  \omega^{\frac{3}{2} + \frac{p}{2}}  \lambda (\ap, p) =   
 \sum_{k=0}^{m-1} 
\frac{\omega^{r(\omega^k .\a, p) + (k+1)p}}
{ \widetilde{L}_0(\omega^k .\ap, p) \widetilde{L}_0 (\omega^{k+1} .\ap, p)}.
\end{equation}
\end{itemize}
\end{thm}

\begin{proof}
\begin{itemize}
\item 1.
 Using formulas
(\ref{det1}) and (\ref{Mk2}) with $k=0$, we get
$$\displaystyle\omega^{s(-p)} L_0(\ap, p)\widetilde{L}_0(\omega.\ap, p) - 
\omega^{s(p)}L_0(\omega.\ap, p)\widetilde{L}_0 (\ap, p) = 
 -\frac{2}{p}\omega^{-1+r(\a)}
$$
and, more generally, for all $ k =0, \cdots, m-1$,
$$\begin{array}{c} 
\displaystyle \omega^{s(-p)} L_0(\omega^k. \ap, p)\widetilde{L}_0(\omega^{k+1}.\ap, p) - 
\omega^{s(p)}L_0(\omega^{k+1}.\ap, p)\widetilde{L}_0 (\omega^k. \ap,
p) = \\
\\
\displaystyle -\frac{2}{p}\omega^{-1+r(\omega^k. \a)}.
\end{array}
$$
We assume  $\a$ generic so that  for all $ k =0, \cdots, m-1$,
$\widetilde{L}_0 (\omega^k. \ap, p) \neq 0$. 
The previous equalities read also:
$$\displaystyle\omega^{s(-p)} 
\frac{L_0(\omega^k. \ap, p)}{\widetilde{L}_0(\omega^{k}.\ap, p)} - 
\omega^{s(p)} \frac{L_0(\omega^{k+1}.\ap, p)}{\widetilde{L}_0
  (\omega^{k+1}.\ap, p)} = 
- \frac{2}{p} 
\frac{\omega^{-1+r(\omega^k. \a)}}{ \widetilde{L}_0(\omega^{k}.\ap, p) \widetilde{L}_0 (\omega^{k+1}. \ap, p)}.
$$
From the holomorphy in $\ap$ of $L_0$ and $\widetilde{L}_0$,
this can be written in the following form, since $\omega^m =e^{2i\pi}$:
$$
\displaystyle \mathfrak{L}
\left(
\begin{array}{c}
\displaystyle \frac{L_0(\ap, p)}{\widetilde{L}_0(\ap, p)}\\
\vdots\\
\displaystyle  \frac{L_0(\omega^{m-1}. \ap, p)}{\widetilde{L}_0(\omega^{m-1}.\ap, p)}
\end{array}
\right)
= 
\displaystyle - \frac{2\omega^{-1}}{p} 
\left(
\begin{array}{c}
\displaystyle \frac{\omega^{r( \a)}}
{ \widetilde{L}_0(\ap, p) \widetilde{L}_0 (\omega. \ap, p)}\\
\vdots\\
\displaystyle \frac{\omega^{r(\omega^{m-1}. \a)}}
{ \widetilde{L}_0(\omega^{m-1}.\ap, p) \widetilde{L}_0 ( \ap, p)}
\end{array}
\right) ,
$$
 where 
$$ 
\mathfrak{L} = 
\left(
\begin{array}{ccccc}
\omega^{s(-p)} & -\omega^{s(p)} &   0   &  \cdots      &  0\\
      0      & \ddots        & \ddots&  \ddots      & \vdots \\
\vdots       &  \ddots       &\ddots &  \ddots      & 0\\
     0       &  \cdots       &   0   &\omega^{s(-p)}  & -\omega^{s(p)}\\
-\omega^{s(p)}&      0        &\cdots &      0       &\omega^{s(-p)}
\end{array}
\right)
$$
is a $m\times m$ circulant matrix  whose determinant is $\omega^{ms(-p)} -
\omega^{ms(p)}$.  This determinant does
not vanish because $s(p)-s(-p) =p$ is not an integer. The inverse of this
matrix is also a circulant matrix, precisely:
 $$\begin{array}{c}
\displaystyle \mathfrak{L}^{-1} = \frac{1}{\omega^{ms(-p)}
  -\omega^{ms(p)}} \times\\
\\
\displaystyle \left(
\begin{array}{ccccc}
\omega^{(m-1)s(-p)} & \omega^{(m-2)s(-p)+s(p)} &\cdots &
\omega^{s(-p)+(m-2)s(p)}&  \omega^{(m-1)s(p)} \\
\omega^{(m-1)s(p)} & \omega^{(m-1)s(-p)} &  &\cdots &
\omega^{s(-p)+(m-2)s(p)}  \\
      & \cdots        & \cdots  &  \cdots      &  \\
      & \cdots        & \cdots  &  \cdots      &  \\
\omega^{(m-2)s(-p)+s(p)} &\cdots & &  \omega^{(m-1)s(p)} &\omega^{(m-1)s(-p)}
\end{array}
\right) .
\end{array}
$$
This yields, since $s(p)-s(-p)=p$,
$$\frac{L_0(\ap, p)}{\widetilde{L}_0(\ap, p)} =-
\frac{2\omega^{-1}\omega^{(m-1)s(-p)}}{p(\omega^{ms(-p)} -\omega^{ms(p)}) } 
\sum_{l=0}^{m-1} \omega^{lp}  \frac{\omega^{r(\omega^l. \a)}}
{ \widetilde{L}_0(\omega^{l}.\ap, p) \widetilde{L}_0 (\omega^{l+1}.\ap, p)}
$$
that is also
$$\frac{L_0(\ap, p)}{\widetilde{L}_0(\ap, p)} = i
\frac{\omega^{-1}\omega^{(m-1)s(-p)}}{p \sin(\pi p) } 
\sum_{l=0}^{m-1} \omega^{lp}  \frac{\omega^{r(\omega^l. \a)}}
{ \widetilde{L}_0(\omega^{l}.\ap, p) \widetilde{L}_0 (\omega^{l+1}.\ap, p)}.
$$

\item 2. We work 
with formulas
(\ref{det1}) and (\ref{Mk2}) with $k=0$, when $p \in
\mathbb{N}^\star$. 
Using also the
 fact that $\lambda(\omega .\ap, p) = \omega^{-p} \lambda (\ap, p)$ 
(see  theorem  \ref{Fuchs0}), we get
$$
\begin{array}{c}
\displaystyle 
\mathfrak{L}
\left(
\begin{array}{c}
\displaystyle \frac{L_0(\ap, p)}{\widetilde{L}_0(\ap, p)}\\
\vdots\\
\displaystyle  
\frac{L_0(\omega^{m-1}. \ap, p)}{\widetilde{L}_0(\omega^{m-1}.\ap, p)}
\end{array}
\right)
= \\
\displaystyle  - \frac{2\omega^{-1}}{p} 
\left(
\begin{array}{c}
\displaystyle \frac{\omega^{r( \a)}}
{ \widetilde{L}_0(\ap, p) \widetilde{L}_0 (\omega. \ap, p)}\\
\vdots\\
\displaystyle \frac{\omega^{r(\omega^{m-1}. \a)}}
{ \widetilde{L}_0(\omega^{m-1}. \ap, p) \widetilde{L}_0 ( \ap, p)}
\end{array}
\right) 
+
\frac{2i\pi}{m}  \omega^{s(-p)}\lambda (\ap, p) \left(
\begin{array}{c}
1\\
\\
\vdots\\
\\
\omega^{-(m-1)p}
\end{array}
\right) 
\end{array}
$$
 where 
$\mathfrak{L}$ is the previous circulant matrix. But now $\det
 (\mathfrak{L})= \omega^{ms(-p)} - \omega^{ms(p)}=0$, since 
$s(p)-s(-p) = p \in \mathbb{N}^\star$. It is straightforward to 
see that $\mathfrak{L}$ is  of rank $m-1$, so that the compatibilty
 condition reads
$$ \det
\left(
\begin{array}{cccccc}
\omega^{s(-p)} & -\omega^{s(p)} &   0   &  \cdots      &  0   & \alpha_0 \\
      0      & \ddots        & \ddots&  \ddots      & \vdots &  \\
\vdots       &  \ddots       &\ddots &  \ddots      & 0 & \vdots\\
     0      &  \cdots       &   0   &\omega^{s(-p)}& -\omega^{s(p)}& \\
 0 &      0        &\cdots &      0   &\omega^{s(-p)} &\alpha_{m-2}\\
-\omega^{s(p)}&      0        &\cdots &      0   & 0 &\alpha_{m-1}\\
\end{array}
\right) =0 
$$
where 
$$ \alpha_k = - \frac{2\omega^{-1}}{p} \frac{\omega^{r(\omega^k .\a)}}
{ \widetilde{L}_0(\omega^k .\ap, p) \widetilde{L}_0 (\omega^{k+1} .\ap, p)}
+\frac{2i\pi}{m}  \omega^{s(-p)}\omega^{-kp}\lambda (\ap, p).
$$
This means that
$$
\begin{array}{c}
\omega^{(m-1)s(-p)}\alpha_{m-1} + \omega^{(m-1)s(p)} \alpha_{m-2} + 
\omega^{(m-2)s(p) +s(-p)} \alpha_{m-3} + \\
\cdots +
\omega^{s(p) + (m-2)s(-p)} \alpha_{0}=0,
\end{array}
$$
that is also, because $s(p)-s(-p)=p$,
$$
\alpha_{m-1} + \omega^{(m-1)p} \alpha_{m-2} + 
\omega^{(m-2)p} \alpha_{m-3} + \cdots +
\omega^{p} \alpha_{0}=0.
$$
We eventually get, since $p$ is an integer:
$$
i\pi p \omega^{\frac{3}{2} + \frac{p}{2}}  \lambda (\ap, p) =   
\frac{\omega^{r(\omega^{m-1} .\a)}}
{ \widetilde{L}_0(\omega^{m-1} .\ap, p) \widetilde{L}_0 (\ap, p)} 
+ \sum_{k=0}^{m-2} 
\frac{\omega^{r(\omega^k .\ap, p) + (k+1)p}}
{ \widetilde{L}_0(\omega^k .\ap, p) \widetilde{L}_0 (\omega^{k+1} \ap, p)}.
$$
\end{itemize}
\end{proof}

Theorem \ref{genial1} induces the following interesting result.

\begin{cor}\label{cestbeau}
The Stokes-Sibuya multiplier $C_0(\a)$ satisfies:
\begin{itemize}
\item when $m=1$, for all $\a \in \mathbb{C}$:
$$C_0 (\a) =-2\cos(\pi p),$$
\item when $m=2$, for all $\ap \in \mathbb{C}$ and  $p \notin -\mathbb{N}$:
$$
C_0 (\a) \widetilde{L}_0 (\omega.\ap, p) = -2i
e^{-i\pi\frac{a_1}{2}}\cos\left(\frac{\pi}{2}(p+a_1)\right) 
\widetilde{L}_0(\ap, p),
$$
\item when $m \geq 3$, for all $\ap \in \mathbb{C}^{m-1}$ 
and  $p \notin -\mathbb{N}$:
$$
\begin{array}{c}
\displaystyle C_0 (\a) \widetilde{L}_0(\omega .\ap, p) =\\
\\
\displaystyle  \omega^{r(\a) -1 +\frac{m}{4}}
  \left( \widetilde{L}_0( \ap, p)
\omega^{-r(\a)+\frac{1}{2} -\frac{m}{4}+ \frac{p}{2}} +
\widetilde{L}_0( \omega^2 .\ap, p)
\omega^{r(\a)-\frac{1}{2} +\frac{m}{4}- \frac{p}{2}} 
 \right).
\end{array}
$$
\end{itemize}
\end{cor}

\begin{proof}
Formula (\ref{SetM}) of proposition \ref{propSetM} with
$k=0$ yields 
$\displaystyle \mathfrak{S}_0 (\a) = M_0 (\ap, p)
  M_{1}^{-1}(\ap, p)$. Using (\ref{Mk2}) and (\ref{Mk20k}), we obtain:
$$
\begin{array}{c}
\displaystyle \left(
\begin{array}{cc} 
C_0(\a) & \widetilde{C}_0 (\a) \\ 
1 & 0 
\end{array} 
\right) = \\
\\
\displaystyle \left(
\begin{array}{cc} 
L_0(\ap, p) & \widetilde{L}_0 (\ap, p)\\
\omega^{s(p)}L_0(\omega. \ap, p)   & 
\omega^{s(-p)} \widetilde{L}_0(\omega. \ap, p)
\end{array} 
\right) 
\left(
\begin{array}{cc} 
\omega^{s(p)}L_0(\omega. \ap, p)   & 
\omega^{s(-p)} \widetilde{L}_0(\omega. \ap, p)\\
\omega^{2s(p)}L_0(\omega^2. \ap, p)   & 
\omega^{2s(-p)} \widetilde{L}_0(\omega^2. \ap, p)
\end{array} 
\right)^{-1}
\end{array}
$$
so that, with (\ref{det1}):
\begin{equation}\label{CC0}
C_0 (\a) = - \frac{p}{2} \omega^{1-r(\omega .\a)} \left( 
\omega^{-p} \frac{L_0(\ap, p)}{\widetilde{L}_0(\ap, p)}  - 
\omega^{p} \frac{L_0(\omega^2 .\ap, p)}{ \widetilde{L}_0(\omega^2 .\ap, p)}
\right) \widetilde{L}_0(\ap, p) \widetilde{L}_0( \omega^2 .\ap, p). 
\end{equation}

We now apply formula (\ref{genial2}) under the assumptions made in  theorem
  \ref{genial1}.
 
\begin{itemize}
\item When $m=1$, $\omega = e^{2i\pi}$ and
$\displaystyle r (\a)  = \frac{1}{4}$, and
formula (\ref{genial2}) of  theorem
  \ref{genial1}  reduces to:
$$
\frac{L_0(p)}{\widetilde{L}_0(p)} =  -
\frac{1}{p \sin(\pi p) \widetilde{L}_0( p) \widetilde{L}_0 ( p)}.
$$
This allows to write (\ref{CC0}) as:
$$C_0 (\a) =- \frac{\sin(2\pi p)}{\sin(\pi p)} = -2\cos(\pi p).$$
This result extends for all $\a \in \mathbb{C}$ by analytic
continuation, since $C_0$ is entire in $\a$. 
\item When $m=2$, we have  $\omega = e^{i\pi}$ and
$\displaystyle r (\a) =  -\frac{a_1}{2}$. \\
Formula
(\ref{genial2}) of  theorem \ref{genial1} becomes:
$$
\frac{L_0(\ap, p)}{\widetilde{L}_0(\ap, p)} = 
\frac{\omega^{-\frac{3p}{2}}}{p \sin(\pi p) } \left(
\frac{\omega^{-\frac{a_1}{2}+p}}
{ \widetilde{L}_0(\ap, p) \widetilde{L}_0 (\omega.\ap, p)}
+  \frac{\omega^{\frac{a_1}{2}+2p}}
{ \widetilde{L}_0(\omega.\ap, p) \widetilde{L}_0 (\ap, p)} 
\right),
$$
or also
$$
\frac{L_0(\ap, p)}{\widetilde{L}_0(\ap, p)} = 
2\frac{\cos\left(\frac{\pi}{2}(p+a_1)\right)}{p \sin(\pi p) } 
\frac{1}
{ \widetilde{L}_0(\ap, p) \widetilde{L}_0 (\omega.\ap, p)}.
$$
This means that  (\ref{CC0}) reads:
$$
C_0 (\a) = -2i
e^{-i\pi\frac{a_1}{2}}\cos\left(\frac{\pi}{2}(p+a_1)\right) \frac{\widetilde{L}_0(\ap, p)}{\widetilde{L}_0 (\omega.\ap, p)}.
$$
The announced result follows  by analytic
continuation for all $\ap \in \mathbb{C}$ and all $p \notin
- \mathbb{N}$, since $C_0$ is entire in $\a$, while $\widetilde{L}_0$
  is holomorphic in $\ap \in \mathbb{C}$ and $p \notin -\mathbb{N}$.
\item When $m>2$, we can write by theorem \ref{genial1}
$$
\begin{array}{c}
\displaystyle C_0 (\a) =
i\frac{\omega^{-\frac{1}{2}}\omega^{-(m-1)\frac{p}{2}-r(\omega .\a)
  }}{2 \sin(\pi p) }\widetilde{L}_0(\ap, p) \widetilde{L}_0( \omega^2
.\ap, p)
 \times \\
\\
\displaystyle \left(
\omega^{-p}\sum_{l=0}^{m-1} \frac{\omega^{lp +r(\omega^l. \a)}}
{ \widetilde{L}_0(\omega^{l}.\ap, p) \widetilde{L}_0 (\omega^{l+1}.\ap, p)}
- \omega^{p} \sum_{l=0}^{m-1} \frac{\omega^{lp +r(\omega^{l+2}. \a)}}
{ \widetilde{L}_0(\omega^{l+2}.\a) \widetilde{L}_0 (\omega^{l+3}.\ap, p)}
\right),
\end{array}
$$
which reads also:
$$
\begin{array}{c}
\displaystyle 
C_0 (\a) =\\
\\
\displaystyle
i\frac{\omega^{-\frac{1}{2}}\omega^{-(m-1)\frac{p}{2}-r(\omega .\a)
  }}{2 \sin(\pi p) } \left( 
\frac{\omega^{-p +r(\a)}}
{ \widetilde{L}_0(\ap, p) \widetilde{L}_0 (\omega .\ap, p)} + 
\frac{\omega^{r(\omega. \a)}}
{ \widetilde{L}_0(\omega .\ap, p) \widetilde{L}_0 (\omega^{2}.\ap, p)}
\right. + \\
\\
\displaystyle
\sum_{l=2}^{m-1} \frac{\omega^{(l-1)p +r(\omega^l. \a)}}
{ \widetilde{L}_0(\omega^{l}.\ap, p) \widetilde{L}_0 (\omega^{l+1}.\ap, p)}
-  \sum_{l=0}^{m-3} \frac{\omega^{(l+1)p +r(\omega^{l+2}. \a)}}
{ \widetilde{L}_0(\omega^{l+2}.\ap, p) \widetilde{L}_0
  (\omega^{l+3}.\ap, p)}\\
\\
\displaystyle 
\left. - \frac{\omega^{(m-1)p +r( \a)}}
{ \widetilde{L}_0(\ap, p) \widetilde{L}_0 (\omega. \ap, p)}
- \frac{\omega^{mp +r(\omega. \a)}}
{ \widetilde{L}_0(\omega. \ap, p) \widetilde{L}_0 (\omega^{2}. \ap, p)}
\right)\widetilde{L}_0(\ap, p) \widetilde{L}_0( \omega^2 .\ap, p).
\end{array}
$$
The right-hand side of this equality simplifies to give
$$
\begin{array}{c}
\displaystyle 
C_0 (\a) =\\
\\
\displaystyle  \omega^{r(\a) -1 +\frac{m}{4}}
  \left( 
\frac{\widetilde{L}_0( \ap, p)}{\widetilde{L}_0(\omega .\ap, p)}
\omega^{-r(\a)+\frac{1}{2} -\frac{m}{4}+ \frac{p}{2}} +
\frac{\widetilde{L}_0( \omega^2 .\ap, p)}{\widetilde{L}_0(\omega .\ap, p)}
\omega^{r(\a)-\frac{1}{2} +\frac{m}{4}- \frac{p}{2}} 
 \right) .
\end{array}
$$
Again, the announced result follows  by analytic
continuation.
\end{itemize}
\end{proof}


\subsection{Third functional equation}

In this subsection, we  study a class of differential equations
($\mathfrak{E}_m$) with higher symmetries. For that purpose, it will 
 be useful to introduce new notations. 

\begin{notation}
For  $m, n \in \mathbb{N}^\star$, we define
$$\an = (0, \cdots, 0, a_n, 0, \cdots, 0, a_{2n}, 0, \cdots, 0,a_{jn},
0\cdots, 0, a_{nm}) \in \mathbb{C}^{nm},$$ 
i.e., 
$$\an = (a_j)_{1\leq j \leq nm} \hspace{3mm} \mbox{so that}
\hspace{3mm} a_j = 0 \hspace{3mm} \mbox{if} \hspace{3mm} j \neq 0 \mod m.
$$ 
For such a $\an$, we also define:
$$\anp := (a_j)_{1\leq j \leq nm-1}$$
and 
$$
\left\{
\begin{array}{l}
\displaystyle \ant := \left(\frac{a_{n}}{n^{\frac{2}{m}}}, 
\frac{a_{2n}}{n^{\frac{4}{m}} }, 
\cdots ,  \frac{a_{n(m-1)}}{ n^{\frac{2(m-1)}{m}}}, -\frac{1}{4} +
\frac{1+4a_{nm}}{4n^2} \right) \in \mathbb{C}^m\\
\\
\displaystyle \antp  := \left(\frac{a_{n}}{n^{\frac{2}{m}}}, 
\frac{a_{2n}}{n^{\frac{4}{m}} }, 
\cdots ,  \frac{a_{n(m-1)}}{ n^{\frac{2(m-1)}{m}}} \right) \in
\mathbb{C}^{m-1}.
\end{array}
\right.
$$
\end{notation}

We shall consider in this subsection the following differential
equation:
$$
x^2 \frac{d^2}{d x^2}\Phi(x, \an)= P_{nm}(x, \an)\Phi(x, \an).
\leqno{(\mathfrak{E}_{nm}^n)}
$$
This equation is a particular case of our main equation
$(\mathfrak{E}_{nm})$, but its higher symmetry will allow us 
to compare its Stokes-Sibuya and  $0\infty$ 
connection matrices with those of $(\mathfrak{E}_{m})$, associated
with the polynomial $P_m(x, \ant)$ of lower order.

We begin with a lemma:

\begin{lem}\label{lemme4}
If $\Phi$ satisfies the differential equation $(\mathfrak{E}_{nm}^n)$
with $n, m \in \mathbb{N}^\star$, then $\Psi$ defined by
$$\Psi (x, \ant) : = x^\frac{n-1}{2n} 
\Phi\left((n^{\frac{2}{m}} x)^{\frac{1}{n}}, \an \right)$$
satifies the differential equation $(\mathfrak{E}_{m})$ with 
$\a = \ant$, that is:
\begin{equation}\label{bof2}
x^2 \frac{d^2}{d x^2}\Psi(x, \ant)= P_m(x, \ant)\Psi(x, \ant).
\end{equation}
\end{lem}

\begin{proof}
We consider the transformation
$$\Psi (x, \ant)  = x^\alpha \Phi(\lambda x^{\frac{1}{n}}, \an)$$
with 
$\displaystyle \alpha = \frac{n-1}{2n}$.
 Then
$$x^2\Psi^{\prime \prime} (x, \ant) = 
\frac{\lambda^2}{n^2} x^{\alpha +\frac{2}{n}}
\Phi^{\prime \prime} (\lambda x^{\frac{1}{n}}, \an) +
\alpha (\alpha -1)x^\alpha\Phi (\lambda x^{\frac{1}{n}}, \an).$$
Assuming that 
$\displaystyle x^2\Phi^{\prime \prime} (x, \an) = P_{nm}(x, \an) 
\Phi (x, \an)$,
one gets
$$x^2\Psi^{\prime \prime} (x, \ant) = 
\left(\frac{1}{n^2} P_{nm}(\lambda x^{\frac{1}{n}}, \an) 
+\alpha (\alpha -1) \right)\Psi (x, \ant).$$
We choose 
$\displaystyle \lambda = n^{\frac{2}{nm}}$ to get the statement.
\end{proof}

\begin{notation}
We note  $C_k^n (\an)$ and 
$\widetilde{C}_k^n (\an)$, $k \in \mathbb{Z}$,  the
 Stokes-Sibuya coefficients associated with equation 
$(\mathfrak{E}_{nm}^n)$.
\end{notation}

The above lemma induces the following corollary:

\begin{cor}\label{Beaucor0}
The Stokes-Sibuya coefficients $C_0^n (\an)$ and 
$\widetilde{C}_0^n (\an)$ associated with equation 
$(\mathfrak{E}_{nm}^n)$ are related to 
the Stokes-Sibuya coefficients $C_0$ and 
$\widetilde{C}_0$ of equation $(\mathfrak{E}_{m})$ by:
\begin{equation}\label{func37bis}
\begin{array}{c}
\displaystyle C_0^n(\an) = 
\omega^{\frac{n-1}{2n}}C_0\left(\ant\right)\\
\\
\displaystyle \widetilde{C}_0^n(\an) =
\omega^{\frac{n-1}{n}}\widetilde{C}_0\left(\ant\right)
\end{array}
\end{equation}
where $\displaystyle \omega = e^{\frac{2i\pi}{m}}$. 
\end{cor}

\begin{proof}
We note $\Phi_0$ 
 the solution of ($\mathfrak{E}_m$) which is characterized by its asymptotics
$$ T\Phi_0(x, \a)=
x^{r_m(\a) } e^{-S_m(x, \a)}\left(1+o(1)\right)$$
at infinity in the sector $\Sigma_0 = \{|x|>0, \, 
 |\arg(x) | < \frac{3\pi}{m}\}$ 
(where $r_m=r$ and $S_m=S$ in  theorem \ref{Sibuya}). 
The Stokes-Sibuya coefficients
 ${C}_0$ and $\widetilde{C}_0$ are defined by:
\begin{equation}\label{func31}
\begin{array}{c}
\Phi_0(\omega^{-1}x, \omega^{-1} .\a) = \\
\\
C_0(\a) \Phi_0(x,\a) + \widetilde{C}_0(\a) 
\Phi_0(\omega x, \omega . \a)
\end{array}
\end{equation}
with $\displaystyle \omega = e^{\frac{2i\pi}{m}}$.
We note  $\Phi_0^n$ its analog for equation 
$(\mathfrak{E}_{nm}^n)$, so that
$$ T\Phi_0^n(x, \an)=
x^{r_{nm}(\an) } e^{-S_{nm}(x, \an)}\left(1+o(1)\right)$$
at infinity in the sector $\Sigma_0^n = \{|x|>0, \, 
 |\arg(x) | < \frac{3\pi}{nm}\}$ 
(where $r_{nm}=r$ and $S_{nm}=S$ in  theorem \ref{Sibuya}), and
\begin{equation}\label{func32}
\Phi_0^n(\omega_n^{-1}x, \omega_n^{-1} .\an )
= 
C_0^n(\an) \Phi_0^n(x,\an) + 
\widetilde{C}_0^n(\an) 
\Phi_0^n(\omega_n x; \omega_n .\an )
\end{equation}
with $\displaystyle \omega_n = e^{\frac{2i\pi}{mn}}$.
 Introducing, with lemma \ref{lemme4}, the function
\begin{equation}\label{func33}
\Psi_0 (x, \ant)  = x^\frac{n-1}{2n} \Phi_0^n
\left((n^{\frac{2}{m}} x)^{\frac{1}{n}}, \an \right)
\end{equation}
we get a solution of $(\mathfrak{E}_{m})$ such that
$$ T\Psi_0 (x, \ant)=n^{\frac{2}{nm}r_{nm}(\an) } 
x^{\frac{1}{n} r_{nm}(\an) +\frac{n-1}{2n}} 
e^{-S_{nm}((n^{\frac{2}{m}} x)^{\frac{1}{n}}, \an)}
\left(1+o(1)\right)$$
at infinity in the sector $\Sigma_0$. One easily checks that
$\displaystyle S_{nm}((n^{\frac{2}{m}} x)^{\frac{1}{n}}, \an) =
S_m(x, \ant)$ and $\displaystyle \frac{1}{n} r_{nm}(\an)  =
r_m(\ant)-\frac{n-1}{2n}$. This means that
\begin{equation}\label{func34}
\Psi_0 (x, \ant)=
n^{\frac{2}{m}r_m \left(\ant\right) - \frac{n-1}{mn}}
\Phi_0\left(x, \ant\right).
\end{equation}
From (\ref{func32}), one observes that
\begin{equation}\label{func35}
\begin{array}{c}
\omega^\frac{n-1}{2n}\left(\omega^{-1}x\right)^\frac{n-1}{2n}\Phi_0^n((n^{\frac{2}{m}}
\omega^{-1} x)^{\frac{1}{n}}, \omega^{-1}_n .\an )
= \\
\\
C_0^n(\an) 
x^\frac{n-1}{2n}\Phi_0^n((n^{\frac{2}{m}} x)^{\frac{1}{n}},\an) +
\omega^{-\frac{n-1}{2n}} \widetilde{C}_0^n(\an) 
\left(\omega x\right)^\frac{n-1}{2n}\Phi_0^n(( n^{\frac{2}{m}}\omega
x)^{\frac{1}{n}}; \omega_n .\an)
\end{array}
\end{equation}
so that, by (\ref{func33}),
\begin{equation}\label{func36}
\begin{array}{c}
\omega^\frac{n-1}{2n}\Psi_0(\omega^{-1} x, \omega^{-1} . \ant)
= \\
\\
C_0^n(\an) \Psi_0( x , \ant) +
\omega^{-\frac{n-1}{2n}} \widetilde{C}_0^n(\an) 
\Psi_0(\omega x , \omega . \ant).
\end{array}
\end{equation}
Using (\ref{func34}), one can compare this last equation with
(\ref{func31}) to get the final result.
\end{proof}

Lemma  \ref{lemme4} can be used also to compare the 
$0\infty$ connection matrices. We shall use the following notations:

\begin{notation}
We note 
$\widetilde{L}_k^n(\anp, p(a_{mn}))$ and 
${L}_k^n(\anp, p(a_{mn}))$, $k \in \mathbb{Z}$, 
the coefficients of the
$0\infty$ connection matrices associated with equation 
$(\mathfrak{E}_{nm}^n)$ with
$\displaystyle p(a_{mn}) = (1+4a_{mn})^\frac{1}{2}$.
\end{notation}

\begin{cor}\label{Beaucor}
When $\displaystyle -\frac{p(a_{mn})}{n} \notin \mathbb{N}$, 
 \begin{equation}\label{func49}
\begin{array}{c}
\displaystyle \widetilde{L}_0^n (\anp, p(a_{mn})) = 
  e^{ \frac{i\pi}{m}(1-\frac{1}{n})}
 n^{ -\frac{2}{m}r \left(\ant\right) 
+\frac{p(a_{mn})}{mn} + \frac{1}{m}-1}
\widetilde{L}_0 \left(\antp, \frac{p(a_{mn})}{n}\right).
\end{array}
\end{equation}
\end{cor}

\begin{proof}
We use the notations of theorem \ref{Fuchs0}. We introduce the
solution $f_1(x, \a)$ of ($\mathfrak{E}_m$) which reads
\begin{equation}\label{func38}
f_1(x, \ap, p(a_m)) = x^{s\left(p(a_m) \right)} g_1(x, \ap, p(a_m)),
\end{equation}
where we assume that $\displaystyle 
p(a_m)=(1+4a_m)^{\frac{1}{2}}\notin -\mathbb{N}$. 
In the same way, we note $f_1^n (x, \anp, p(a_{mn}))$ the solution of 
$(\mathfrak{E}_{nm}^n)$ which
can be written as
\begin{equation}\label{func39}
 f_1^n (x, \anp, p(a_{mn})) = x^{s\left(p(a_{mn})\right)}
g^n(x, \anp, p(a_{mn}))
\end{equation}
under the condition $p(a_{mn}) \notin -\mathbb{N}$.
Following lemma \ref{lemme4}, we define
\begin{equation}\label{func40}
F_1 (x, \ant)  = x^\frac{n-1}{2n} f_1^n
\left((n^{\frac{2}{m}} x)^{\frac{1}{n}}, \anp,  p(a_{mn}) \right)
\end{equation}
which is solution of ($\mathfrak{E}_m$) with parameter $\ant$. 
One easily checks that, necessarily,
\begin{equation}\label{func41}
F_1 (x, \ant)  =
n^{\frac{2}{mn}s\left(p(a_{mn})\right)} 
f_1\left(x, \antp, \frac{p(a_{mn})}{n}\right).
\end{equation}
In other words,
\begin{equation}\label{func44}
f_1\left(x, \antp, \frac{p(a_{mn})}{n}\right) = 
 n^{-\frac{2}{mn}s\left(p(a_{mn})\right)} 
x^\frac{n-1}{2n} f_1^n
\left((n^{\frac{2}{m}} x)^{\frac{1}{n}}, \anp,  p(a_{mn}) \right).
\end{equation}
Note that this equality allows to extends analytically  $ f_1^n
(x, \anp, p(a_{mn}))$ for $\displaystyle \frac{p(a_{mn})}{n} \notin
-\mathbb{N}$, and this translates to the $\widetilde{L}_k^n (\anp,
p(a_{mn}))$ as well.
 
We consider the $0\infty$ connection matrices 
$M_1$ and $M_1^n$ associated with 
($\mathfrak{E}_m$) and $(\mathfrak{E}_{nm}^n)$ respectively. We have
\begin{equation}\label{func42}
\begin{array}{c}
\widetilde{L}_1 (\ap, p(a_m)) = \\
-\frac{1}{p(a_m)}
\left(
f_1(x, \ap, p(a_m)) \Phi_0^\prime (x, \a) -
 f_1^\prime(x, \ap, p(a_m)) \Phi_0 (x,  \a)
\right)
\end{array}
\end{equation}
and
\begin{equation}\label{func43}
\begin{array}{c}
\widetilde{L}_1^n (\anp, p(a_{mn})) = \\
-\frac{1}{p(a_{mn})}
\left(
f_1^n(x, \anp, p(a_{mn})) {\Phi_0^n}^\prime (x, \an) -
 {f_1^n}^\prime(x, \anp, p(a_{mn})) \Phi_0^n (x, \an)
\right)
\end{array}
\end{equation}
where $\Phi_0^n$ has been defined in the proof of corollary 
\ref{Beaucor0}. By (\ref{func33}) and (\ref{func34}), we know that
\begin{equation}\label{func45} 
\Phi_0\left(x, \ant \right) =
\displaystyle 
n^{-\frac{2}{m}r \left(\ant\right) + \frac{n-1}{mn}}
 x^\frac{n-1}{2n} \Phi_0^n
\left((n^{\frac{2}{m}} x)^{\frac{1}{n}}, \an \right).
\end{equation}
Equation (\ref{func42}) together with (\ref{func44}) and
 (\ref{func45}) yields:
\begin{equation}\label{func46}
\begin{array}{c}
\widetilde{L}_1 \left( \antp, \frac{p(a_{mn})}{n} \right) = 
-n^{-\frac{2}{m}r \left(\ant \right) 
-\frac{p(a_{mn})}{mn} + \frac{1}{m} 
}\times \\
\\
\frac{1}{p(a_{mn})}
\left(
f_1^n((n^{\frac{2}{m}} x)^{\frac{1}{n}}, \anp, p(a_{mn}))
  {\Phi_0^n}^\prime ((n^{\frac{2}{m}} x)^{\frac{1}{n}}, \an) - \right.\\
\\
\left.
 {f_1^n}^\prime((n^{\frac{2}{m}} x)^{\frac{1}{n}}, \anp, p(a_{mn})) 
\Phi_0^n ((n^{\frac{2}{m}} x)^{\frac{1}{n}}, \an)
\right).
\end{array}
\end{equation}
Comparing (\ref{func46}) with (\ref{func43}), we  get
 \begin{equation}\label{func47}
\widetilde{L}_1^n (\anp, p(a_{mn})) = 
n^{\frac{2}{m}r \left(\ant \right) 
+\frac{p(a_{mn})}{mn} - \frac{1}{m} }
\widetilde{L}_1 \left( \antp, \frac{p(a_{mn})}{n} \right).
\end{equation}
Using formula (\ref{Mk20k}) of theorem \ref{0infty0}, 
we eventually find, since $\displaystyle 
r(\a)+r(\omega .\a) = 1-\frac{m}{2}$,
 \begin{equation}\label{func48}
\widetilde{L}_0^n (\anp, p(a_{mn}))= 
e^{\frac{i\pi}{m}(1-\frac{1}{n})}
n^{-\frac{2}{m}r \left(\ant\right) 
+\frac{p(a_{mn})}{mn} + \frac{1}{m}-1}
\widetilde{L}_0 \left( \antp, \frac{p(a_{mn})}{n} \right).
\end{equation}
\end{proof}


\section{Some applications}\label{sec6}

\subsection{Application for a special class}

Some simplifications occur when $\ap = 0$, allowing to get the
following proposition:

\begin{prop}\label{casrestreint}
We consider $(\mathfrak{E}_{m})$ on restriction to  $\ap =0$. Then
$$\mathfrak{S}_0 (0, a_m) =
\left(
\begin{array}{cc} 
\displaystyle  2 e ^{-\frac{i\pi}{m}} 
\cos \left(\pi \frac{p}{m} \right) &  - e ^{-\frac{2i\pi}{m}}\\ 
1 & 0 
\end{array} 
\right)$$
where $ \displaystyle p = (1+4a_m)^{1/2}$. Furthermore, 
for  $\displaystyle \frac{p}{m} 
\notin \mathbb{Z}$, the $0\infty$ connection matrix $M_0$ is
given by
$$\displaystyle M_0 (0, p) =
\left(
\begin{array}{cc} 
\displaystyle e ^{\beta_m (-p)}\frac{\omega^{-\frac{1}{2}}}{\sqrt{m\pi} } 
\Gamma\left(- \frac{p}{m} \right)  & 
\displaystyle e ^{\beta_m (p)}\frac{\omega^{-\frac{1}{2}}}{\sqrt{m\pi} }\Gamma \left( \frac{p}{m}\right) \\
\displaystyle \omega^{s(p)}e ^{\beta_m (-p)}
\frac{\omega^{-\frac{1}{2}}}{\sqrt{m\pi} } \Gamma\left(- \frac{p}{m} \right) & 
\displaystyle \omega^{s(-p)} e ^{\beta_m (p)}\frac{\omega^{-\frac{1}{2}}}{\sqrt{m\pi} }\Gamma \left( \frac{p}{m}\right)
\end{array} 
\right) 
$$
where $\displaystyle s(p) = \frac{1+p}{2}$, 
while $\beta_m (p)$ is an odd function, entire in $p$, such
that for all $k\ \in \mathbb{N}^\star$,
$$e ^{\beta_m (km)} = \pm m^{k}.$$
\end{prop}

\begin{rem}\label{996}
We shall see in a moment (\S \ref{cacestm1}) by other means  that 
$\displaystyle 
\widetilde{L}_0 (\a) = - e^{i\pi p} \frac{\Gamma (p)}{\sqrt{\pi}}$
when $m=1$. Moreover $\displaystyle r (\a) = \frac{1}{4}$ for $m=1$. 
Applying corollary \ref{Beaucor} with $\ap = 0$, we  deduce that, 
$$\widetilde{L}_0 (0 , p ) =  e^{-\frac{i\pi}{m}} m^{\frac{p}{m}}
e^{i\pi \frac{p}{m}} 
\frac{\Gamma (\frac{p}{m})}{\sqrt{m \pi}}$$
while
$${L}_0 (0 , p) = e^{-\frac{i\pi}{m}} m^{-\frac{p}{m}}
e^{-i\pi \frac{p}{m}} 
\frac{\Gamma (-\frac{p}{m})}{\sqrt{m \pi}}.$$
\end{rem}

\begin{proof}
We note that when $\ap = 0$, then 
 $\omega . \a = \a$ and $\displaystyle r ( \a) = \frac{1}{2}
- \frac{m}{4} $. This has two consequences.

Firstly, when  applying corollary \ref{cestbeau}, one immediately
gets, for $p \notin -\mathbb{N}$,
$$\displaystyle {C}_0 (0, a_m) = 2 e ^{-\frac{i\pi}{m}} 
\cos \left(\pi \frac{p}{m} \right).$$
Since  $C_0 (\a)$ is an entire function in $\a$,  
the above formula extends to all $a_m$,
by analytic continuation.

Secondly, formula (\ref{genial2}) of  theorem
  \ref{genial1} reduces into
$$\frac{1}{L_0(0, p) \widetilde{L}_0(0, p)} = 
- \omega p \sin \left(\pi \frac{p}{m} \right),$$
which resembles the Euler reflection formula $\displaystyle 
\frac{\pi}{\Gamma(z) \Gamma(-z)} = 
-  z \sin \left(\pi z \right)$. Since by proposition \ref{rem41},  the
restriction to $p \notin \mathbb{Z}$ of the function
$L_0(0, p)$ ({\sl resp.} $\widetilde{L}_0 (0, p)$) has a 
 meromorphic continuation in $p$, with at most simple
poles at $p \in m\mathbb{N}$ ({\sl resp.} $-p \in m\mathbb{N}$), we can
write
$$ L_0 (0, p) = \alpha (p) \frac{\omega^{-\frac{1}{2}}}{\sqrt{m\pi} } 
\Gamma\left(- \frac{p}{m} \right)$$
and 
$$  \widetilde{L}_0 (0, p) = \frac{1}{\alpha(p)}
\frac{\omega^{-\frac{1}{2}}}{\sqrt{m\pi} }\Gamma \left( \frac{p}{m}\right)$$  
where $\alpha (p)$ is a nowhere vanishing entire function of $p$. In
other words, 
$$\alpha (p)= e^{-\beta_m(p)}$$
 with $\beta_m (p)$ an  entire function. Furthermore, 
by theorem \ref{0infty0} again,  we know that $
 L_0(\ap, -p)=\widetilde{L}_0(\ap, p)$. This means that
 $\beta_m$ can be chosen as an odd function.

Still by theorem  \ref{0infty0}, we know that $\widetilde{L}_0(\ap,
p)$ extends analytically to $p \in \mathbb{N}^\star$.
 Moreover, when
$$p = km, \, \, \, k \in \mathbb{N}^\star,$$
then formula (\ref{genial2bis}) of theorem \ref{genial1} gives 
$$\widetilde{L}_0^2(0, p) \, |_{p = km} = (-1)^{k+1}
\frac{m\omega^{-1}}{\pi p\lambda(0, p)}\, |_{p = km}.$$
By remark \ref{rem50}, we know that
$$\lambda (0, p)\, |_{p = km} = \frac{(-1)^{k+1}}{m^{2k-1}k\Gamma(k)^2},$$
and therefore,
$$\widetilde{L}_0^2(0, p)\, |_{p = km} = m^{2k}
\frac{\omega^{-1} \Gamma^2(k)}{m \pi},$$
that is
$$\widetilde{L}_0 (0, p) \, |_{p = km} = 
\pm m^{k} \frac{\omega^{-1/2} \Gamma(k)}{\sqrt{m \pi}}.$$
\end{proof}

\subsection{Application when $m=2$ and consequences.}

We consider the case $m=2$, so that $\omega = e^{i\pi}$ and
$\displaystyle r (\a) =  -\frac{a_1}{2}$. 

On the one hand, corollary \ref{CORFUNCT} implies
$$
C_0(\a) C_1( \a) + \widetilde{C}_0 (\a) + \widetilde{C}_1 (\a)= 
-2\cos(\pi p) 
$$
with $C_1( \a)  = C_0(\omega. \a)$, 
where,  by  (\ref{qhCk}) of theorem \ref{SibuyaCoeff}, 
$$
\widetilde{C}_0 (\a) = e^{-i\pi a_1}, \hspace{5mm} \widetilde{C}_1 (a)
= \widetilde{C}_0 ( \omega .\a) = e^{i\pi a_1}.
$$
This means that
$$
C_0(\a) C_0(\omega .\a)  = -4\cos\left(\frac{\pi}{2}(p+a_1)\right) 
\cos\left(\frac{\pi}{2}(p-a_1)\right).
$$
 On the other hand, we know by corollary \ref{cestbeau} that (for $\a$
 generic): 
\begin{equation}\label{cestbeaubis}
C_0 (\a) = -2i
e^{-i\pi\frac{a_1}{2}}\cos\left(\frac{\pi}{2}(p+a_1)\right) 
\frac{\widetilde{L}_0(\ap, p)}{\widetilde{L}_0 (\omega. \ap, p)}.
\end{equation}
Also, by formula (\ref{genial2bis}) of theorem \ref{genial1}, we have, when
$p \in \mathbb{N}^\star$ and for $\ap = a_1$ generic,
$$
\pi p    \lambda (\ap, p) =
2\frac{\cos \left( \frac{\pi}{2}(p-a_1)\right)}
{ \widetilde{L}_0(\ap, p) \widetilde{L}_0 (\omega. \ap, p)}\, |_{p \in
  \mathbb{N}^\star}. 
$$
By remark \ref{rem50}, 
$\displaystyle
\left\{
\begin{array}{l}
\displaystyle 
\lambda (\ap, p) \, |_{p = 0 \mod 2} =
\frac{(-1)^{p+1}}{p\Gamma(p)^2}\prod_{k=1}^{p/2} (a_1 +2k-1)(a_1
-2k+1)\\
\\
\displaystyle 
\lambda (\ap, p) \, |_{p = 1 \mod 2} =
\frac{(-1)^{p+1}}{p\Gamma(p)^2}a_1\prod_{k=1}^{(p-1)/2} (a_1 +2k)(a_1
-2k),\\
\end{array}
\right. 
$ so that
$$\widetilde{L}_0(\ap, p) 
\widetilde{L}_0 (\omega. \ap, p)\, |_{p \in \mathbb{N}^\star}
=$$
$$
\left\{
\begin{array}{l}
\displaystyle 
2(-1)^{\frac{p+2}{2}}\frac{\cos \left( \frac{\pi}{2}a_1\right)\Gamma(p)^2}
{\pi \prod_{k=1}^{p/2} (a_1 +2k-1)(a_1
-2k+1)}, \hspace{2mm} p \mbox{  even } \\
\\
\displaystyle
2(-1)^{\frac{p+1}{2}}\frac{\sin \left( \frac{\pi}{2}a_1\right)\Gamma(p)^2}
{\pi a_1\prod_{k=1}^{(p-1)/2} (a_1 +2k)(a_1
-2k)}, \hspace{2mm} p \mbox{ odd. } \\
\end{array}
\right.
$$
This can be also written as
\begin{equation}\label{pasbanal}
\displaystyle
\widetilde{L}_0(\ap, p) 
\widetilde{L}_0 (\omega. \ap, p)\, |_{p \in \mathbb{N}^\star}
= -2^{-p+1} e^{i \pi p} \frac{\Gamma(p)^2}{\Gamma(\frac{p}{2} + \frac{a_1}{2}+
  \frac{1}{2}) \Gamma(\frac{p}{2} - \frac{a_1}{2}+
  \frac{1}{2})}.
\end{equation}
At this point, we can use
 the following lemma, whose easy proof is left to the
reader:

\begin{lem}
When $p+a_1+1=-2N$ with $N \in \mathbb{N}$, 
then, for $p \notin -\mathbb{N}^\star$,  
$$f_1(x, a_1, p) = x^{s(p)} e^{-x}\sum_{n=0}^N
\frac{\Gamma(p+1)Q_n(a_1, p)}{\Gamma(n+p+1)} x^n$$
where the $Q_n(a_1, p) \in \mathbb{C}[a_1, p]$ are defined by:
$$\displaystyle 
\left\{
\begin{array}{l}
Q_0(a_1, p) =1 \\
(a_1 +p+1 + 2n)Q_{n}(a_1, p) -(n+1)Q_{n+1}(a_1, p)=0, \, \, n \geq 0.
\end{array}
\right.
$$
In particular,
$\displaystyle f_1(x, a_1, p)  = (-1)^N 2^N\frac{\Gamma(p+1)}{\Gamma(N+p+1)}
\Phi_0 (x, \a)$.
\end{lem}

This lemma implies that 
$$\widetilde{L}_0 (\omega. \ap, p) = 0 \, \, \, \mbox{ when }
 \, \, \,
\left\{ 
\begin{array}{c}
p+a_1+1 \in -2\mathbb{N} \\
p \notin -\mathbb{N}^\star
\end{array}
\right. ,
$$
that is also
\begin{equation}\label{10121}
\widetilde{L}_0 (\ap, p) = 0 \, \, \, \mbox{ when }
 \, \, \,
\left\{ 
\begin{array}{c}
p-a_1+1 \in -2\mathbb{N} \\
p \notin -\mathbb{N}^\star .
\end{array}
\right. 
\end{equation}
Since the right-hand side of
(\ref{pasbanal}) has only simple zeros when  $p+a_1+1 \in
-2\mathbb{N}$, we can write:
\begin{equation}\label{101}
\begin{array}{l}
\displaystyle \widetilde{L}_0(\ap, p)\, |_{p \in \mathbb{N}^\star} = 
-i 2^{-\frac{p-1}{2}}e^{i \pi \frac{p}{2}}\frac{\Gamma(p)}
{ \Gamma(\frac{p}{2} - \frac{a_1}{2}+ \frac{1}{2})}
e^{\beta(a_1,p)}\, |_{p \in \mathbb{N}^\star} \\
 \mbox{with}
\hspace{3mm} \beta(-a_1,p) = -\beta(a_1,p).
\end{array}
\end{equation}
Now when $p \notin \mathbb{Z}$, the coefficient $L_0(\ap, p)$ can be
derived from formula (\ref{genial2}) of  theorem \ref{genial1}. This
gives:
\begin{equation}\label{qcq}
L_0(\ap, p)\widetilde{L}_0(\omega. \ap, p) = 2
\frac{\cos\left(\frac{\pi}{2}(p+a_1)\right)}{p \sin(\pi p) }.
\end{equation}
We recall also that $L_0(\ap, p)$ can be derived from
$\widetilde{L}_0(\ap, p)$
just by changing $p$ into $-p$. Using (\ref{10121}) and 
the known analytic properties of 
$L_0(\ap, p)$ and  $\widetilde{L}_0(\ap, p)$ described by proposition 
\ref{rem41},  equation (\ref{qcq}) shows that equation 
(\ref{101}) can be extended 
to all $(\ap, p) \in \mathbb{C}^2$ with $p \notin -\mathbb{N}^\star$, 
$$\widetilde{L}_0(\ap, p) = 
-i 2^{-\frac{p-1}{2}}e^{i \pi \frac{p}{2}}\frac{\Gamma(p)}
{ \Gamma(\frac{p}{2} - \frac{a_1}{2}+ \frac{1}{2})}
e^{\beta(a_1,p)} 
$$
where  $\beta$ can be chosen as  an entire function satisfying
$$\beta(-a_1,p) = -\beta(a_1,p) \hspace{5mm} \mbox{and} \hspace{5mm} 
\beta(a_1,- p) = \beta(a_1,p).
$$
Finally, formula \ref{cestbeaubis} reduces into:
$$
C_0 (\a) = -2i
e^{-i\pi\frac{a_1}{2}}\cos\left(\frac{\pi}{2}(p+a_1)\right) 
\frac{\Gamma(\frac{p}{2} + \frac{a_1}{2}+ \frac{1}{2})}
{\Gamma(\frac{p}{2} - \frac{a_1}{2}+ \frac{1}{2})}e^{2\beta(a_1,p)}.
$$

To summarize:

\begin{prop}\label{proppourm2}
We assume $m=2$. Then the Stokes-Sibuya mutiplier $C_0$ may be written
as
\begin{equation}\label{aconfirmer}
C_0 (\a) = -2i
e^{-i\pi\frac{a_1}{2}}\cos\left(\frac{\pi}{2}(p+a_1)\right) 
\frac{\Gamma(\frac{p}{2} + \frac{a_1}{2}+ \frac{1}{2})}
{\Gamma(\frac{p}{2} - \frac{a_1}{2}+ \frac{1}{2})}e^{2\beta(a_1,p)}
\hspace{3mm} \mbox{and} \hspace{3mm}
\widetilde{C}_0 (\a) = e^{-i\pi a_1}
\end{equation}
where $\beta$ is an entire function satisfying $
\beta(a_1,p) = \beta(a_1,-p) = -\beta(-a_1,p)$.\\
Moreover, the coefficients of the $0\infty$ connection matrix $M_0$
of theorem  \ref{0infty0} satisfy, for $p \notin \mathbb{Z}$,
\begin{equation}\label{genial21bis}
\left\{
\begin{array}{l}
\displaystyle
\widetilde{L}_0(a_1, p) = 
-i 2^{-\frac{p-1}{2}}e^{i \pi \frac{p}{2}}\frac{\Gamma(p)}
{ \Gamma(\frac{p}{2} - \frac{a_1}{2}+ \frac{1}{2})}
e^{\beta(a_1,p)}\\
\\
\displaystyle
L_0(a_1, p)\widetilde{L}_0(\omega a_1, p) = 2
\frac{\cos\left(\frac{\pi}{2}(p+a_1)\right)}{p \sin(\pi p) }.
\end{array}
\right.
\end{equation}
\end{prop}

\begin{rem}
The above proposition is interesting since, for instance, it already provides
the location of the zeroes of $C_0$ and of the other Stokes-Sibuya
coefficients. However, one can be more precise, using the Whittaker special
functions. We shall see in a moment (\S \ref{cacestm2}) that 
$$\beta(a_1,p) = -a_1.$$
\end{rem}

With this remark and corollaries \ref{Beaucor0} and \ref{Beaucor},
 proposition \ref{proppourm2} implies the following
consequences:

\begin{cor}\label{Cor56}
We consider the differential equation
$$
x^2 \frac{d^2}{d x^2}\Phi= \left(x^{2n} + a_nx^n
+a_{2n} \right)\Phi
\leqno{(\mathfrak{E}_{2n}^n)}
$$
where $n \in \mathbb{N}^\star$. Then,
$$\left\{
\begin{array}{l}
\displaystyle
C_0^n(\an)=2e^{-\frac{i\pi}{2n}}  
e ^{ -i\pi\frac{a_n}{2n}} 2^{-\frac{a_n}{n}} 
\frac{\Gamma(\frac{p}{2n}+\frac{a_n}{2n}+\frac{1}{2})}{\Gamma(\frac{p}{2n}-\frac{a_n}{2n}+\frac{1}{2})}
\cos\left((\frac{p}{2n}+\frac{a_n}{2n})\pi\right)\\
\\
\widetilde{C}_0^n (\an) = -e^{-\frac{i\pi }{n}}e^{-i\pi \frac{a_n}{n}}
\end{array}
\right. 
$$
where $\displaystyle  p = (1+4a_{2n})^\frac{1}{2}$. Moreover, when 
$p \notin -n\mathbb{N}$,
$$
\displaystyle
\widetilde{L}_0^n (\anp, p) = e^{-\frac{i \pi }{2n}}
 e^{i \pi \frac{p}{2n}} 
\left( \frac{n}{2}\right)^{ \frac{a_n}{2n}+\frac{p}{2n} -\frac{1}{2}}
\frac{ \Gamma( \frac{p}{n} ) }{\Gamma(\frac{p}{2n}-\frac{a_n}{2n}+\frac{1}{2})}
.$$
\end{cor}

\subsection{Application when $m \geq 3$.}

As a matter of fact, 
we shall only discuss the cases $m=3$ and $m=4$ to show what kind of
information we can extract from our analysis.

\subsubsection{The case $m=3$.}

In a sense,  $m=3$ is the first interesting case, since no special
function solution of $(\mathfrak{E}_{3})$
is known. \\
Here we have  $\omega = e^{\frac{2i\pi}{3}}$ and
$\displaystyle r (\a)  = -\frac{1}{4}$ is a constant function.

We first apply corollary \ref{CORFUNCT}, to get a functional relation
between the Stokes-Sibuya multipliers:
$$
C_0(\a) C_1(\a)C_2(\a) + \widetilde{C}_0(\a) C_2(\a) + 
\widetilde{C}_1(\a)  C_0(\a)  
+\widetilde{C}_2(\a)  C_1(\a)=  -2\cos(\pi p) 
$$
where,  by  (\ref{qhCk}) of theorem \ref{SibuyaCoeff}, 
$$
\widetilde{C}_0 (\a) = \widetilde{C}_1 (\a) = \widetilde{C}_2 (\a) =
e^{\frac{i\pi}{3}}.
$$
Applying now corollary \ref{cestbeau}, we find, for all $\ap \in
\mathbb{C}^2$ and $p \notin -\mathbb{N}$:
\begin{equation}\label{m32}
\widetilde{L}_0(\omega . \ap, p) C_0 (\a) =\omega^{-\frac{1}{2}}
  \left( 
\widetilde{L}_0( \ap, p) \omega^{ \frac{p}{2}} +
\widetilde{L}_0( \omega^2 .\ap, p) \omega^{- \frac{p}{2}} 
 \right).
\end{equation}

We concentrate on the case $p \in \mathbb{N}^\star$. 
By formula (\ref{genial2bis}) of theorem \ref{genial1}, we get
\begin{equation}\label{m35}
\begin{array}{c}
\displaystyle i\pi p  \omega^{\frac{7}{4} + \frac{p}{2}}  \lambda (\ap , p) 
\widetilde{L}_0 (\ap , p)
\widetilde{L}_0 (\omega .\ap , p) \widetilde{L}_0 (\omega^2 .\ap , p)
=   \\
\\
\displaystyle
\widetilde{L}_0 (\omega .\ap , p) + 
\omega^{ p}\widetilde{L}_0 (\omega ^2 .\ap , p) +
\omega^{2p} \widetilde{L}_0 (\ap , p).
\end{array}
\end{equation}
We now add the assumption that  $\ap$ has been chosen so that
$$\lambda (\ap , p)\widetilde{L}_0 (\ap , p)\widetilde{L}_0 (\omega .\ap , p) 
\widetilde{L}_0 (\omega^2 .\ap , p)=0.$$
 Using the remark that  $\widetilde{L}_0 (\omega .\ap , p)=0$ implies
$\widetilde{L}_0 (\ap , p)
 \widetilde{L}_0 (\omega^2 .\ap , p) \neq 0$ necessarily (otherwise, one of the
two  $0\infty$ connection matrices $M_0$ or $M_1$ is not invertible, which
 is absurd), equations   (\ref{m35}) and (\ref{m32}) imply that:
$$
 C_0 (\a) =  C_1 (\a) =  C_2 (\a) =
-\omega^{-\frac{1}{2}-\frac{3p}{2} }= (-1)^{p+1}e^{-\frac{i\pi}{3}}.
$$

We summarize our results:

\begin{prop}\label{pourmeq3}
We assume $m=3$. Then the Stokes-Sibuya multiplier $C_0(\a)$ satisfies
the functional equation
\begin{equation}\label{m39}
C_0(\a) C_0(\omega .\a)C_0(\omega^2 .\a) + 
e^{\frac{i\pi}{3}} \left(  C_0(\a) + 
 C_0(\omega .\a) + C_0(\omega^2 .\a) \right)=  -2\cos(\pi p) 
\end{equation}
 with $p = (1+4a_3)^\frac{1}{2}$ and $\omega = e^{\frac{2i\pi}{3}}$, whereas
\begin{equation}\label{m40}
\widetilde{C}_0 (\a) = e^{\frac{i\pi}{3}}.
\end{equation}
Moreover, when $\displaystyle a_3 = \frac{p^2-1}{4}$ with $p \in
\mathbb{N}^\star$, then
$$\lambda (\ap, p)\widetilde{L}_0 (\ap, p)\widetilde{L}_0 (\omega .\ap, p) 
\widetilde{L}_0 (\omega^2 .\ap, p)=0$$
is equivalent to  $C_0$ being a constant, precisely 
$$
C_0 = (-1)^{p+1}e^{-\frac{i\pi}{3}}.
$$
\end{prop}

We note that proposition \ref{pourmeq3} can be
  derived from corollary \ref{CORFUNCT2} when  $\lambda (\ap, p)=0$, while
the particular case $a_1=a_2=0$ is given by proposition
\ref{casrestreint}. 

For a given  $p \in \mathbb{N}^\star$, the case
  $\lambda (\ap, p)=0$ can be seen as an isomonodromic deformation
  condition, since both the monodromy at the origin and the Stokes
  structure are fixed. We get:

\begin{cor}  For $m=3$ and $p \in \mathbb{N}^\star$, 
the condition $\lambda (\ap, p)  = 0$ is an isomonodromic deformation
condition.  
\end{cor}

By computing $\lambda (\ap, p)$ (see proposition \ref{Calambda}), 
one obtains for example, from proposition 
\ref{pourmeq3}: 
\begin{itemize}
\item If $p=1$, then $\lambda(\ap, p) =a_2$. Therefore, for all 
$a_1 \in \mathbb{C}$,
$$C_0 (a_1, 0, 0) =e^{-\frac{i\pi}{3}}.$$
By a Tschirnhaus transformation, this case is equivalent to the Airy
equation. This means also that 
$\widetilde{L}_0 (a_1,0,1) =
\widetilde{L}_0 (0,0,1)$
 so that by  remark \ref{996},
$$\displaystyle \widetilde{L}_0 (a_1,0,1) =  
e^{-\frac{i\pi}{3}} 3^{\frac{1}{3}}
e^{i\pi \frac{1}{3}} 
\frac{\Gamma (\frac{1}{3})}{\sqrt{3 \pi}},$$ 
while
$\displaystyle {L}_0 (a_1,0,1) =  e^{-\frac{i\pi}{3}} 3^{-\frac{1}{3}}
e^{-i\pi \frac{1}{3}} 
\frac{\Gamma (-\frac{1}{3})}{\sqrt{3 \pi}}$.
\item If $p=2$, then $\displaystyle \lambda(\ap, p) =-\frac{a_2^2}{2} +
  \frac{a_1}{2}$. We deduce that, for all $a_2 \in \mathbb{C}$,
$$C_0 (a_2^2, a_2, \frac{3}{4}) =-e^{-\frac{i\pi}{3}}.$$
\item If $p=3$, then $\displaystyle \lambda(\ap, p) =\frac{a_2^3}{12} -
  \frac{a_2a_1}{3} + \frac{1}{3} $. Thus, for all $a_2 \in \mathbb{C}^\star$,
$$C_0 (\frac{4+a_2^3}{4a_2}, a_2, 2) =e^{-\frac{i\pi}{3}}.$$
\end{itemize}

Since $\lambda (\ap, p)$ can be computed exactly for all 
fixed $\displaystyle p \in \mathbb{N}^\star$, it is natural to 
try to get more informations from equation (\ref{m35}). 
The result is a little bit disappointing, as we now explain.

We assume that
\begin{equation}\label{m60} 
\widetilde{L}_0 (\ap, p)\widetilde{L}_0 (\omega .\ap, p) 
\widetilde{L}_0 (\omega^2 .\ap, p) \neq 0.
\end{equation}
We note that, by remark \ref{996}, 
$$\widetilde{L}_0 (0, p) =  e^{-\frac{i\pi}{3}} 3^{\frac{p}{3}}
e^{i\pi \frac{p}{3}} 
\frac{\Gamma (\frac{p}{3})}{\sqrt{3 \pi}}$$
so that hypothesis (\ref{m60}) is valid for $\ap$ in a
neighbourhood of the origin. We can write (\ref{m35}) as
\begin{equation}\label{m61}  
y(\ap, p) +\omega^{p} y(\omega .\ap, p) + \omega^{2p} y(\omega^2 .\ap,
p)
=
-\pi p  \omega^{1 - \frac{p}{2}}  \lambda (\ap, p) 
\end{equation}
 with
$$y(\ap, p) = \frac{1}{\widetilde{L}_0 (\ap, p)
\widetilde{L}_0 (\omega .\ap, p)}.$$
Equation (\ref{m61}) can be thought of as a non-homogenous second order linear
 $q$-difference equation. Unfortunately, we are in the worst
 situation, when $q$ is a root of unity, so that solving (\ref{m61})
 gives very little information. Indeed, we first observe that 
$\displaystyle -\frac{\pi p}{3}  \omega^{1 - \frac{p}{2}}
\lambda (\ap, p)$ is a particular solution of (\ref{m61}), because 
$\lambda (\omega .\ap, p) = \omega^{-p} \lambda (\ap, p)$. 
Therefore, by linearity of
(\ref{m61}), one only needs to solve  the
homogenous equation
\begin{equation}\label{m62}  
y(\ap, p) +\omega^{p} y(\omega .\ap, p) + \omega^{2p} y(\omega^2 .\ap, p)
=0
\end{equation}
in  the space
$\mathbb{C}\{a_1,a_2\}$.
Writing $\displaystyle y(\ap, p) = \sum_{k,l=0}^\infty b_{k,l} a_1^k
a_2^l$, equation (\ref{m62}) is equivalent to (since $\omega^3=e^{2i\pi}$):
$$
\sum_{k,l=0}^\infty (1+\omega^{p + k+2l} + \omega^{2p + 2k+l}) b_{k,l} a_1^k
a_2^l=0.
$$
Thus, $y \in \mathbb{C}\{a_1,a_2\}$ is solution of (\ref{m62})
provided that  $ b_{k,l} =0$ when $p+k+2l = 0 \mod 3$. This corresponds
to a vector-space of infinite dimension !

To end this subsection, we mention \cite{Sla96} for the numerical
computations of the $0\infty$ connection matrices.

\subsubsection{the case $m=4$} 

When $m=4$, we have
$\displaystyle r(\a) =-\frac{1}{2} -\frac{1}{2}a_2+\frac{1}{8}a_1^2$
so that, by theorem \ref{SibuyaCoeff},
$$\widetilde{C}_0 (\a) = i e^{-\frac{i\pi}{2}(a_2 -
  \frac{1}{4}a_1^2)}.$$
Also, by corollary \ref{CORFUNCT} :
\begin{equation}\label{Functrelm4}
\begin{array}{c}
C_0(\a) C_1(\a)C_2(\a)C_3(\a) + \widetilde{C}_0(\a) C_2(\a)C_3(\a) + 
\widetilde{C}_1(\a) C_0(\a)C_3(\a) +\\
\\
\widetilde{C}_2(\a) C_0(\a)C_1(\a) +\widetilde{C}_3(\a) C_1(\a)C_2(\a)
+\widetilde{C}_0(\a)\widetilde{C}_2(\a)+\widetilde{C}_1(\a)\widetilde{C}_3(\a)=
-2\cos(\pi p).
\end{array} 
\end{equation}

We already know, 
by corollary \ref{Cor56} applied with
$n=2$, that 
\begin{equation}\label{weber}
\displaystyle
C_0(0, a_2, 0, a_{4})=
e ^{ -\frac{i\pi}{4}(a_2 +1)} 2^{-\frac{a_2}{2}+1} 
\frac{\Gamma(\frac{p}{4}+\frac{a_2}{4}+\frac{1}{2})}{\Gamma(\frac{p}{4}-\frac{a_2}{4}+\frac{1}{2})}
\cos\left(\frac{\pi}{4}(p+a_2)\right)
\end{equation}
with $p = \displaystyle (1+4a_4)^\frac{1}{2}$.

Propositions similar to  proposition 
\ref{pourmeq3} can be obtained for every $m \geq 3$. In particular for
$m=4$, we show what happens for values of $\a$ such that 
$\displaystyle a_4 = \frac{p^2-1}{4}$ with $p \in
\mathbb{N}^\star$ and $\lambda (\ap, p)=0$. Evaluating the last row of
the product $ \mathfrak{S}_0 (\a) \mathfrak{S}_1 (\a) \cdots 
\mathfrak{S}_{3} (\a)$ and applying  corollary
\ref{CORFUNCT2}, one gets
$$
\begin{array}{c}
\left(
\begin{array}{cc}
\cdots  & \cdots \\
\left( C_1(\a)C_2(\a) + \widetilde{C}_1(\a)\right)C_3(\a)
  +C_1(\a) \widetilde{C}_2(\a) &
\left(C_1(\a)C_2(\a) + \widetilde{C}_1(\a)\right)\widetilde{C}_3(\a)
\end{array}
\right)\\
=(-1)^{p+1}
\left(
\begin{array}{cc}
1 & 0 \\
0 & 1
\end{array}
\right).
\end{array}
$$
We have $\omega = e^{\frac{i\pi}{2}}$ so that 
$\displaystyle \widetilde{C}_2(\a) = 
\widetilde{C}_0(\omega^2.\a) = \widetilde{C}_0(\a)$.
 We thus get
$$ C_0(\a)C_1(\a)  + \widetilde{C}_0 (\a)=
(-1)^{p+1}\widetilde{C}_0^{-1}(\a)
$$
and
$$
\left( C_0(\a)C_1(\a) + \widetilde{C}_0(\a)\right)C_2(\a)
  +C_0(\a) \widetilde{C}_1(\a)=0.
$$
Since 
$\displaystyle \widetilde{C}_0(\a) \widetilde{C}_1(\a) = -1$, we obtain
$$C_0(\a) = (-1)^{p+1}C_2(\a), \hspace{5mm} 
 C_0(\a)C_1(\a) = (-1)^{p+1}\widetilde{C}_0^{-1}(\a) -
\widetilde{C}_0 (\a).$$
Computing $\lambda (\ap, p)$, one obtains for instance:
\begin{itemize}
\item If $p=1$, then $\lambda(\ap, p) =a_3$. Therefore, for all 
$(a_1, a_2)  \in \mathbb{C}^2$,
$$C_0 (a_1, a_2, 0, 0) = C_2 (a_1, a_2, 0, 0), $$
$$
 C_0(a_1, a_2, 0, 0)C_1(a_1, a_2, 0, 0) = -2i \cos \left( \frac{\pi}{2}
(a_2 - \frac{1}{4}a_1^2)\right).
$$
This case corresponds to the Weber equation. By a Tschirnhaus
transformation, one can use equation (\ref{weber}) to get
$$
\begin{array}{c}
\displaystyle
C_0(a_1, a_2, 0, 0)=\\
\displaystyle
2 e^{-\frac{i\pi}{4}}  
e ^{ -i\pi\frac{1}{4}(a_2- \frac{1}{4}a_1^2)} 2^{-\frac{a_2}{2}+ 
\frac{1}{8}a_1^2} 
\frac{\Gamma(\frac{a_2}{4}- \frac{1}{16}a_1^2+
\frac{3}{4})}{\Gamma(-\frac{a_2}{4}+\frac{1}{16}a_1^2+\frac{3}{4})}
\cos\left(\frac{\pi}{4}(a_2- \frac{1}{4}a_1^2+1)\right).
\end{array}
$$
By the Euler reflection formula and the duplication formula for the
Gamma function, one gets the usual well-known formula
(cf. \cite{Sib75}).

Comparing this result with (\ref{weber}), it is tempting to conjecture that:
$$
\begin{array}{c}
\displaystyle
C_0(a_1, a_2, 0, a_4)=  \\
\displaystyle
e ^{ -\frac{i\pi}{4}(a_2- \frac{1}{4}a_1^2 +1)} 2^{-\frac{a_2}{2}+ 
\frac{1}{8}a_1^2+1} 
\frac{\Gamma(\frac{p}{4}+\frac{a_2}{4}- \frac{1}{16}a_1^2+
\frac{1}{2})}{\Gamma(\frac{p}{4} -\frac{a_2}{4}+\frac{1}{16}a_1^2+\frac{1}{2})}
\cos\left(\frac{\pi}{4}(p+a_2- \frac{1}{4}a_1^2)\right).
\end{array}
$$
(This satisfies the functional relation (\ref{Functrelm4})).
\item If $p=2$, then $\displaystyle \lambda(\ap, p) =-\frac{a_3^2}{2} +
  \frac{a_2}{2}$. We deduce that, for all $(a_1, a_3) \in \mathbb{C}^2$,
$$\begin{array}{c}
C_0 (a_1,a_3^2, a_3, \frac{3}{4}) = - C_2(a_1,a_3^2, a_3,
\frac{3}{4}), \\
\\
 C_0(a_1,a_3^2, a_3, \frac{3}{4})C_1(a_1,a_3^2, a_3, \frac{3}{4}) = 
-2 \sin \left( \frac{\pi}{2} (a_3^2 - \frac{1}{4}a_1^2)\right).
\end{array}
$$
\item If $p=3$, then $\displaystyle \lambda(\ap, p) =\frac{a_3^3}{12} -
  \frac{a_2a_3}{3} + \frac{a_1}{3} $. Thus, for all $(a_2,a_3) \in
  \mathbb{C}^2$,
$$\begin{array}{c}
C_0 (-\frac{a_3^3}{4}+a_2a_3, a_2, a_3, 2) = 
C_2 (-\frac{a_3^3}{4}+a_2a_3, a_2, a_3, 2) , \\
\\
 C_0(-\frac{a_3^3}{4}+a_2a_3, a_2, a_3, 2)
C_1(-\frac{a_3^3}{4}+a_2a_3, a_2, a_3, 2) = \\
-2i \cos \left( \frac{\pi}{2}
(a_2 - \frac{1}{4}(-\frac{a_3^3}{4}+a_2a_3)^2)\right).
\end{array}
$$
\end{itemize}


\appendix

\protect \boldmath

\section{Using special functions}\label{sec7}

\subsection{Example 1~: a normal form of Heun's equation }\label{cacestm1}

\protect \unboldmath

We consider the equation
$$x^2\Phi''=(x+a)\Phi.  \leqno{(\mathfrak{E}_1)}$$
This is the simplest case, when $m=1$. In this case, the Stokes-Sibuya
connection matrix $\mathfrak{S}_0$ is given by proposition 
\ref{casrestreint}. This proposition provides also 
the $0\infty$ connection matrix $M_0$, 
up to  an odd entire function of $p=(1+4a)^{\frac{1}{2}}$, which we
are going to compute here thanks to  the fact that the above
normal form $(\mathfrak{E}_1)$  of Heun's equation reduces to a modified
Bessel equation. Indeed, setting
\begin{equation}\label{changeofvar}
  \left \lbrace 
\begin{array}{l}
t=2x^{1/2} \\
\\
\Psi(t)=x^{-1/2} \Phi(x)
\end{array} 
\right. 
\end{equation}
equation  $(\mathfrak{E}_1)$ is converted into the equation~:
\begin{equation}\label{modbessel}
t^2\Psi''(t)+t\Psi'(t)-(t^2+p^2)\Psi(t)=0, \hspace{5mm}
p=(1+4a)^{\frac{1}{2}}, 
\end{equation}
which is a modified Bessel equation of parameter $p$.
Thus, we can use the well-known 
special functions
associated with the modified Bessel equation.

We assume  that $p=(1+4a)^{\frac{1}{2}} \notin \mathbb{Z} $.

With the notations of theorem \ref{Fuchs0}, one easily gets the
fundamental system of solutions $(f_1, f_2)$ of $(\mathfrak{E}_1)$ 
in the form:
\begin{equation}\label{attheorigin}
 \left\lbrace \begin{array}{l}
\displaystyle f_1(x,p) = 
\Gamma(p+1) \sqrt{x} I_p(2\sqrt{x}), \quad 
I_p(t)= \left(\frac{t}{2}\right)^p \sum_{n=0}^{+\infty} 
\frac{1}{n! \Gamma(n+p+1)} \left(\frac{t}{2}\right)^{2n} \\
\\
\displaystyle f_2(x,p) = 
\Gamma(-p+1) \sqrt{x} I_{-p}(2\sqrt{x}), \quad 
I_{-p}(t)= \left(\frac{t}{2}\right)^{-p} \sum_{n=0}^{+\infty} 
\frac{1}{n! \Gamma(n-p+1)} \left(\frac{t}{2}\right)^{2n}.
\end{array} \right. 
\end{equation} 
where $I_p$ (respectively $I_{-p}$) is the modified Bessel function (or
Bessel function of imaginary argument) of order $p$ 
(respectively of order $-p$), see \cite{Olv74}.

\begin{rem}\label{rem109}
We recall that the functions $I_p$ and $I_{-p}$ are very closely 
connected to the Bessel functions of the first kind $J_p$ and
$J_{-p}$ by  (see \cite{Olv74})~:
$$ \left\lbrace \begin{array}{l}
\displaystyle I_p(t)=e^{-\frac{i\pi p}{2}} J_p(it) \\
\\
\displaystyle I_{-p}(t)=e^{\frac{i\pi p}{2}} J_{-p}(it) 
\end{array} \right. $$
\end{rem}

Now, thanks to theorem \ref{Sibuya}, there exists an unique solution $\Phi_0$ 
of $(\mathfrak{E}_1)$, asymptotic at infinity to 
 $T\Phi_0(x,a) =e^{-2\sqrt{x}} x^{\frac{1}{4}} \phi_0(x,a)$ with $\phi_0\in
\mathbb{C}[a][[x^{-\frac{1}{2}}]]$ in the sector 
$-3\pi<arg(x)<3\pi$. Precisely
\begin{equation}\label{Phi0pourm1}
\begin{array}{c}
 \displaystyle T\Phi_0(x,a)=
e^{-2\sqrt{x}} x^{\frac{1}{4}}\left(1+\sum_{n=1}^{+\infty}
\frac{(4p^2-1) \ldots (4p^2-(2n-1)^2)}{n! 16^n
  x^{\frac{n}{2}}}\right)\\
 \mbox{ in }   \Sigma_0 = \{ \mid \arg(x) \mid < 3\pi \}.
\end{array}
\end{equation}
Also, by lemma \ref{lemme1} and theorem \ref{SibuyaCoeff}, we have a 
fundamental system of solutions $(\Phi_0,\Phi_{1})$ of
$(\mathfrak{E}_1)$, where $\Phi_{1}$ is characterized by the following
asymptotic expansion at infinity~:
\begin{equation}\label{Phi1pourm1}
\begin{array}{c}
 \displaystyle T\Phi_{1}(x,a)=e^{2\sqrt{x}} \omega^{\frac{1}{4}}x^{\frac{1}{4}}\left(1+\sum_{n=1}^{+\infty} (-1)^n
\frac{(4p^2-1) \ldots (4p^2-(2n-1)^2)}{n! 16^n x^{\frac{n}{2}}}\right) \\
 \mbox{ in }   \Sigma_1 = \{ \mid \arg(x) +2\pi\mid < 3\pi \}
\end{array}
\end{equation}
where $\omega = e^{2i\pi}$.
As we shall see, these functions $\Phi_0$ and $\Phi_{1}$ 
can be expressed with the
MacDonald functions $K^{(1)}_p$ and
$K^{(2)}_p$.

The MacDonald functions $K^{(1)}_p (t)$ and $K^{(2)}_p (t)$ 
are analytic functions in the variable $t$ for $t$ not equal to zero;
they are linearly
independent ($W(K^{(1)}_p,K^{(2)}_p)=\frac{\pi}{t}$) and
solutions of the modified Bessel equation (\ref{modbessel}).
These functions are derived from the Hankel functions by the
following relationships~:
\begin{equation}\label{withHankel}
\left\lbrace \begin{array}{l}
\displaystyle K^{(1)}_p(t)=\frac{1}{2}i\pi e^{\frac{pi\pi}{2}} H^{(1)}_p(it)\\
\\
\displaystyle K^{(2)}_p(t)=\frac{1}{2}\pi e^{-\frac{pi\pi}{2}} H^{(2)}_p(it).
\end{array} \right. 
\end{equation}
Furthermore, they admit respectively an asymptotic expansion $TK^{(1)}_p$ and
$TK^{(2)}_p$ when $t$ tends to infinity of the form  (see \cite{Olv74})~:
\begin{equation}\label{MacDo1}
\begin{array}{c}
 \displaystyle TK^{(1)}_p(t)=\left(\frac{\pi}{2t}\right)^{\frac{1}{2}}
 e^{-t} \left(1+\sum_{n=1}^{+\infty} \frac{(4p^2-1)
  \cdots (4p^2-(2n-1)^2)}{n! 8^n t^n}\right)\\
\\
 \mbox{ in }   \widetilde{\Sigma}_0 = \{ \mid \arg(t) \mid < \frac{3\pi}{2} \}
\end{array}
\end{equation}
and
\begin{equation}\label{MacDo2}
\begin{array}{c}
\displaystyle TK^{(2)}_p(t)=\left(\frac{\pi}{2t}\right)^{\frac{1}{2}}
e^{t} \left(1+\sum_{n=1}^{+\infty} (-1)^n \frac{(4p^2-1)
  \cdots (4p^2-(2n-1)^2)}{n! 8^n t^n}\right)\\
\\
 \mbox{ in }   
\widetilde{\Sigma}_1 = \{ \mid \arg(t) +\pi\mid < \frac{3\pi}{2} \}.
\end{array}
\end{equation}
Using (\ref{changeofvar}) and (\ref{modbessel}), we deduce by uniqueness
of $\Phi_0$ (resp. $\Phi_1$), comparing (\ref{MacDo1})
(resp. (\ref{MacDo2})) with (\ref{Phi0pourm1})
(resp. (\ref{Phi1pourm1})), 
 that
\begin{equation}\label{Phi0McDo}
\Phi_0(x,a) = \frac{2}{\sqrt{\pi}}\sqrt{x}K^{(1)}_p(2\sqrt{x})
\end{equation}
and 
\begin{equation}\label{Phi1McDo}
\Phi_1(x,a) = \frac{2}{\sqrt{\pi}}\omega^{\frac{1}{4}} 
\sqrt{x}K^{(2)}_p(2\sqrt{x}).
\end{equation}
Recalling the connection formulas (see \cite{Olv74}),
$$ \left\lbrace \begin{array}{l}
\displaystyle J_p(t)=\frac{1}{2}(H^{(1)}_p(t)+H^{(2)}_p(t))\\
\\
\displaystyle J_{-p}(t)=\frac{1}{2}(e^{i\pi p}H^{(1)}_p(t)+
e^{-i\pi p}H^{(2)}_p(t)),
\end{array} \right. $$
we deduce with remark \ref{rem109} and (\ref{withHankel}):
\begin{equation}\label{conn1}
 \left\lbrace \begin{array}{l}
\displaystyle I_p(t)=
\frac{e^{-i\pi p}}{i\pi}K^{(1)}_p(t)+\frac{1}{\pi}K^{(2)}_p(t)\\
\\
\displaystyle I_{-p}(t)=
\frac{e^{i\pi p}}{i\pi}K^{(1)}_p(t)+\frac{1}{\pi}K^{(2)}_p(t).
\end{array} \right. 
\end{equation}
Putting (\ref{attheorigin}), (\ref{Phi0McDo}), (\ref{Phi1McDo}) and (\ref{conn1}) together, we obtain
\begin{equation}\label{cestpas}
 \left(
\begin{array}{c}
f_1\\
\\
f_2
\end{array} 
\right)(x,p) = 
 \left(
\begin{array}{cc}
\displaystyle -i\frac{\Gamma (p+1)}{2\sqrt{\pi}} e^{-i\pi p} &
\displaystyle -  i\frac{\Gamma (p+1)}{2\sqrt{\pi}} \\
\\
 \displaystyle -i\frac{\Gamma (-p+1)}{2\sqrt{\pi}} e^{i\pi p} &
\displaystyle -  i\frac{\Gamma (-p+1)}{2\sqrt{\pi}}
\end{array} 
\right)
 \left(
\begin{array}{c}
\Phi_0\\
\\
\Phi_1
\end{array} 
\right)(x,a)
\end{equation} 
where the matrix on the right-hand side of this equality is the inverse of
the $0\infty$ connection matrix $M_1$ (cf. formula (\ref{lesMk})). By
proposition \ref{easyprop} we deduce:
$$M_0 (p) =  \left(
\begin{array}{cc}
\displaystyle -e^{-i\pi p}\frac{\Gamma (-p)}{\sqrt{\pi}}  &
\displaystyle - e^{i\pi p} \frac{\Gamma (p)}{\sqrt{\pi}} \\
\\
 \displaystyle \frac{\Gamma (-p)}{\sqrt{\pi}}  &
\displaystyle   \frac{\Gamma (p)}{\sqrt{\pi}}
\end{array} 
\right).$$

\begin{rem}
This result is consistent with proposition \ref{casrestreint} and remark
\ref{996}. 
\end{rem}


\protect \boldmath

\subsection{Example 2~: a normal form of Whittaker's equation}\label{cacestm2}
  
\protect \unboldmath

We now focus on the equation
$$x^2\Phi''=(x^2+a_1 x+a_2)\Phi. \quad \leqno{(\mathfrak{E}_2)}$$
This equation reduces to the Whittaker equation.
Indeed, the transformation 
$$  \left \lbrace \begin{array}{l}
x=\frac{t}{2} \\
\phi(t)=\Phi(x)
\end{array} \right. $$
converts equation $(\mathfrak{E}_2)$ into~:
$$\phi''(t)=(\frac{1}{4}+\frac{a_1}{2t}+\frac{a_2}{t^2})\phi(t)$$
which is the Whittaker equation of parameters  
$$k=-\frac{a_1}{2} \hspace{5mm} \mbox{and} \hspace{5mm}
n= \frac{p}{2}= (\frac{1}{4}+a_2)^{\frac{1}{2}}.$$
In what follows, we shall make a heavy use of the known properties of 
the special functions associated with the Whittaker equation, see for instance 
\cite{Olv74}.

\subsubsection{Study near the origin}

We assume here that $p= (1+4a_2)^{\frac{1}{2}}
  \notin \mathbb{Z}$, i.e $2n \notin \mathbb{Z}$. The fundamental
  system of solutions $(f_1, f_2)$ of theorem \ref{Fuchs0} can be
  written as follows:
\begin{equation}\label{lesfpourm2}
 \left\lbrace \begin{array}{c}
\displaystyle f_1(x,a_1,p) = 2^{-n-\frac{1}{2}}M_{k,n}(2x) \quad \text {where} \quad\\
\displaystyle  M_{k,n}(t)= e^{-\frac{t}{2}}
t^{n+\frac{1}{2}}M(n-k+\frac{1}{2},2n+1,t),
\quad 
M(\alpha,c,t)=\sum_{s=0}^{+\infty} \frac{(\alpha)_s}{(c)_s} 
\frac{t^s}{s!} \\
\\
\displaystyle f_2(x,a_1,p) = 2^{n-\frac{1}{2}}N_{k,n}(2x) 
\quad \text {where} \quad\\
\displaystyle  N_{k,n}(t)= e^{-\frac{t}{2}}
t^{n+\frac{1}{2}}N(n-k+\frac{1}{2},2n+1,t),
 \quad N(\alpha,c,t)=t^{1-c}M(1+\alpha-c,2-c,t). 
 \end{array} \right. 
\end{equation}
with the Pochhammer's notation:
$$ \left\lbrace \begin{array}{l}
(\alpha)_0=1 \\
(\alpha)_s=\alpha(\alpha+1)\ldots(\alpha+s-1).
 \end{array} \right. $$ 

\begin{rem}
The function $M_{k,n}$ is called a Whittaker function and
$M(\alpha,c,t)$ is called the Kummer function (which is an entire function in
$t$).
\end{rem}

\subsubsection{Study at infinity}

In  theorem \ref{Sibuya},  the
solution $\Phi_0$  of $(\mathfrak{E}_2)$ $\Phi_0$ can be characterized
by its asymptotics. Since
$$r (\a) = -\frac{a_1}{2}=k \hspace{5mm} \mbox{and} \hspace{5mm}
\omega=e^{i\pi},$$
we have 
  $T\Phi_0(x,\a)=x^k e^{-x} \phi_0(x,\a)$ with $\phi_0(x, \a) \in
\mathbb{C}[a_1,a_2][[x^{-1}]]$ with constant term $1$, in the sector 
$-\frac{3\pi}{2}< arg(x) < \frac{3\pi}{2}$.\\
In the same way, $\Phi_1$ is characterized by its asymptotics,
$T\Phi_1(x,\a) =\omega^{-k} x^{-k} e^x \phi_1(x,\a)$ where $\phi_1(x,\a) \in
\mathbb{C}[a_1,a_2][[x^{-1}]]$ with constant term $1$,
 in the sector $-\frac{3\pi}{2}< arg(x)+\pi < \frac{3\pi}{2}$.\\
In fact, these two functions $\Phi_0$ and $\Phi_1$ can be expressed in terms
of the functions $U$ and $V$ of the confluent hypergeometric equation: 
$$ t\frac{d^2 f}{d t^2}+(c-t)\frac{d f}{d t}-\alpha f=0.$$

\begin{prop}
$$ \left\lbrace \begin{array}{l} 
\displaystyle \Phi_0(x,\a)=2^{-k}W_{k,n}(2x) \quad \text{where} \quad
W_{k,n}(t)=e^{-\frac{t}{2}}t^{n+\frac{1}{2}}U(n-k+\frac{1}{2},2n+1,t)\\
\displaystyle  \text{and} \quad U(\alpha,c,t)\sim t^{-\alpha}\sum_{s=0}^{+\infty}
\frac{(-1)^s(1+\alpha-c)_s}{s!t^s} \text{ in the sector } \mid
arg(t) \mid < \frac{3\pi}{2}\\ 
\\
\displaystyle \Phi_{1}(x,\a)=i2^ke^{n\pi i}V_{k,n}(2x)\quad \text{where} \quad
V_{k,n}(t)=e^{-\frac{t}{2}}t^{n+\frac{1}{2}}V(n-k+\frac{1}{2},2n+1,t)\\
\displaystyle \text{and} \quad V(\alpha,c,t)\sim e^t 
(e^{i\pi}t)^{\alpha-c}\sum_{s=0}^{+\infty}
\frac{(c-\alpha)_s(1-\alpha)_s}{s!t^s} \text{ in the sector }  \mid
arg(t) + \pi \mid < \frac{3\pi}{2}\\ 
\end{array} \right. $$
\end{prop}

\begin{rem}
The function $W_{k,n}$ is also called a Whittaker function, see \cite{Olv74}.
\end{rem}

\subsubsection{Connection formulas}  
We recall the following connection formula (see \cite{Olv74}):
$$ \displaystyle M(\alpha,c,t)= \Gamma(c)\Big(\frac{e^{-\alpha\pi
    i}}{\Gamma(c-\alpha)}U(\alpha,c,t)+\frac{e^{(c-\alpha)\pi
    i}}{\Gamma(\alpha)}V(\alpha,c,t)\Big)$$
Therefore, 
$$ \displaystyle M_{k,n}(2x)=\frac{-i e^{(k-n)\pi
    i}\Gamma(2n+1)}{\Gamma(n+k+\frac{1}{2})}W_{k,n}(2x)+\frac{i
  e^{(k+n)\pi i}\Gamma(2n+1)}{\Gamma(n-k+\frac{1}{2})}V_{k,n}(2x)$$
which means that
\begin{equation}\label{pourf1}
 \displaystyle f_1(x,a_1, p)=-ie^{i \pi (k-n)}\frac{ 2^{k-n}
   \Gamma(2n+1)}{\sqrt{2}\Gamma(n+k+\frac{1}{2})}\Phi_0(x,\a)+
e^{i \pi k}\frac{2^{-k-n}
    \Gamma(2n+1)}{\sqrt{2}\Gamma(n-k+\frac{1}{2})}\Phi_1(x,\a).
\end{equation}

Furthermore, thanks to the connection formula (see \cite{Olv74}) 
$$ \displaystyle
U(\alpha,c,t)=\frac{\Gamma(1-c)}{\Gamma(1+\alpha-c)}M(\alpha,c,t)-\frac{\Gamma(c)\Gamma(1-c)}{\Gamma(\alpha)\Gamma(2-c)}N(\alpha,c,t),$$
we deduce that :
$$ \displaystyle
N(\alpha,c,t)=\Big(\frac{\Gamma(\alpha)\Gamma(2-c)e^{-\alpha\pi
    i}}{\Gamma(1+\alpha-c)\Gamma(c-\alpha)}-\frac{\Gamma(\alpha)\Gamma(2-c)}{\Gamma(c)\Gamma(1-c)}\Big)U(\alpha,c,t)+\frac{e^{(c-\alpha)\pi i}\Gamma(2-c)}{\Gamma(1+\alpha-c)}V(\alpha,c,t)$$
so that :
$$ \displaystyle
N_{k,n}(2x)=\Big(\frac{\Gamma(n-k+\frac {1}{2})\Gamma(1-2n)e^{-(n-k+\frac{1}{2})\pi
    i}}{\Gamma(-n-k+\frac{1}{2})\Gamma(n+k+\frac{1}{2})}-\frac{\Gamma(n-k+\frac{1}{2})\Gamma(1-2n)}{\Gamma(2n+1)\Gamma(-2n)}\Big)W_{k,n}(2x)$$

$$\displaystyle +
\Big(\frac{e^{(n+k+\frac{1}{2})\pi i}\Gamma(1-2n)}{\Gamma(-n-k+\frac{1}{2})}\Big)V_{k,n}(2x) $$
 i.e :
$$
\begin{array}{c}
 \displaystyle f_2(x,a_1,p)=\\
\\
\displaystyle \Big(\frac{-i2^{k+n} e^{(k-n)\pi
    i}\Gamma(1-2n)\Gamma(n-k+\frac{1}{2})}{\sqrt{2}\Gamma(n+k+\frac{1}{2})\Gamma(-n-k+\frac{1}{2})}-\frac{2^{k+n}\Gamma(1-2n)\Gamma(n-k+\frac{1}{2})}{\sqrt{2}\Gamma(2n+1)\Gamma(-2n)}\Big)\Phi_0(x,\a)\\
\\
 +\displaystyle \Big(\frac{2^{n-k}e^{k\pi
    i}\Gamma(1-2n)}{\sqrt{2}\Gamma(-n-k+\frac{1}{2})}\Big)\Phi_1(x,\a).
\end{array}
$$
which reads also:
\begin{equation}\label{pourf2}
\begin{array}{c}
 \displaystyle f_2(x,a_1,p)= 
-i e^{i \pi (n+k)} \frac{2^{k+n}\Gamma(1-2n)}{\sqrt{2}\Gamma(-n+k +
  \frac{1}{2} )}\Phi_0(x,\a)
 +\displaystyle e^{ i\pi k}\frac{2^{n-k}\Gamma(1-2n)}{\sqrt{2}\Gamma(-n-k+\frac{1}{2})}\Phi_1(x,\a).
\end{array}
\end{equation}

Formulas (\ref{pourf1}) and (\ref{pourf2}) yield the inverse of the $0\infty$
connection matrix $M_1$. Going back to the variables $p$ and $a_1$,
one  gets:

\begin{equation}
 M_1^{-1}(a_1,p) = 
\left(
\begin{array}{cc} 
\displaystyle -ie^{i \pi (-\frac{a_1}{2}-\frac{p}{2})}\frac{ 2^{-\frac{a_1}{2}-\frac{p}{2}}
   \Gamma(p+1)}{\sqrt{2}\Gamma(\frac{p}{2}-\frac{a_1}{2}+\frac{1}{2})} &
\displaystyle  e^{-i \pi \frac{a_1}{2}}\frac{2^{\frac{a_1}{2}-\frac{p}{2}}
    \Gamma(p+1)}{\sqrt{2}\Gamma(\frac{p}{2}+\frac{a_1}{2}+\frac{1}{2})}  \\
       \\ & \\                                        
 \displaystyle
-i e^{i \pi (-\frac{a_1}{2}+\frac{p}{2})} \frac{2^{-\frac{a_1}{2}+\frac{p}{2}}\Gamma(1-p)}{\sqrt{2}\Gamma(-\frac{p}{2}-\frac{a_1}{2} +
  \frac{1}{2} )} &
 \displaystyle e^{ -i\pi\frac{a_1}{2}}\frac{2^{\frac{a_1}{2} +\frac{p}{2}}\Gamma(1-p)}{\sqrt{2}\Gamma(-\frac{p}{2}+\frac{a_1}{2}+\frac{1}{2})}
 \end{array} \right). 
\end{equation}

Using 
proposition \ref{easyprop} we deduce that

\begin{equation}\label{finalpourMo}
 M_0(a_1,p) = 
\left(
\begin{array}{cc} 
 \displaystyle 
-i e^{-i \pi \frac{p}{2}} \frac{2^{-\frac{a_1}{2}+\frac{p}{2}+1}\Gamma(-p)}{ \sqrt{2}\Gamma(-\frac{p}{2}-\frac{a_1}{2}+\frac{1}{2})} &
\displaystyle -ie^{i \pi \frac{p}{2}}\frac{2^{-\frac{a_1}{2}-\frac{p}{2}+1}
    \Gamma(p)}{\sqrt{2}\Gamma(\frac{p}{2}-\frac{a_1}{2}+\frac{1}{2})}  \\
 \\ & \\     
\displaystyle 
  \frac{2^{\frac{a_1}{2}+\frac{p}{2}+1}\Gamma(-p)}{\sqrt{2}\Gamma(-\frac{p}{2}+\frac{a_1}{2} + \frac{1}{2} )} &
\displaystyle  \frac{ 2^{\frac{a_1}{2}-\frac{p}{2}+1}
   \Gamma(p)}{\sqrt{2}\Gamma(\frac{p}{2}+\frac{a_1}{2}+\frac{1}{2})}
 \end{array} \right). 
\end{equation}

\begin{rem}
This result is consistent with theorem \ref{0infty0}. Also, with the
notations of theorem \ref{0infty0}, we have found
$$\widetilde{L}_0 (a_1, p) = \displaystyle -ie^{i \pi \frac{p}{2}}\frac{2^{-\frac{a_1}{2}-\frac{p}{2}+1}
    \Gamma(p)}{\sqrt{2}\Gamma(\frac{p}{2}-\frac{a_1}{2}+\frac{1}{2})}.
$$
In particular, when $a_1=0$  we get, using the Legendre
    duplication formula for the Gamma function:
$$\widetilde{L}_0 (0,p) = \displaystyle 2^\frac{p}{2}e^{i \pi \frac{p}{2}}\frac{(-i)
  }{\sqrt{2 \pi}}\Gamma(\frac{p}{2}), $$
a result which agrees also with proposition \ref{casrestreint} and remark
\ref{996}.
\end{rem}

By formula (\ref{SetM}) of theorem \ref{FUNCT}, we have 
$\displaystyle \mathfrak{S}_0 (\a) =  M_0 (a_1,p)
  M_{1}^{-1}(a_1,p)  
$ 
and the result extends to $2n \in \mathbb{Z}$ by analytic
continuation, since $\mathfrak{S}_0$ is entire.  We eventually get:

\begin{prop}
We assume $m=2$. Then, for all $\a=(a_1, a_2) \in \mathbb{C}^2$, the
Stokes-Sibuya connection matrix $ \mathfrak{S}_0$ is given by
 $$ \mathfrak{S}_0 (\a)= 
\left(
\begin{array}{cc} 
\displaystyle -2i e ^{ -i\pi\frac{a_1}{2}} 2^{-a_1} 
\frac{\Gamma(\frac{p}{2}+\frac{a_1}{2}+\frac{1}{2})}{\Gamma(\frac{p}{2}-\frac{a_1}{2}+\frac{1}{2})}
\cos\left((\frac{p}{2}+\frac{a_1}{2})\pi\right)&
      e ^{-i\pi a_1} \\
 & \\
     \displaystyle 1 & \displaystyle 0 
\end{array} 
\right) $$
where 
$$p=  (1+4a_2)^{\frac{1}{2}}.$$
Moreover, when $p \notin \mathbb{Z}$, the $0\infty$ connection matrix
$M_0$ is given by formula (\ref{finalpourMo}).
\end{prop}

\end{document}